\documentclass[final]{elsarticle}

\journal{}
\usepackage[USenglish]{babel}
\usepackage{lineno}
\usepackage{pifont}
\usepackage{hyperref}
\usepackage{natbib}
\usepackage{geometry}
\usepackage{graphicx}
\usepackage{amsmath}
\usepackage{amsfonts}
\usepackage{amssymb}
\usepackage{mathrsfs}
\usepackage{multirow}
\usepackage{cases}
\usepackage{subfig}
\usepackage{booktabs}
\usepackage{floatrow}
\usepackage{xcolor}
\usepackage{algorithm} %for algorithms

\let\today\relax
\makeatletter
\def\ps@pprintTitle{%
    \let\@oddhead\@empty
    \let\@evenhead\@empty
    \def\@oddfoot{\footnotesize\itshape
         {Submitted preprint} \hfill\today}%
    \let\@evenfoot\@oddfoot
    }
\makeatother

\usepackage{algpseudocode} %for \begin{algorithm}
\newfloatcommand{capbtabbox}{table}[][\FBwidth]

\begin{document}

\begin{frontmatter}

\title{Towards solving large--scale topology optimization problems with buckling constraints at the cost of linear analyses}

% Authors
\author[dtu]{Federico Ferrari\fnref{ft1}}
 \ead{feferr@mek.dtu.dk}

\author[dtu]{Ole Sigmund\fnref{ft2}\corref{cor1}}
 \ead{sigmund@mek.dtu.dk}

% Positions
\fntext[ft1]{Post Doctoral Fellow}
\fntext[ft2]{Full Professor of Solid Mechanics}

% Corresponding Authors
\cortext[cor1]{Corresponding author}

% Address
\address[dtu]{Department of Mechanical Engineering, Technical University of Denmark, Kongens Lyngby, (DK)}

\vspace{12pt}

\textbf{The official version of the paper can be downloaded at \url{https://doi.org/10.1016/j.cma.2020.112911}}

\begin{abstract}
This work presents a multilevel approach to large--scale topology optimization accounting for linearized buckling criteria. The method relies on the use of preconditioned iterative solvers for all the systems involved in the linear buckling and sensitivity analyses and on the approximation of buckling modes from a coarse discretization. The strategy shows three main benefits: first, the computational cost for the eigenvalue analyses is drastically cut. Second, artifacts due to local stress concentrations are alleviated when computing modes on the coarse scale. Third, the ability to select a reduced set of important global modes and filter out less important local ones. As a result, designs with improved buckling resistance can be generated with a computational cost little more than that of a corresponding compliance minimization problem solved for multiple loading cases. Examples of 2D and 3D structures discretized by up to some millions of degrees of freedom are solved in Matlab to show the effectiveness of the proposed method. Finally, a post--processing procedure is suggested in order to reinforce the optimized design against local buckling.
\end{abstract}


\end{frontmatter}

\section{Introduction}
 \label{Sec:Intro}

% Introductory hook: what we want to do and which problems we have to overcome
This paper describes the benefits of using a multilevel approximation method for computing buckling modes in the context of large--scale Topology Optimization (TO). In particular, we show that besides remarkable computational savings, some issues arising when buckling criteria are accounted for in TO may be alleviated.

% General discussion about state-of-art and difficulties encountered by other researchers
Topology Optimization is rapidly spreading to engineering practice as a promising, powerful tool for the conceptual design of whole components, or for highly detailed microstructures and architected material \cite{clausen-etal_16a, thomsen-etal_18a}. Therefore, large--scale applications have become a hot topic \cite{aage-etal_17a, chin-kennedy_16a} and this may be credited to the leverage effect of rapidly increasing computational capabilities and the emergence of new manufacturing techniques.

However, there are still important issues to be overcome when considering geometric nonlinearities, as required from a stability analysis. In such situations, the computational effort due to analysis substantially increases and the optimization problem becomes much harder, showing poor conditioning and many local, possibly non--physical, solutions. Therefore, even the simplest approach to stability optimization, based on linearized buckling, is far from being a customary and free--from--issues practice in TO. This puts a severe limitation on the dissemination of TO to engineering practice, as realistic structures must meet stability requirements, whereas they could even be weakened if optimized without accounting for this phenomenon \cite{thompson-hunt_74a}. Hence, there is still a substantial gap between the scale of compliance--based TO problems and those accounting for buckling, yet some works have appeared on this topic in recent years. Dunning et al. \cite{dunning-etal_16a} optimized a 3D structure for minimum mass with 144,000 design variables, by using a robust eigenvalue solver capable of dealing with clustered eigenvalues \cite{ovtchinnikov_08a}. However, their approach still needs (at least) one factorization of the full system matrix, which is undesirable for large scale problems. While studying mass and compliance minimization for the Common Research Model wing, Chin and Kennedy \cite{chin-kennedy_16a} considered buckling constraints on 2D panels. They reported the huge sensitivity of the obtained designs to the number of constrained buckling modes. Bian and Feng \cite{bian-feng_17a} proposed an assembly--free iterative solver for the eigenvalue problem coupled with a voxelization--based discretization, leading to low memory requirements and parallelization capabilities, and solved problems in the order of $2$ to $5\cdot 10^{5}$ DOFs.

% Motivate why the presence of many local buckling modes is an issue and why we want to keep these modes out from the optimization process
The following issues are systematically pointed out from these and other works \cite{bruyneel-etal_08a, ferrari-sigmund_19a}:
\begin{enumerate}
 \item High computational cost due to repeated solution of large eigenvalue problems;
 \item High sensitivity of the results to the set of buckling modes considered in the formulation;
 \item Activation and clustering of many buckling modes as the optimization progresses.
\end{enumerate}

The issues above are partly interlaced, as the need to account for a large set of buckling modes increases the computational cost of many standard eigensolvers \cite{knyazev_01a, zhou-saad_08a}. This is even worsened by the appearance of artificial and/or not physically meaningful deformations associated with low eigenvalues \cite{neves-etal_95a, book:bendsoe-sigmund_2004}. Buckling modes appearing in low--density regions are a classical issue for TO formulations involving eigenvalues and several approaches are available to identify and filter them \cite{ma-etal_93a, pedersen_00a, gao-ma_15a}. However, even some localized buckling modes appearing in solid regions may be physically meaningless or undesirable to take into account within the optimization process, at least for two reasons. First, some highly localized instabilities are due to stress concentrations and singularities linked to geometrical irregularities, and these phenomena are eventually worsened by discretization effects. Second, physical but still local modes, such as the failure of single bars, may be inconvenient to consider, especially as the design space becomes larger and the structural layout more intricate. In this case, as hierarchical structures form (c.f. \cite{thomsen-etal_18a, ferrari-sigmund_19a}), and many thin bars appear, the number of active modes grows rapidly, making it unfeasible to include them all in the optimization process.

% How can we achieve this goal with our approach? Describe the main innovative features of the proposed approach;
We aim at developing a methodology which efficiently takes into account the most global buckling modes only, driving the optimization towards a structure having global stability. Subsequently, tiny and slender features which may undergo local buckling may be fixed in a post--processing phase by means of a local reinforcement. A similar effect is often used in truss TO, where buckling of individual bars is ignored or handled by separate constraints on element compressive stresses \cite{khot-etal_76a, achtziger_99a, achtziger_99b}.

% Connect what is proposed in the present work with what was done for natural vibrations. Give also a feel of what are the specific novel features of the proposed method when applied to buckling problems;
Our goal can be accomplished with inspiration from another multilevel concept outlined in \cite{ferrari-etal_18a}, in the context of dynamic eigenvalue problems. The strategy there was to compute an approximation to vibration modes from a coarse discretization, project them on the (much finer) discretization where the optimization takes place, and use them as drivers for a harmonic response (linear) problem. The method was given a physical motivation, replacing the eigenvalue problem with a frequency response one, and similarities between the two were discussed by the authors \cite{andreassen-etal_17a, ferrari-etal_18a}. A multilevel concept is exploited also here, with the main focus on cheaply computing a satisfactory approximation to some buckling modes without ever solving an eigenvalue problem on the fine discretization. However, in contrast to \cite{ferrari-etal_18a}, the approximate buckling modes and their associated load factors estimated by means of the Rayleigh quotient, are now directly used to run the optimization. We emphasize that this multilevel approach helps filtering the aforementioned unphysical artifacts and localized buckling modes originating on the fine grid, thus significantly simplifying and speeding up the optimization process.

% Outline of the paper
The outline of the paper is as follows. In \autoref{Sec:FormulationOpt} we set up the problem and describe the steps for the multilevel approximation of buckling modes and load factors. In \autoref{Sec:2Dexample-2BarsTruss} we present 2D results showing the potential of the method to produce buckling resistant designs with very low computational cost. Then, the designs are carefully discussed in \autoref{Sec:AnalysisDesigns2BS} and the ability of our approach to overcome some artificial effects is discussed. A post--processing strategy to reinforce the design against local buckling is proposed in \autoref{sSecSpuriousLocalizedModes}. A fairly large 3D example is shown in \autoref{Sec:3DcantileverMassMinimization} and conclusions are drawn in \autoref{Sec:Conclusions}, including a discussion about open issues.

\begin{algorithm}[t]
 \caption{Linearized Buckling Analysis (LBA)}
  \label{alg:EigenvalueBucklingAnalysis}
   \begin{algorithmic}[1]
        \State Select a reference load vector $\mathbf{f} \in \mathbb{R}^{n}$
        \State \textbf{LA :} Compute the equilibrium displacement $\mathbf{u} = K\left[ \mathbf{x} \right]^{-1} \mathbf{f}$
        \State Set up the stress stiffness matrix $G\left[ \mathbf{x}, \mathbf{u}\left( \mathbf{x} \right) \right] \in {\rm Sym}^{n\times n}$
        \State \textbf{EA :} Compute the pairs $\left( \lambda_{i}, \boldsymbol{\varphi}_{i} \right) \in \mathbb{R}\times \mathbb{R}^{n}$, $i = 1, \ldots, r$ by solving
        \begin{equation}
         \label{eq:eigenvalueProblem0}
         \left( K\left[ \mathbf{x} \right] + \lambda G\left[ \mathbf{x}, \mathbf{u}\left( \mathbf{x} \right) \right] \right) \boldsymbol{\varphi} = \mathbf{0} \: , \qquad \boldsymbol{\varphi} \neq \mathbf{0}
        \end{equation}
   \end{algorithmic}
\end{algorithm}

\section{Setting and methods}
 \label{Sec:FormulationOpt}

We consider a continuum body $\Omega \subset \mathbb{R}^{d}$, $d = \{2, 3\}$, and its discretization $\Omega_{1} = \cup^{m}_{e=1}\Omega^{e}$, obtained through a uniform and regular grid of $m$ elements $\Omega_{e}$. Hereafter $\Omega_{1}$ will be referred to as the \emph{fine} discretization and $n$ denotes the total number of Degrees of Freedom (DOFs).

Let $\mathbf{\hat{x}} = \{ \hat{x}_{e} \}^{m}_{e = 1}$ be the vector of design variables. We consider a three field approach to impose a length scale on the design \cite{lazarov-etal_16a}. Physical variables $x_{e}$ are given by the relaxed Heaviside projection ($\eta\in\left[ 0, 1 \right]$ and $\beta\in\left[ 1, \infty \right)$) \cite{wang-etal_11a}

\begin{equation}
 \label{eq:projectionOperator}
  x_{e}\left( \tilde{x}_{e}, \eta, \beta \right) =
  \frac{\tanh(\beta\eta) + \tanh(\beta(\tilde{x}_{e} - \eta))}
  {\tanh(\beta\eta) + \tanh(\beta(1-\eta))}
\end{equation}
where $\tilde{x}_{e}=\tilde{x}_{e}(\hat{x}_{e})$ is obtained thorough a standard density filter \cite{bourdin_01a}, with radius $r_{\rm min}$.

The global elastic stiffness matrix $K = K[\mathbf{x}]$ and the global stress stiffness matrix $G = G[\mathbf{x}, \mathbf{u}(\mathbf{x})]$, which depends on both design variables and displacements $\mathbf{u}$, are assembled from the elemental ones $K_{e}[x_{e}] = E_{\kappa}(x_{e})K_{e0}$ and $G_{e} = E_{\sigma}(x_{e})G_{e0}[\mathbf{u}_{e}(x_{e})]$, respectively. These depend on the two different interpolations of the Young modulus \cite{book:bendsoe-sigmund_2004}

\begin{equation}
 \label{eq:interpolationSIMP2-SIMP}
  \begin{aligned}
   E_{\kappa}\left( x_{e} \right) & =
    E_{0} + x^{p}_{e} \left( E_{1} - E_{0} \right) \\
   E_{\sigma}\left( x_{e} \right) & =         x^{p}_{e} E_{1}
  \end{aligned}
\end{equation}
where the contrast in coefficients is $E_{1}/E_{0} = 10^{6}$ and the Poisson ratio is fixed to $\nu = 0.3$.

The element matrices $K_{e0}$ and $G_{e0}$ are obtained using incompatible finite elements (i.e. 6--DOFs Wilson quadrilaterals in 2D and 11--DOFs hexahedra in 3D) whose description and  implementation details can be found e.g. in \cite{wilson-etal_73a, wilson-ibrahimbergovic_90a, book:pian-wu2005}. The benefits of such elements applied to buckling problems have been observed \cite{ferrari-sigmund_19a} and we remark that, for an elementwise constant material distribution, these are equivalent to some mixed elements \cite{pian-cheng_82a, pian-sumihara_84a}.

Accounting for geometrical non--linearities \cite{book:deborst2012}, buckling under the applied load $\mathbf{f}$ is described by the modes $\boldsymbol{\varphi}_{i}\in\mathbb{R}^{n}$ and the associated Buckling Load Factors (BLFs) $\lambda_{i}\in\mathbb{R}$, $i = 1\ldots n$. Thus, the fundamental BLF, which can be characterized by the Rayleigh quotient \cite{book:washizu}

\begin{equation}
 \label{eq:minBLFcharRayleigh}
  \lambda_{1}\left( \mathbf{x}, \mathbf{u} \right) = 
   \min\limits_{\mathbf{v}\in\mathbb{R}^{n}, \mathbf{v}\neq\mathbf{0}}
     R\left( \mathbf{x}, \\ \mathbf{v} \right) := - \frac{\mathbf{v}^{T}K\left[ \mathbf{x} \right]\mathbf{v}}
 {\mathbf{v}^{T}G\left[ \mathbf{x}, \mathbf{u} \right]\mathbf{v}}
\end{equation}
provides an approximate measure of the stability of the discretized system. In the following we consider the buckling modes to be normalized such that $\boldsymbol{\varphi}^{T}_{i} K[\mathbf{x}]\boldsymbol{\varphi}_{j} = \delta_{ij}$.

For a simple eigenvalue $\lambda_{i}$ the sensitivity w.r.t. each variable $x_{e}$ is expressed as \cite{rodrigues-etal_95a}

\begin{equation}
 \label{eq:SensitivityLambda}
  \frac{\partial \lambda_{i}}{\partial x_{e}} = \boldsymbol{\varphi}^{T}_{i} \left( \frac{\partial K}{\partial x_{e}} + \lambda_{i} \frac{\partial G}{\partial x_{e}} \right)
  \boldsymbol{\varphi}_{i} - \lambda_{i}\mathbf{z}^{T}_{i}\frac{\partial K}{\partial x_{e}}\mathbf{u}
\end{equation}
where $\mathbf{z}_{i}$ solves the adjoint system

\begin{equation}
 \label{eq:AdjointSystem}
 K\mathbf{z}_{i} = \boldsymbol{\varphi}^{T}_{i}\left( \nabla_{\mathbf{u}}G \right)\boldsymbol{\varphi}_{i}
\end{equation}
and the chain rule must then be applied to recover the filter and projection dependence \eqref{eq:projectionOperator}.

\subsection{Approximation of the buckling modes by the multilevel procedure}
\label{sSec:multilevelBucklingModes}

In a standard \textit{nested} TO approach the BLFs and buckling modes are computed at each optimization step through a Linearized Buckling Analysis (see \autoref{alg:EigenvalueBucklingAnalysis}), consisting of a linear analysis (LA) and an eigenvalue analysis (EA). The latter represents the main computational burden in the LBA, rapidly increasing with the number of DOFs. Moreover, a large and growing number of eigenpairs may be required in order to consider all the active buckling modes \cite{neves-etal_95a, bruyneel-etal_08a, ferrari-sigmund_19a}, further increasing the cost of each analysis.

Therefore, we propose to take advantage of the multilevel discretization used for setting up the multigrid preconditioner when performing the LA \cite{amir-etal_14a} in order also to cheaply compute an approximation to the buckling modes. Let $\Omega_{\ell}$ be the coarsest discretization, $\Omega_{j}$ an intermediate one and $I^{j}_{j + 1}$, $I^{j + 1}_{j}$ the interpolation and restriction operators between two consecutive levels \cite{book:briggs00}. The procedure is summarized in the following steps

\begin{enumerate}
 \item Solve the coarse scale eigenvalue problem
 
  \begin{equation}
   \label{eq:eigenvalueProblemCoarseScale}
    \left( K^{\ell}\left[ \mathbf{x} \right] + \lambda^{\ell} G^{\ell}\left[ \mathbf{x}, \mathbf{u}\left( \mathbf{x} \right) \right] \right)\boldsymbol{\psi}^{\ell} = \mathbf{0} \, , \qquad \boldsymbol{\psi}^{\ell} \neq \mathbf{0}
  \end{equation}
  where $K^{\ell}\left[\mathbf{x}\right]$ and $G^{\ell}\left[\mathbf{x}, \mathbf{u}\left( \mathbf{x} \right)\right]$ are obtained from the fine scale operators through Galerkin projection \cite{book:briggs00} (we omit the superscript $\ell = 1$ when referring to quantities on the fine scale);
  \item The set of $q$ lowest coarse scale modes, say $\Psi^{\ell} = \{\boldsymbol{\psi}^{\ell}_{i}\}^{q}_{i=1}$ is projected on $\Omega_{1}$ through each $\Omega_{j}$, by means of the iteration $\Psi^{j} = I^{j}_{j + 1}\Psi^{j + 1}$;
  \item Once on the finest scale $\Omega_{1}$, the following linear problem is solved computing an approximation of the fine scale modes $\tilde{\Phi} = \{ \tilde{\boldsymbol{\varphi}}_{i} \}^{q}_{i=1}$

  \begin{equation}
   \label{eq:FineScaleBVP}
    K\left[ \mathbf{x} \right] \tilde{\Phi} =
    G\left[ \mathbf{x}, \mathbf{u}\left( \mathbf{x} \right) \right] \Psi
  \end{equation}
 \item Finally, the corresponding BLFs are calculated as
 
\begin{equation}
 \label{eq:FineScaleRQ}
  \{\tilde{\lambda}_{i}\}^{q}_{i=1} = R(\tilde{\Phi}) = 
  - \frac{\tilde{\Phi}^{T}K\left[ \mathbf{x} \right]\tilde{\Phi}}
  {\tilde{\Phi}^{T}G\left[ \mathbf{x}, \mathbf{u}\left( \mathbf{x} \right) \right]\tilde{\Phi}}
\end{equation}
\end{enumerate}

The idea is now to use the pairs $(\tilde{\lambda}_{i}, \tilde{\boldsymbol{\varphi}}_{i})$, in place of the fine scale eigenpairs $( \lambda_{i}, \boldsymbol{\varphi}_{i} )$, to run the optimization, thus eliminating the need for solving any eigenvalue equation on $\Omega_{1}$. Therefore, the pairs $(\tilde{\lambda}_{i}, \tilde{\boldsymbol{\varphi}}_{i})$ are used within the sensitivity expressions \eqref{eq:SensitivityLambda} and \eqref{eq:AdjointSystem}. This is formally not consistent, as $( \tilde{\lambda}_{i}, \tilde{\boldsymbol{\varphi}}_{i} )$ is not an eigenpair and we have the following residual

\begin{equation}
 \label{eq:eigResidual}
  \mathbf{y} = (K\left[ \mathbf{x} \right] + \tilde{\lambda}_{1} G\left[ \mathbf{x}, \mathbf{u}\left( \mathbf{x} \right) \right]) \tilde{\boldsymbol{\varphi}}_{1}
\end{equation}

The formally consistent sensitivity expression, which requires the solution of an extra adjoint problem, is given in \ref{Sec:App-NumericalProcedures}. However, here we treat \eqref{eq:FineScaleRQ} as yet another approximation to the exact Rayleigh quotients and the contribution associated with \eqref{eq:eigResidual} is not accounted for.

The procedure above is rooted in the method originally described in \cite{ferrari-etal_18a}, where an analogy with Preconditioned Inverse Iteration (PInvIt) \cite{golub-ye_00a, neymeyr_01a} was pointed out. Steps 1 and 2 were recognized as a cheap way for computing a very good initial guess for the inverse iteration on the fine scale \cite{ferrari-etal_18a}. Here the focus is on more than numerical details. The pivotal advantage of the proposed method lies in the filtering of some artificial and/or highly localized buckling modes, but still preserving the quality of ``global'' ones. By computing only these latter, the optimization runs towards a structure with improved global stability without overly enlarging the set of constrained modes.

For the sake of brevity and to quickly focus on the core results, we do not discuss numerics any deeper. The interested reader may find some further details and validations in \ref{Sec:App-NumericalProcedures}.

%% GEOMETRICAL SKETCH AND OPTIMIZED DESIGNS
\begin{figure}[tb]
 \centering
  \subfloat[]{
   \includegraphics[scale = 0.275, keepaspectratio]
   {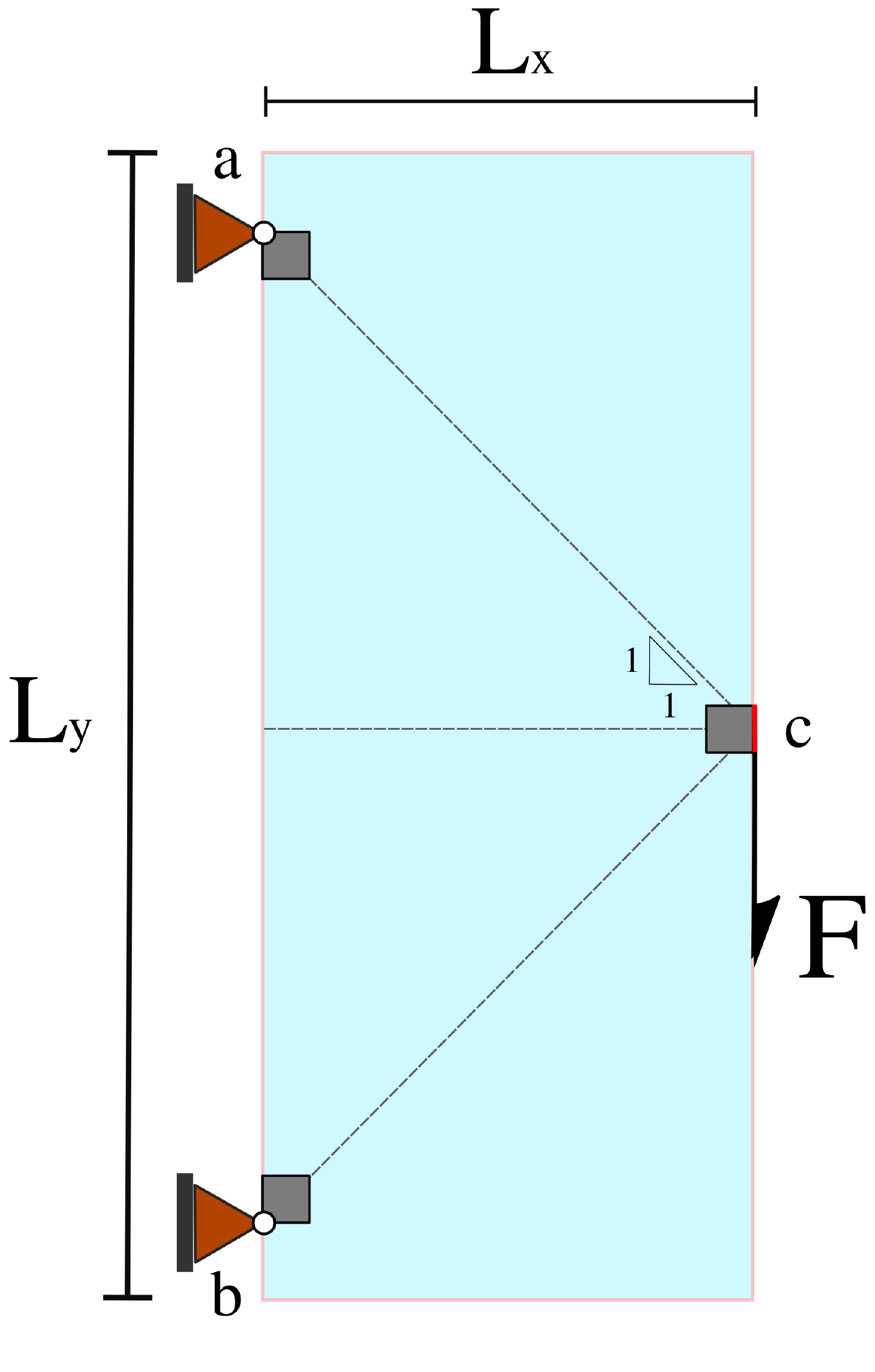}} \qquad
  \subfloat[$\Omega_{1} = 840 \times 360$]{
   \includegraphics[scale = 0.3, keepaspectratio]
   {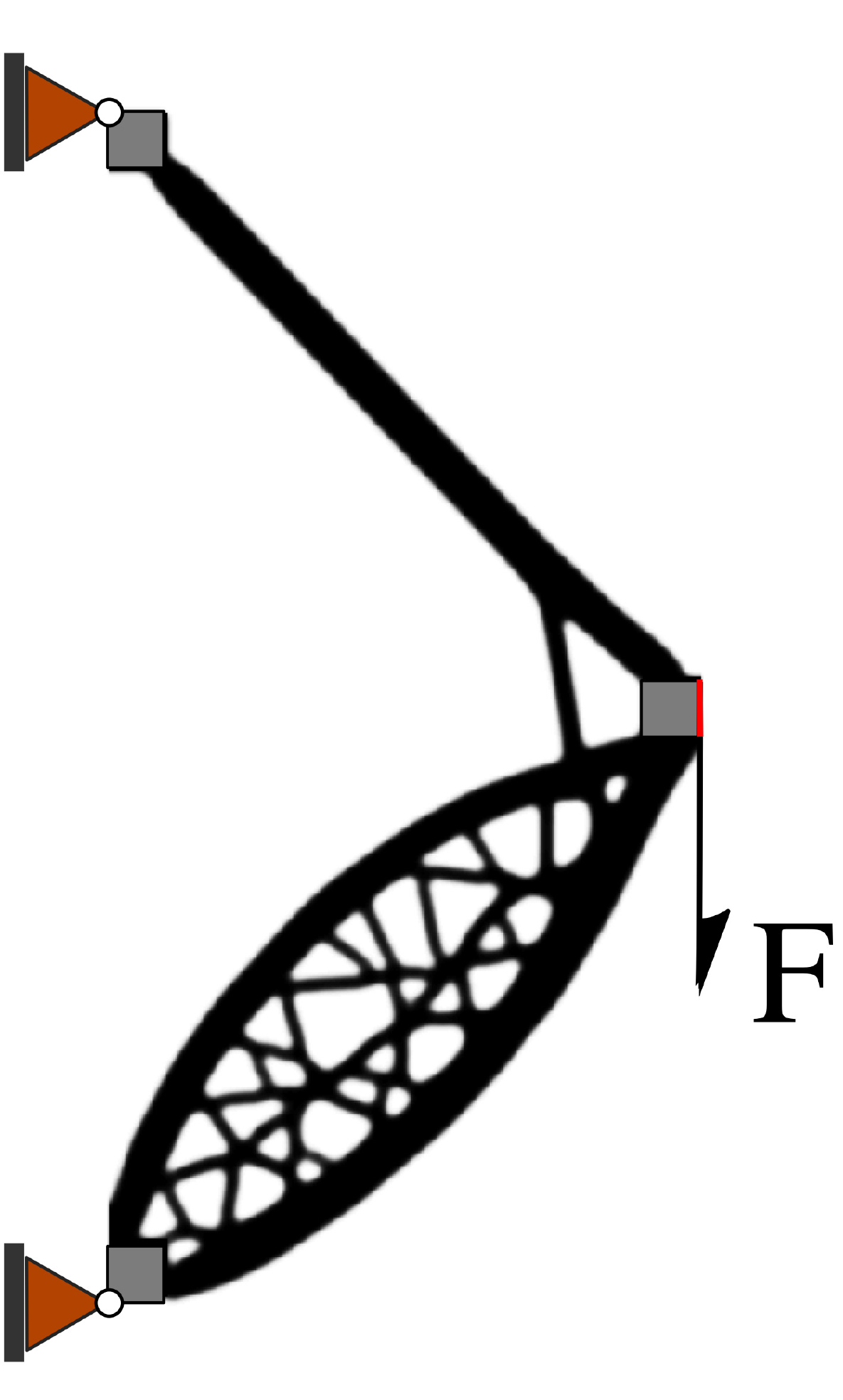}}
  \qquad
  \subfloat[$\Omega_{1} = 1680 \times 720$]{
   \includegraphics[scale = 0.3, keepaspectratio]
   {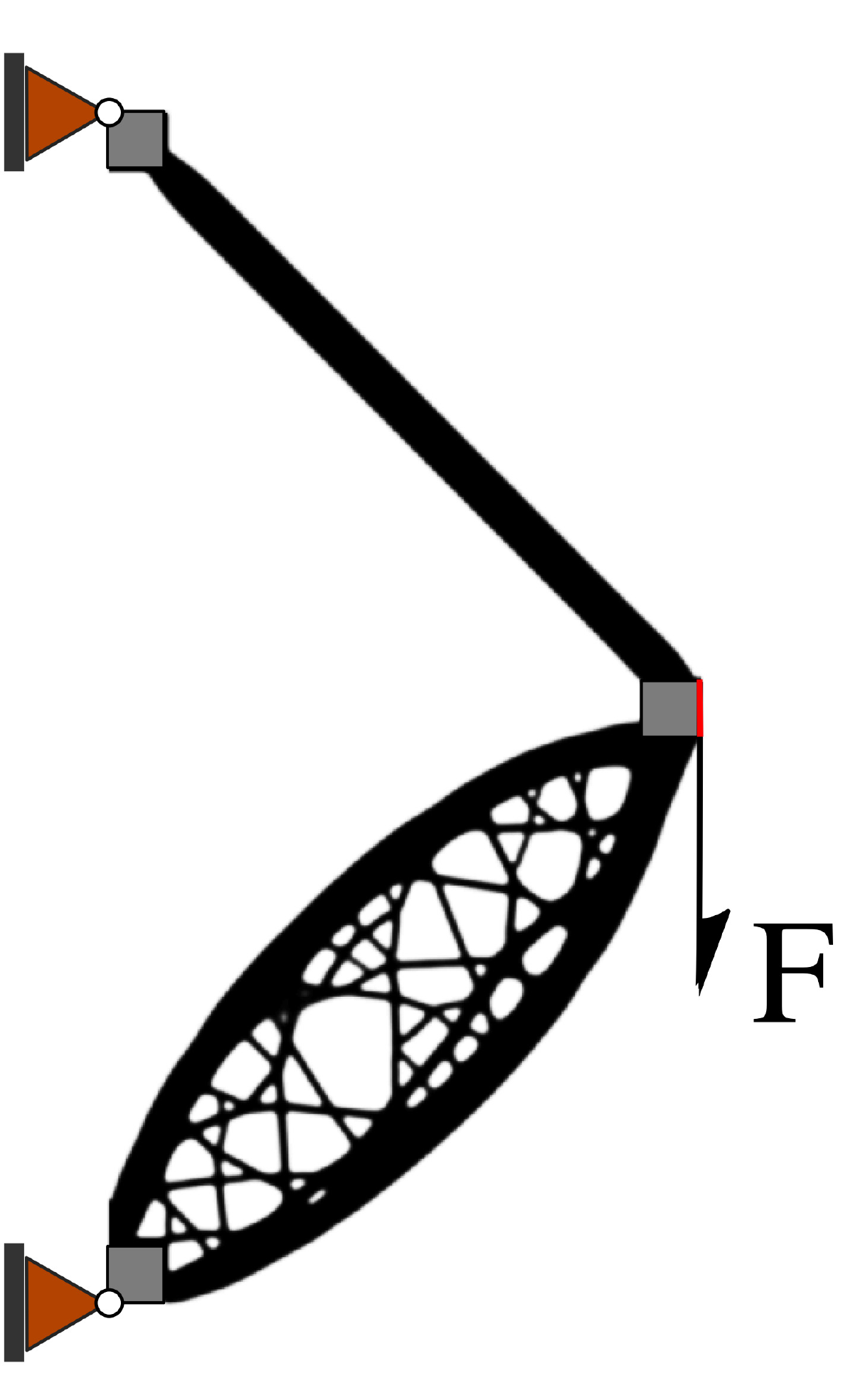}}
 \caption{Geometry of the two--bar frame example introduced in \cite{ferrari-sigmund_19a} (a) and minimum compliance designs, corresponding to $\hat{f} = 0.16$ and $\bar{\lambda} = 1.0$, on two different discretizations (b), (c)}
 \label{fig:NumericalTest-2DExample}
\end{figure}

%% OPTIMIZATION HISTORIES AND EVOLUTION OF THE BUCKLING LOAD FACTORS FOR THE MIN_CPL AND MIN_VOL PROBLEMS
\begin{figure}[tb]
 \centering
  \subfloat[]{
   \includegraphics[scale = 0.45, keepaspectratio]
   {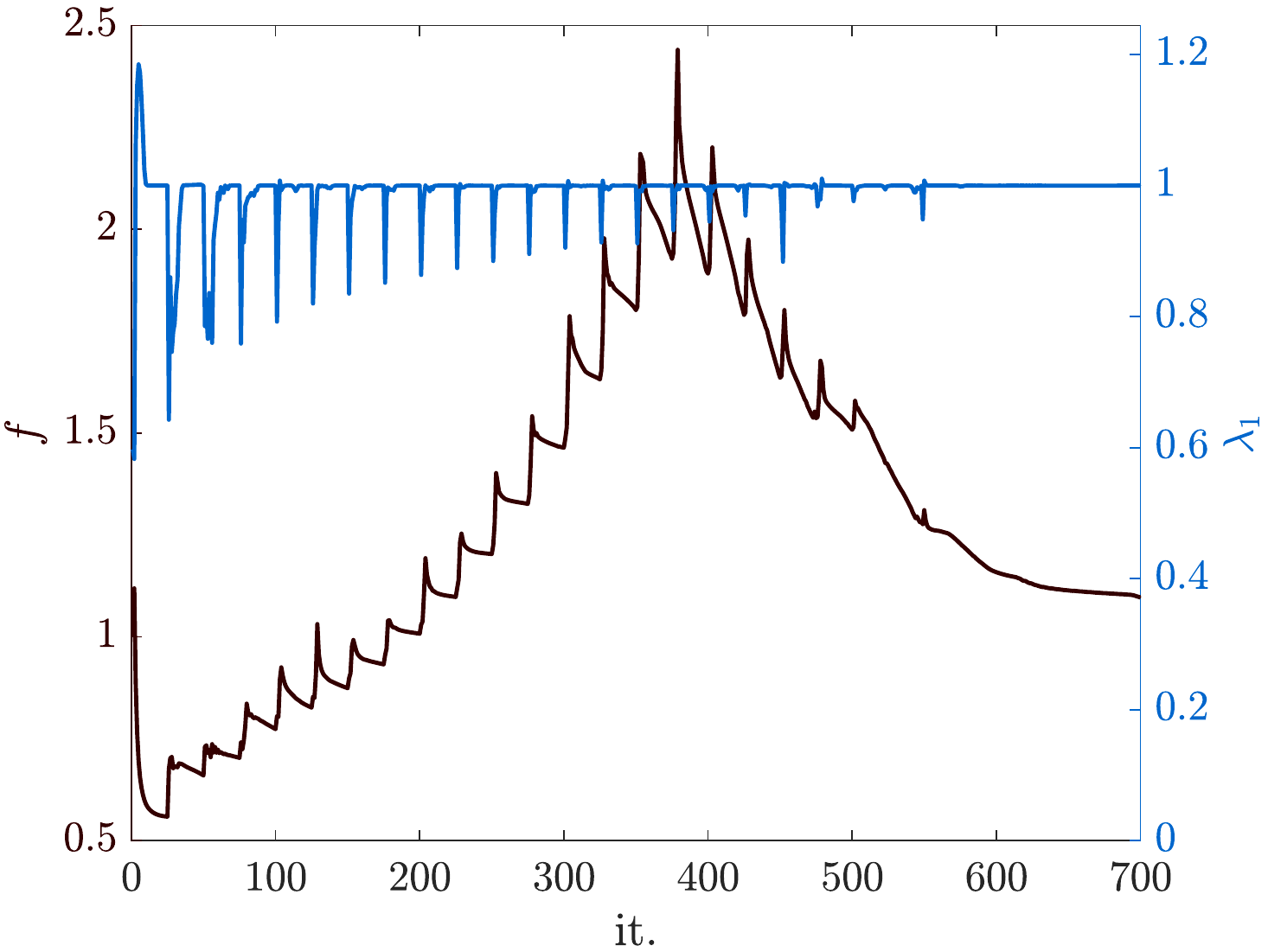}} \qquad
  \subfloat[]{
   \includegraphics[scale = 0.45, keepaspectratio]
   {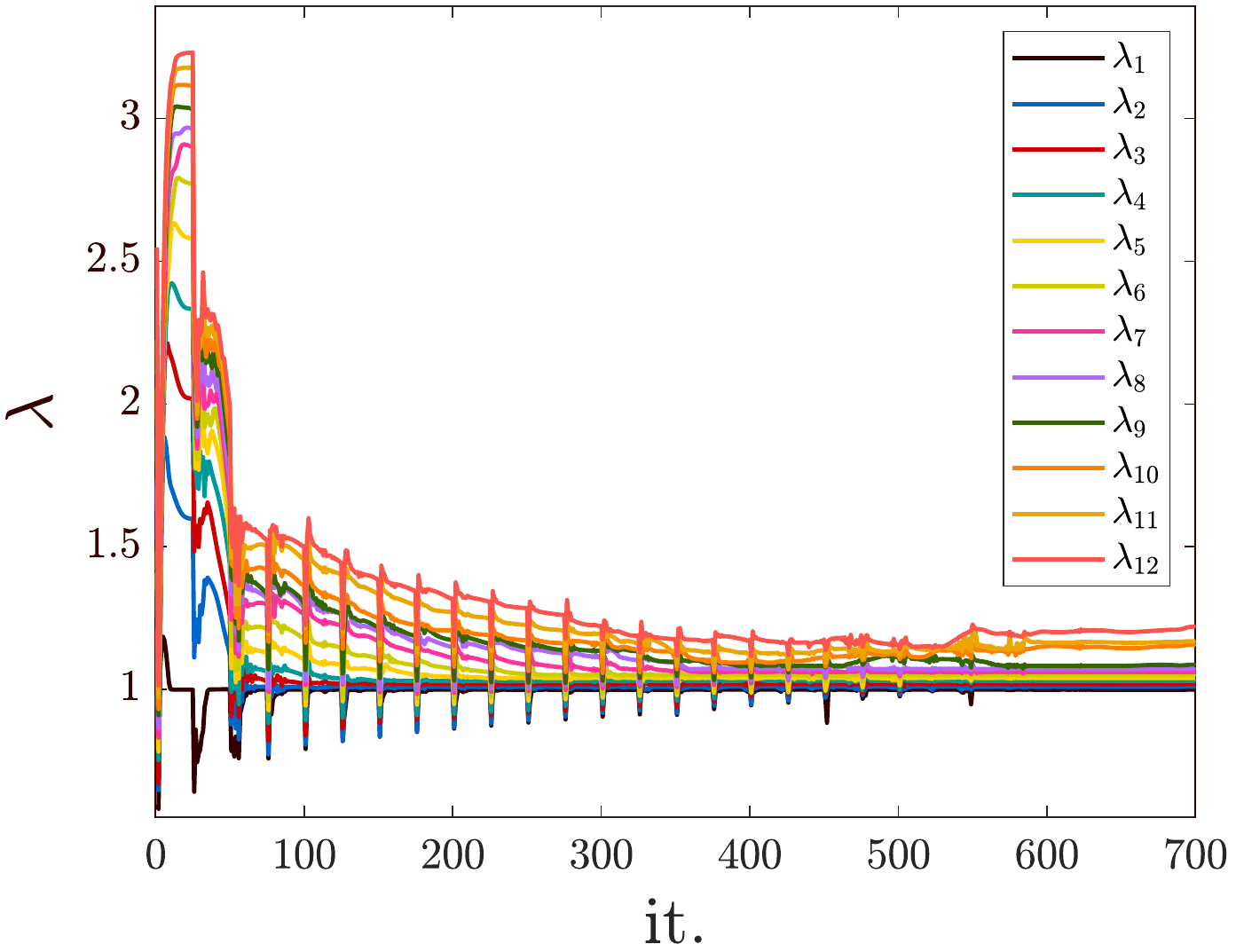}}
 \caption{Optimization histories for problem $\mathcal{P}_{1}$ solved on the discretization $\Omega_{1} = 840\times 360$. (a) shows the evolution of the normalized compliance (black curve) and of the fundamental BLF (blue curve). (b) shows the evolution of the BLFs corresponding to the constrained buckling modes}
 \label{fig:NumericalTest-2DExample-convergenceCurves-M4}
\end{figure}

%% TABLE AND FIGURE ABOUT COMPUTATIONAL TIME
\begin{figure}[t]
  \begin{minipage}[t]{0.575\linewidth}
   \vspace{0pt}
    \centering
     \begin{tabular}{c|cccccc}
      \hline\noalign{\smallskip}
                     & \multicolumn{3}{c}{$\Omega_{1} = 840 \times 360$} & \multicolumn{3}{c}{$\Omega_{1} = 1680 \times 720$} \\
      \noalign{\smallskip}\hline
              $\ell$ & \texttt{tLBA}(s) & $\texttt{sF}$ & $\texttt{eR}$ & \texttt{tLBA}(s) & $\texttt{sF}$ & $\texttt{eR}$ \\
                     \hline \\
                   1  & 120.0 & 1.0 & 0.981 & 649.6 &  1.0 & 0.979 \\
                   2  &  47.7 & 2.5 & 0.644 & 257.3 &  2.5 & 0.639 \\
                   3  &  15.1 & 8.0 & 0.587 &  69.1 &  9.4 & 0.568 \\
                   4  &     - &    - &     - & 42.5 & 15.3 & 0.512 \\
                     \noalign{\smallskip}\hline
     \end{tabular}
   \end{minipage}
   \begin{minipage}[t]{0.4\linewidth}
    \vspace{0pt}
     \centering
     \includegraphics[scale = 0.425, keepaspectratio]
         {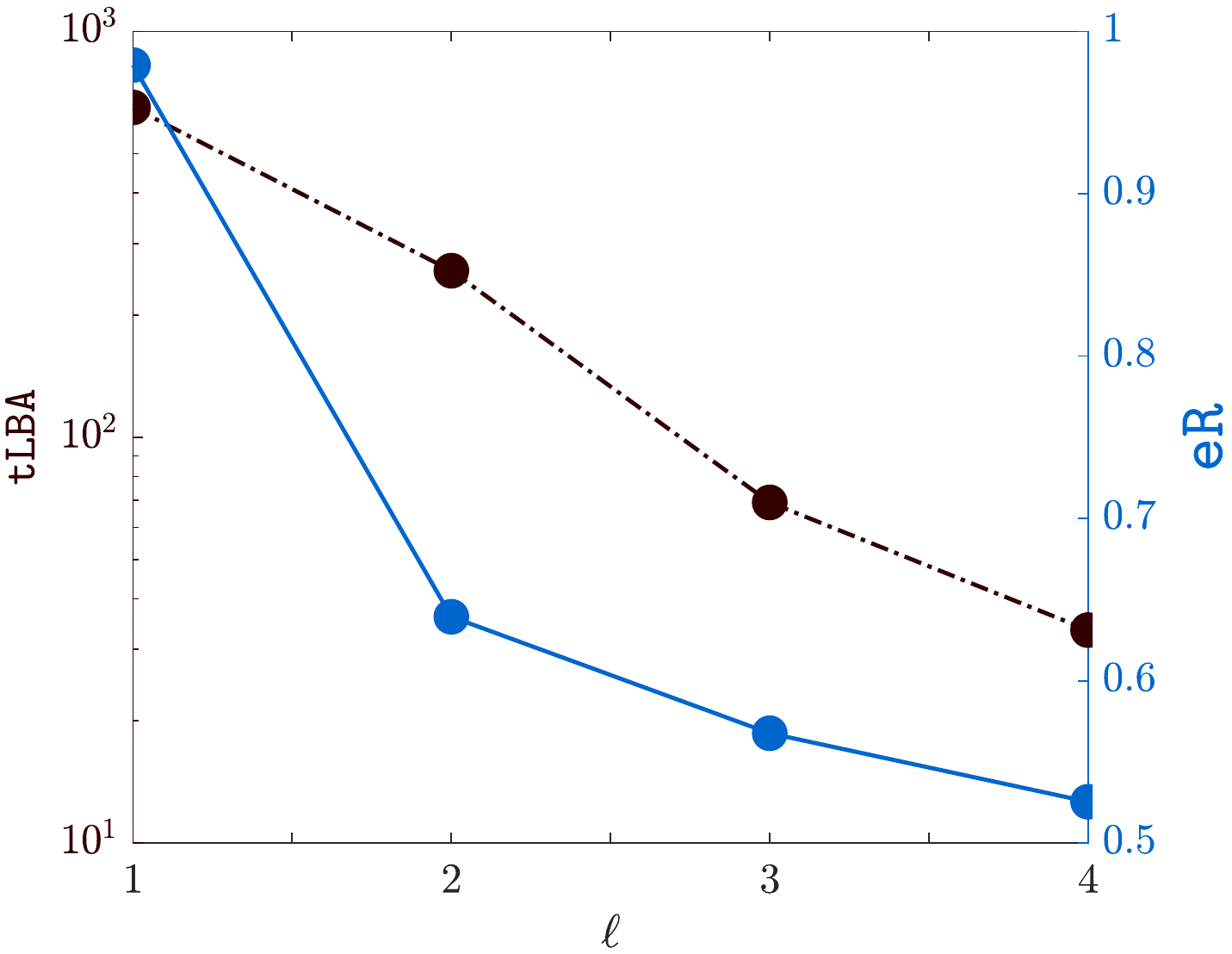}
   \end{minipage}%
\caption{Average time spent for a Linearized Buckling Analysis (\texttt{tLBA}) referred to the examples discussed in \autoref{Sec:2Dexample-2BarsTruss}. The plot refers to the discretization $\Omega_{1} = 1680 \times 720$. For $\ell = 1$ the LA and EA are performed on the fine discretization while $\ell > 1$ defines the coarse discretization for the multilevel procedure. $\texttt{sF} =\texttt{tLBA}(\ell)/\texttt{tLBA}(\ell = 1)$ are the saving factors and $\texttt{eR}$ represents the weight of the EA compared to the overall $\texttt{tLBA}$.}
\label{tab:tableComputational}
\end{figure}

%% Numerical results for the minimum compliance design on the discretization 840x360
\begin{table}
 \centering
  \begin{tabular}{ccl|cccc}
   \hline\noalign{\smallskip}
                     & \multicolumn{2}{c}{\textbf{Opt.}} & \multicolumn{4}{c}{\textbf{Post-Proc. (BW)}} \\
   \noalign{\smallskip}\hline
                 $i$ & $\lambda_{i}$  & $\delta_{i}$ ($10^{-4}$) & $\ell = 1$ & $\ell = 2$ & $\ell = 3$ & $\ell = 4$ \\
                    \hline \\
                   1 & 1.000 &      - & 1.001        & 1.002        & 1.005        & 1.021        \\
                   2 & 1.015 &   1.08 & 1.125        & 1.130        & 1.139        & 1.207        \\
                   3 & 1.026 &   1.46 & 1.134        & 1.262        & 1.271        & 1.328        \\
                   4 & 1.038 &   1.01 & 1.256        & 1.314        & 1.324        & 1.387        \\
                   5 & 1.045 &   0.52 & 1.264        & 1.331        & 1.351        & 1.430        \\
                   6 & 1.064 &   0.71 & 1.297        & 1.339        & 1.371        & 1.444        \\
                   7 & 1.071 &   0.77 & 1.306        & 1.359        & 1.383        & 1.486        \\
                   8 & 1.097 &   1.52 & 1.317        & 1.363        & 1.418        & 1.501        \\
                   9 & 1.111 &  43.90 & 1.329        & 1.401        & 1.423        & 1.544        \\
                  10 & 1.116 & 631.12 & 1.333        & 1.405        & 1.457        & 1.563        \\
                  11 & 1.166 & 643.31 & 1.349        & 1.410        & 1.488        & 1.586        \\
                  12 & 1.221 & 1,058  & 1.353        & 1.442        & 1.501        & 1.614        \\
                    \noalign{\smallskip}\hline
 \end{tabular}
\caption{Eigenvalues from the optimization (columns \textbf{Opt.}) and post--processed BW design for the structure from \autoref{fig:NumericalTest-2DExample}(b). The coarse level used in the optimization is $\ell = 3$. The relative difference between the compliance obtained from optimization, and those computed on the post--processed BW design for different $\ell$ is always below $10^{-6}$}
\label{tab:tableA}
\end{table}

\begin{table}
 \centering
  \begin{tabular}{ccl|ccccc|ccrc}
   \hline\noalign{\smallskip}
                    & \multicolumn{2}{c}{\textbf{Opt.}} & \multicolumn{4}{c}{\textbf{Post-Proc. (BW)}} &
                    & \multicolumn{2}{c}{\textbf{P--R}} \\
    \noalign{\smallskip}\hline
                $i$ & $\lambda_{i}$  & $\delta_{i}$ ($10^{-5}$) & $\ell = 1$ & $\ell = 2$ & $\ell = 3$ & $\ell = 4$ &
                & $\mathfrak{L}$ & $\Delta\lambda_{i}/\lambda_{i}$ $\%$ \\
                   \hline \\
                   1 & 0.997 &     - & 0.999        & 1.001        & 1.002       & 1.005  &&   1.0 & $\approx$ 0 \\
                   2 & 1.010 &  6.21 & 1.001        & 1.099        & 1.170       & 1.181  && 212.6 & 0.08 \\
                   3 & 1.021 &  7.28 & 1.066        & 1.166        & 1.235       & 1.249  && 184.4 & 0.11 \\
                   4 & 1.032 & 11.76 & 1.154        & 1.224        & 1.253       & 1.279  &&   3.1 & 0.31 \\
                   5 & 1.047 & 56.58 & 1.164        & 1.228        & 1.266       & 1.307  && 166.6 & 0.14 \\
                   6 & 1.060 & 90.63 & 1.216        & 1.240        & 1.293       & 1.321  &&  41.7 & 2.28 \\
                   7 & 1.094 & 320.4 & 1.229        & 1.263        & 1.305       & 1.356  &&   3.9 & 2.19 \\
                   8 & 1.098 & 251.3 & 1.233        & 1.289        & 1.311       & 1.394  &&   2.0 & 1.57 \\
                   9 & 1.117 & 338.9 & 1.260        & 1.304        & 1.334       & 1.404  && 150.6 & 1.65 \\
                  10 & 1.136 & 420.5 & 1.285        & 1.327        & 1.360       & 1.416  &&   8.8 & 2.02 \\
                  11 & 1.139 & 334.8 & 1.298        & 1.348        & 1.377       & 1.435  &&   4.6 & 3.09 \\
                  12 & 1.162 & 458.5 & 1.302        & 1.357        & 1.381       & 1.459  &&  13.1 & 2.0 \\
                  \noalign{\smallskip}\hline
 \end{tabular}
\caption{Eigenvalues from the optimization (columns \textbf{Opt.}) and post--processed BW design for the structure from \autoref{fig:NumericalTest-2DExample}(c). The coarse level used in the optimization is $\ell = 4$. The relative difference between the compliance obtained from optimization, and those computed on the post--processed BW design for different $\ell$ is always below $10^{-6}$. The columns below \textbf{P--R} display the TV measure of the modal strain energy, according to \eqref{eq:TVsed}, and the change in the BLFs after the thickening operation $\mathcal{R}[\mathbf{x}]$ \eqref{eq:postProcessingDilation}}
\label{tab:tableB}
\end{table}

\section{Design of a 2D structure for minimum compliance}
 \label{Sec:2Dexample-2BarsTruss}

Let us consider the geometry sketched in \autoref{fig:NumericalTest-2DExample} (a), originally discussed in \cite{ferrari-sigmund_19a}. Points $a$ and $b$ are hinged and a downward load, having total magnitude $|F| = 2 \cdot 10^{-2}$, is spread over a length of $L_{x}/10$ near points $c$. Square regions near these three points, with dimension $L_{x}/10$, are set to be solid during the optimization. The values of Young's moduli used in \autoref{eq:interpolationSIMP2-SIMP} are $E_{1} = 1$ and $E_{0}=10^{-6}$ over the design domain, while $E_{p} = 10^{3}$ for the prescribed solid regions (to alleviate problems with local stress concentrations at load and supports).

We address the compliance minimization problem for a maximum volume fraction $\bar{f} = 0.16$ and minimum BLF $\bar{\lambda} = 1.0$

\begin{equation}
 \label{eq:optProblemMinCpl}
 \mathcal{P}_{1}
  \begin{cases}
   & \min\limits_{\hat{\mathbf{x}}\in [ 0, 1 ]^{m}}
   J\left( \mathbf{x} \right) = \mathbf{u}^{T} K\left[ \mathbf{x} \right] \mathbf{u} \\
   {\rm s.t.} & \min\limits_{i\in\mathcal{B}}\tilde{\lambda}_{i}
   \geq \bar{\lambda} \\
   & V\left( \mathbf{x} \right) \leq
   \bar{f} |\Omega_{h}|
 \end{cases}
\end{equation}

We recall that, due to linearity, the compliance can be evaluated in any reference state and it is convenient to consider the state $\left(\mathbf{f}, \mathbf{u}\right)$ already solved for in connection with the linearized buckling analysis (Step 2 in \autoref{alg:EigenvalueBucklingAnalysis}). The buckling constraint has been implemented with a bound formulation \cite{bendsoe_89a}, imposing a small gap between eigenvalues in order to prevent their complete coalescence. Specifically, the constraint $\min_{i\in\mathcal{B}}\tilde{\lambda}_{i}$ is replaced by the set

\[
 \alpha^{i}\tilde{\lambda}_{i}/\bar{\lambda} - 1 \leq 0 \ , \qquad i\in\mathcal{B}
\]
where $\alpha = 0.99$, and the lowest 12 buckling modes are considered within the optimization. However, eigenpairs up to the $24^{\rm th}$ are still computed in this test problem for monitoring purposes. No substantial differences have been observed if instead considering aggregation of these constraints (with e.g. $p$--norm or Kreisselmeier--Steinhauser \cite{kreisselmeier-steinhauser_79a} functions), provided that the aggregation parameter is chosen high enough \cite{ferrari-sigmund_19a}.

The optimization problem is run for 700 steps, increasing the penalization $p$ from $1$ to $6$ each 25 steps, with $\Delta p = 0.25$. The filter radius is $r_{\rm min} = 8h$ and the projection parameters are fixed to $\eta = 0.5$ and $\beta = 6$. The Method of Moving Asympotes \cite{svanberg_87a} is used to update the design variables $\hat{\mathbf{x}}$. \autoref{fig:NumericalTest-2DExample} (b) and (c) show the optimized designs corresponding to the two different fine discretizations $\Omega_{1} = 840 \times 360$ ($6.07\cdot 10^{5}$ DOFs) and $\Omega_{1} = 1680 \times 720$ ($2.424\cdot10^{6}$ DOFs). The coarse levels for the multilevel procedure are set to $\ell = 3$ and $\ell = 4$, with cut in the DOFs number of 16 and 64 times, respectively.

The effect of a finer discretization, increasing the design freedom, is clearly seen with a more complex distribution of thinner bars, for the same value of $\bar{\lambda}$ and $\bar{f}$. The optimization progress is shown in \autoref{fig:NumericalTest-2DExample-convergenceCurves-M4} (a, b), referring to the structure of \autoref{fig:NumericalTest-2DExample} (b). The initial compliance and fundamental BLF are $J_{(0)} = 2.501 \cdot 10^{-3}$ and $\lambda_{1(0)} = 0.597$, respectively and we see that in the beginning of the optimization $\lambda_{1}$ is increased to quickly meet the buckling constraints. The compliance is also reduced, as we are considering $p = 1$. Then, as $p$ is raised, the compliance increases as the optimizer strives to fullfill the buckling constraint. The distinct jumps in the compliance and BLFs evolution curves correspond to increases of the penalization parameters, while in the last 150 steps, when $p = 6$, the value of $\lambda_{1}$ is stable. The optimized design has a compliance of $J_{(700)} = 2.743 \cdot 10^{-3}$, about $7.9 \%$ higher than the initial one, and both the volume and buckling constraints are active.

From \autoref{fig:NumericalTest-2DExample-convergenceCurves-M4} (b) we clearly notice the activation of more and more BLFs as the optimization progresses and the following quantity \cite{ferrari-sigmund_19a}

\begin{equation}
 \label{eq:coalescingMeasure}
  \delta_{i} = \lambda_{i}/\lambda_{1} - \alpha^{(i-1)}
  \ , \qquad i = 2, \ldots, |\mathcal{B}|
\end{equation}
can be used to quantify this coalescing phenomenon. At the end of the optimization, coefficients $\delta_{2}$ to $\delta_{8}$ are below $10^{-4}$ and therefore the corresponding modes can be considered active. Moreover, $\lambda_{9}$ to $\lambda_{12}$ are also very close to the active set (see \autoref{tab:tableA}).

Similar observations apply to the design corresponding to the finer discretization $\Omega_{1} = 1680\times 720$ in \autoref{fig:NumericalTest-2DExample}(c). The compliance for this design is $J_{(700)} = 2.789\cdot 10^{-3}$ and the value of the 12 lowest BLFs is reported in the second column of \autoref{tab:tableB}. We recognize that the more complicated structural pattern is associated with the coalescence of more BLFs and the activation of more buckling modes. Coefficients $\delta_{2}$ to $\delta_{6}$ are below $10^{-4}$ and those up to $\delta_{12}$ are below $10^{-2}$.

Now, refer to \autoref{tab:tableComputational} and the Table therein in order to discuss computational savings. The time spent performing the Linearized Buckling Analysis ($\texttt{tLBA}$) directly on the fine scale $\Omega_{1}$, using direct solvers for both the Linear Analysis (LA) and the Eigenvalue Analysis (EA), is compared to that for the multilevel approach (starting from $\ell > 1$). The simulations have been performed with a laptop equipped with an Intel(R) Core(TM) i7-5500U@2.40GHz CPU, 15GB of RAM and Matlab 2018b, running in serial mode.

As the multilevel approximation starts from a coarser mesh, the computational savings become apparent. For the case $\Omega_{1} = 840\times 360$, choosing $\ell = 3$ we cut the computational time by 8 times and for the case $\Omega_{1} = 1680\times 720$ and $\ell = 4$ the cut is reaching 15 times. Moreover, and most importantly, we can back up our main claim: the cost for obtaining the approximation of the fine scale buckling modes approaches that of solving the linear system. Indeed, from the ratio $\texttt{eR} = \texttt{tEA}/\texttt{tLBA}$, we see how the cost for the EA and that for the LA become almost equal (\texttt{eR} approaches 0.5) as $\Omega_{\ell}$ becomes coarser (see plot and Table in \autoref{tab:tableComputational}).

%% BUCKLING MODES DISTRIBUTION ON DIFFERENT SCALES FOR THE DESIGN OBTAINED ON M4
\begin{figure}[tb]
 \centering
  \includegraphics[scale = 0.115, keepaspectratio]
   {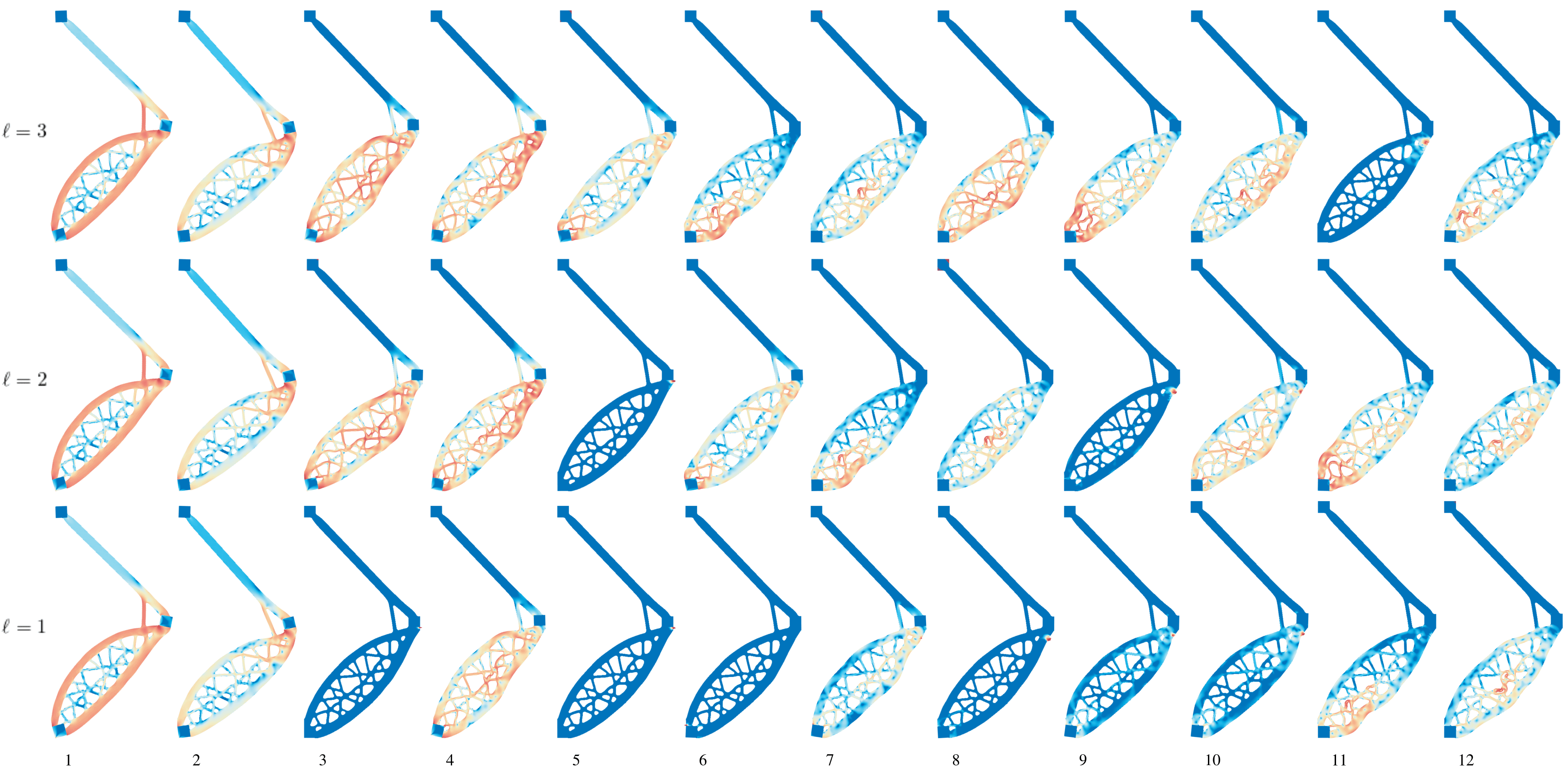}
 \caption{Buckling modes $\boldsymbol{\varphi}_{i}$ computed by the multilevel procedure starting from different coarse levels $\ell$ for the structure of \autoref{fig:NumericalTest-2DExample}(b) projected to a completely Black and White (BW) design. $\ell = 3$ is the coarse level used for running the optimization and $\ell = 1$ refers to a full LBA on the finest discretization. Coloring (blue to red) corresponds to the strain energy density distribution (low to high), in logarithmic scale}
 \label{fig:NumericalTest-2DExample-modesComparison-M4}
\end{figure}

%% BUCKLING MODES DISTRIBUTION ON DIFFERENT SCALES FOR THE DESIGN OBTAINED ON M8
\begin{figure}[tb]
 \centering
  \includegraphics[scale = 0.11, keepaspectratio]
   {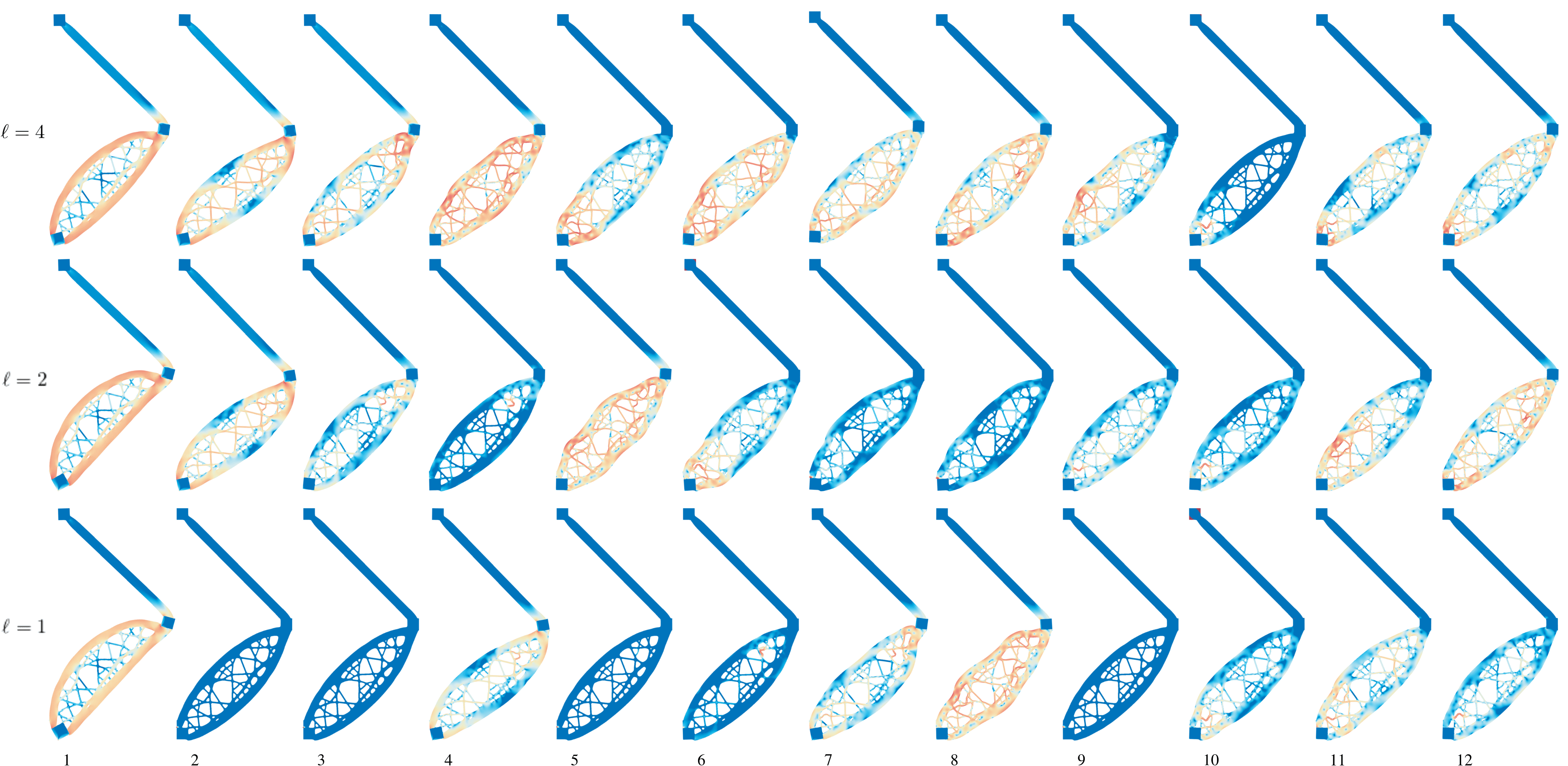}
 \caption{Buckling modes $\boldsymbol{\varphi}_{i}$ computed by the multilevel procedure starting from different coarse levels $\ell$ for the structure of \autoref{fig:NumericalTest-2DExample}(c) projected to a completely Black and White (BW) design. $\ell = 4$ is the coarse level used for running the optimization and $\ell = 1$ refers to a full LBA on the finest discretization. Coloring (blue to red) corresponds to the strain energy density distribution (low to high), in logarithmic scale}
 \label{fig:NumericalTest-2DExample-modesComparison-M8}
\end{figure}

\section{Post--processing of the obtained designs}
 \label{Sec:AnalysisDesigns2BS}

To validate the designs obtained in \autoref{Sec:2Dexample-2BarsTruss}, we perform a full linearized buckling analysis on $\Omega_{1}$ using direct solvers for both the linear and eigenvalue equations. The designs show some grayscale, which can be quantified by the non--discreteness measure \cite{sigmund_07a}

\begin{equation}
 \label{eq:discretenessMeasure}
  m_{\rm  nd} = \frac{4}{m}\mathbf{x}^{T}\left(1 - \mathbf{x}\right) \approx 2.3\%
\end{equation}

To rule out any effect due to grayscales, we first recover a completely Black and White (BW) design by means of a sharp Heaviside projection with $\eta = 0.5$. We stress that even if the projection operation does not produce any noticeable change neither in the topology, nor in the compliance values, it does affect the BLFs, their distribution and the associated buckling modes (see \autoref{tab:tableA} and \autoref{tab:tableB}). This is expected, as buckling response is generally very sensitive to structural modifications, and this feature is even sharpened for an optimized design.

%% MAC CORRELATION FACTORS COMPUTED FOR THE TWO OPTIMIZED DESIGNS AND TAKING THE FINE GRID MODES AS REFERENCE ONES
\begin{figure*}[tb]
 \centering
  \subfloat[Design of \autoref{fig:NumericalTest-2DExample}(b)]{
   \includegraphics[scale = 0.425, keepaspectratio]
   {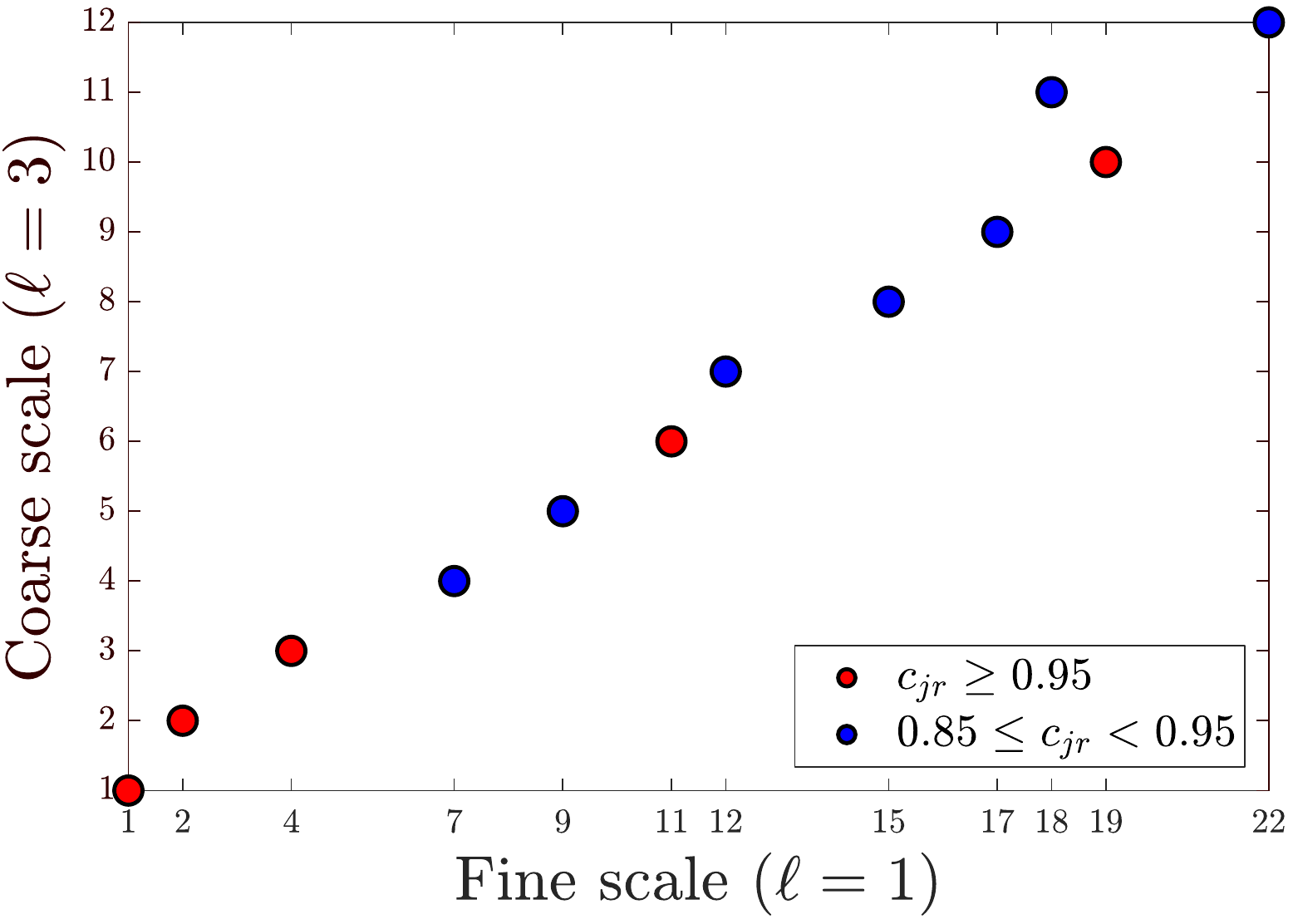}} \quad
  \subfloat[Design of \autoref{fig:NumericalTest-2DExample}(b)]{
   \includegraphics[scale = 0.425, keepaspectratio]
   {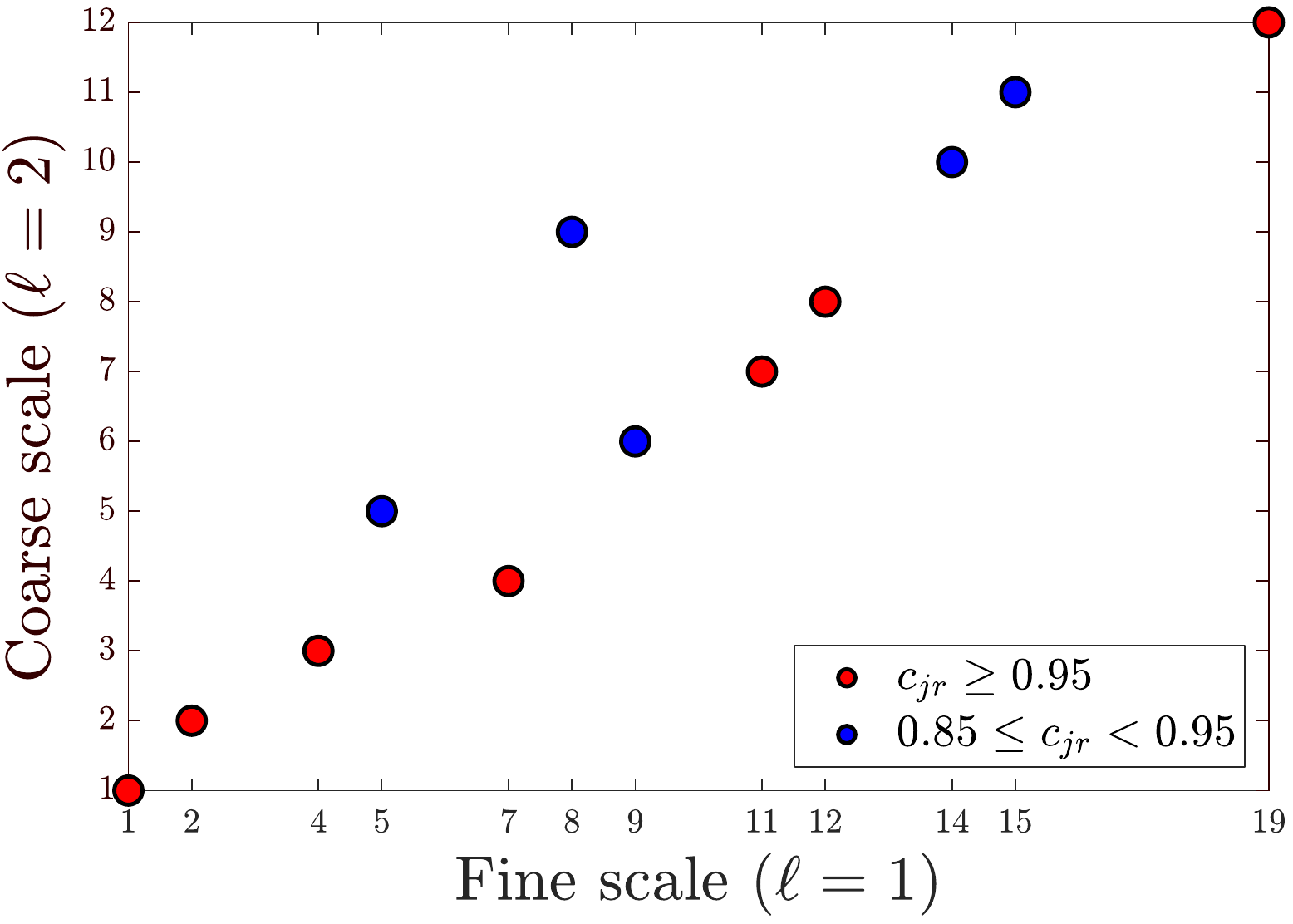}} \\
  \subfloat[Design of \autoref{fig:NumericalTest-2DExample}(c)]{
   \includegraphics[scale = 0.425, keepaspectratio]
   {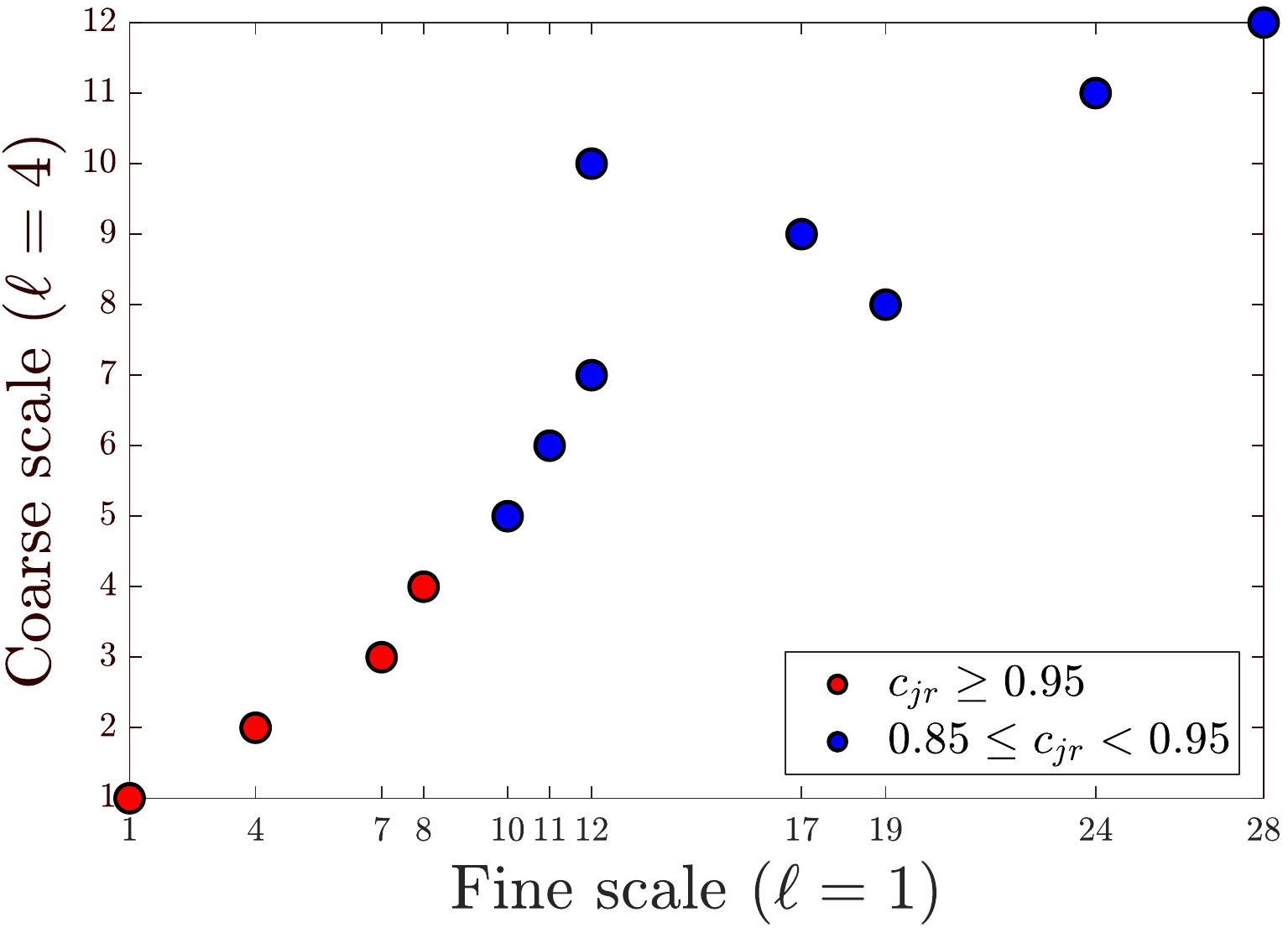}} \quad
  \subfloat[Design of \autoref{fig:NumericalTest-2DExample}(c)]{
   \includegraphics[scale = 0.425, keepaspectratio]
   {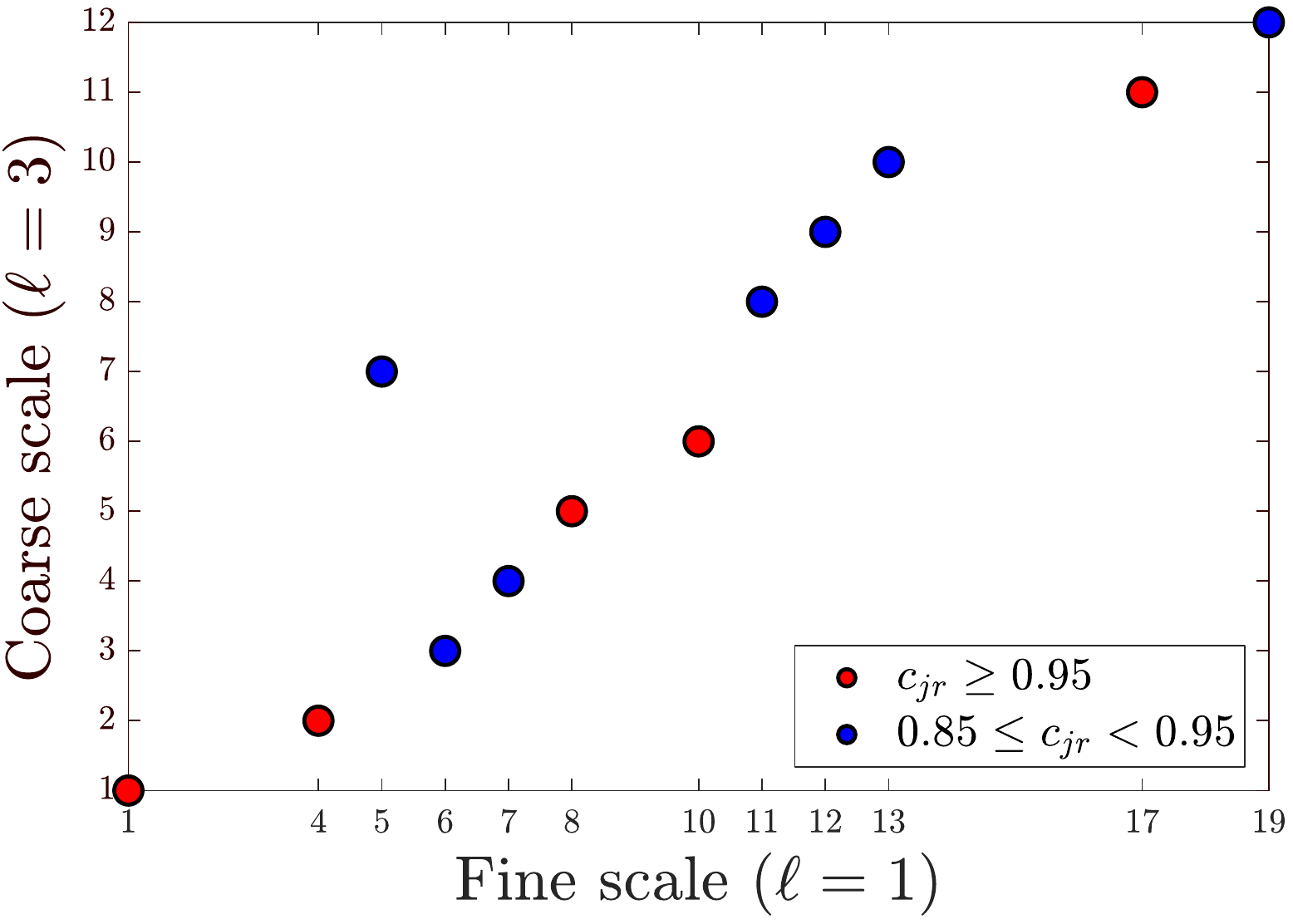}}
 \caption{Distribution of coefficients $c_{jr}$ defined by \eqref{eq:macCorrelationIndex} for the two designs of \autoref{fig:NumericalTest-2DExample}(b,c). Values $c_{jr} > 0.95$ are shown in red and values $c_{jr} \in [0.85, 0.95]$ are shown in blue, while $c_{jr} < 0.85$ are not represented.}
 \label{fig:NumericalTest-2DExample-macPlot-M4}
\end{figure*}

%% ACCURACY MEASURE PLOT FOR MODES APPROXIMATED STARTING FROM SEVERAL COARSE GRIDS
\begin{figure}[t]
 \centering
    \subfloat[Design of \autoref{fig:NumericalTest-2DExample}(b)]{
     \includegraphics[scale = 0.45, keepaspectratio]
   {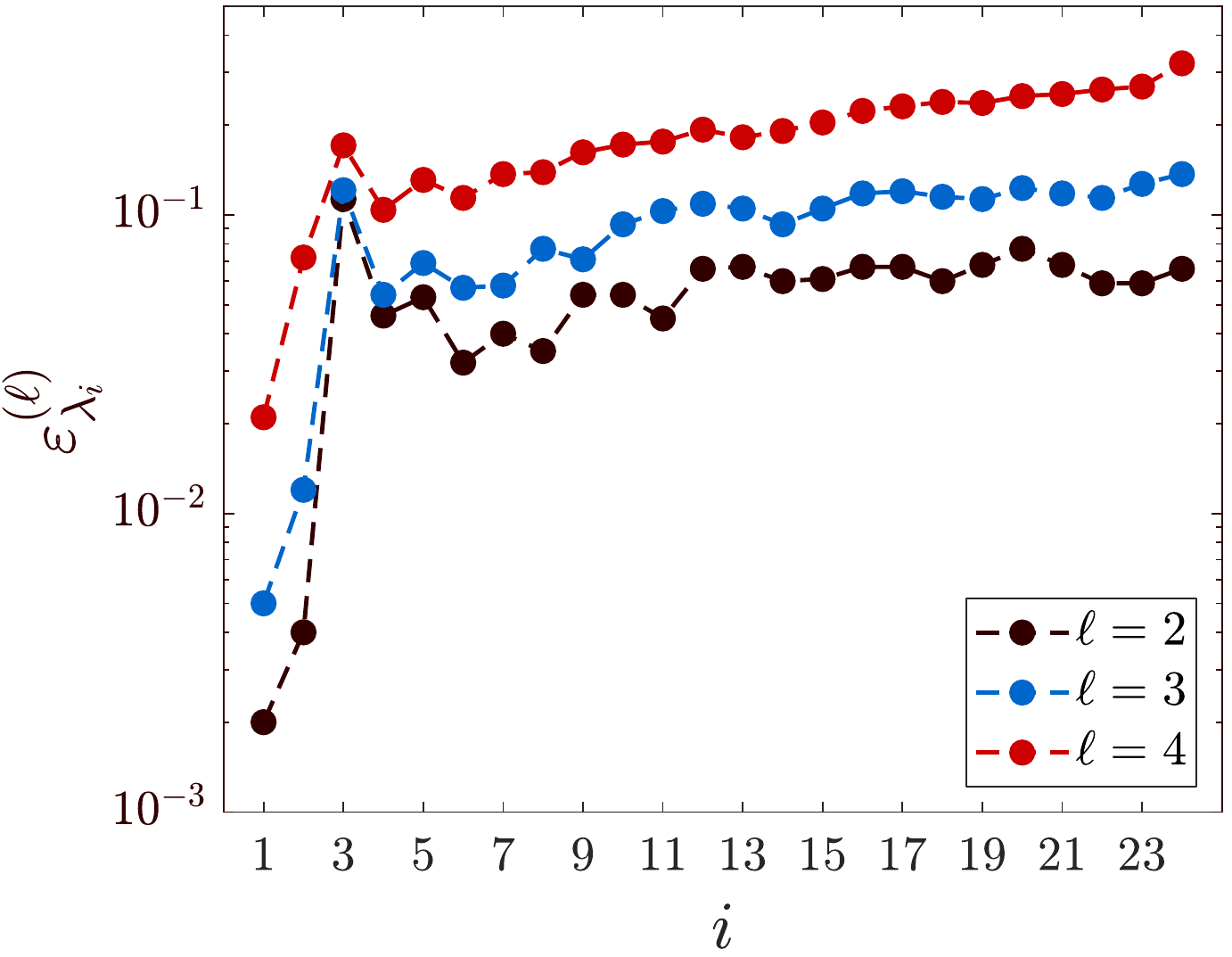}} \qquad
    \subfloat[Design of \autoref{fig:NumericalTest-2DExample}(b)]{
     \includegraphics[scale = 0.45, keepaspectratio]
   {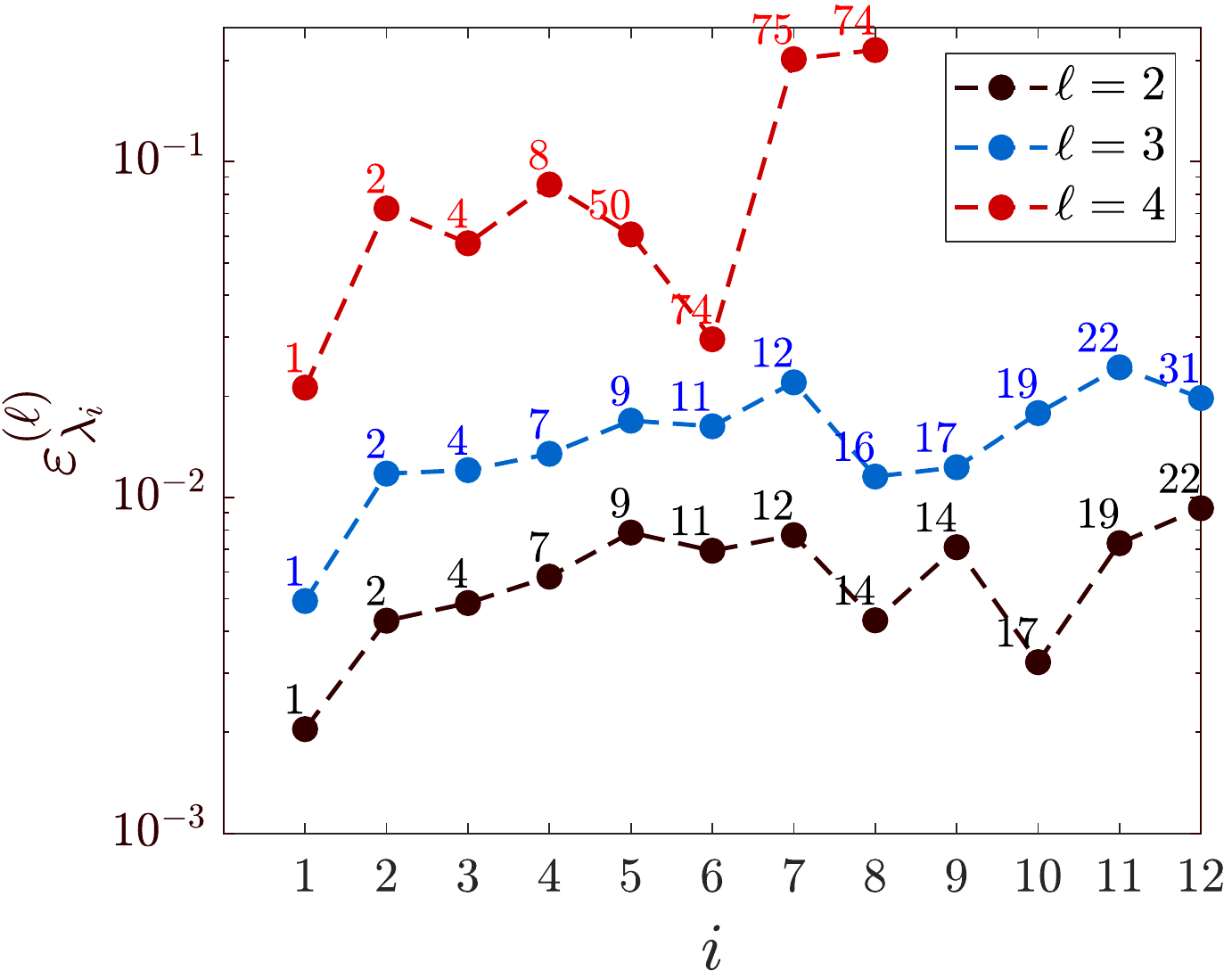}} \\
    \subfloat[Design of \autoref{fig:NumericalTest-2DExample}(c)]{
     \includegraphics[scale = 0.45, keepaspectratio]
   {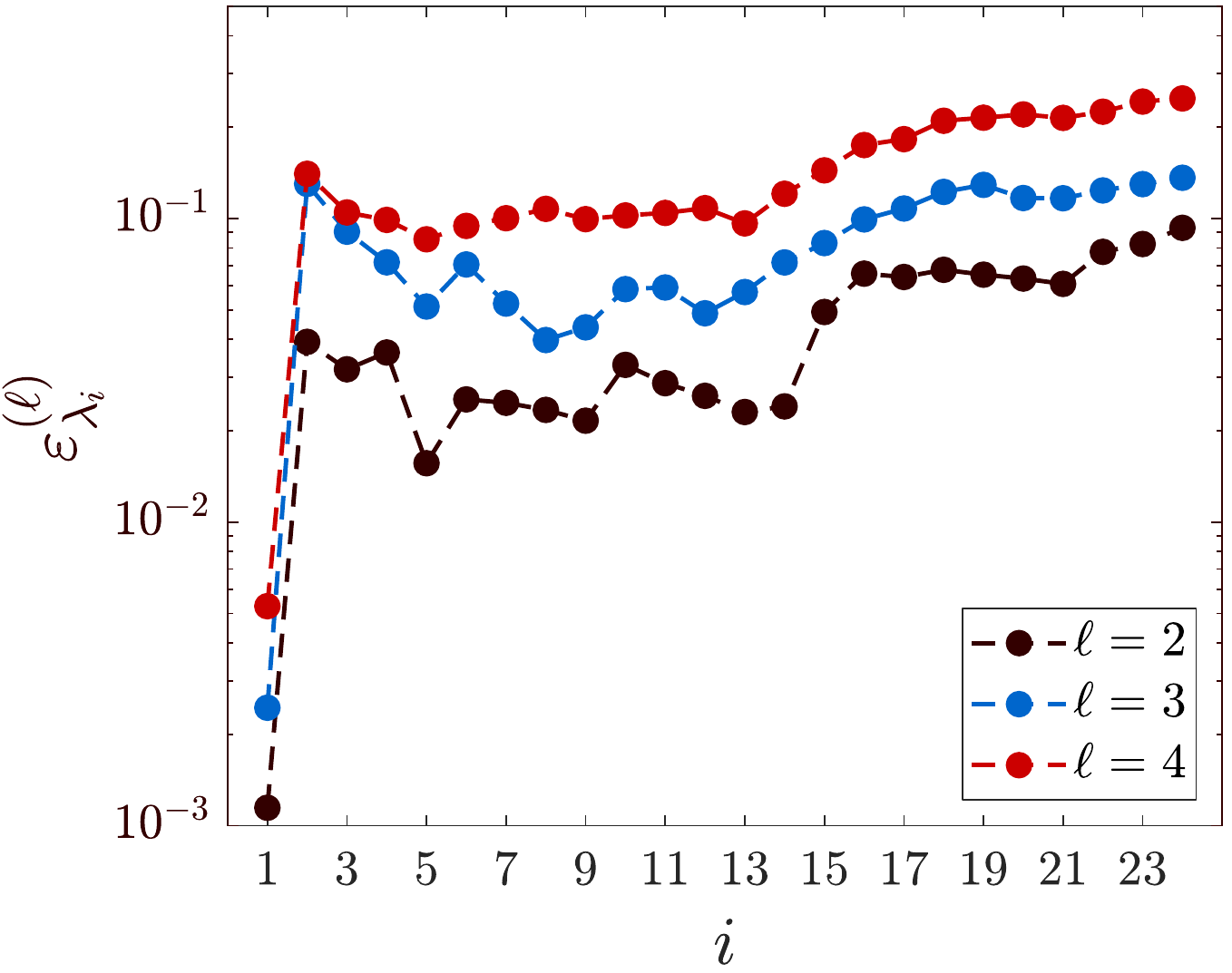}} \qquad
    \subfloat[Design of \autoref{fig:NumericalTest-2DExample}(c)]{
     \includegraphics[scale = 0.45, keepaspectratio]
   {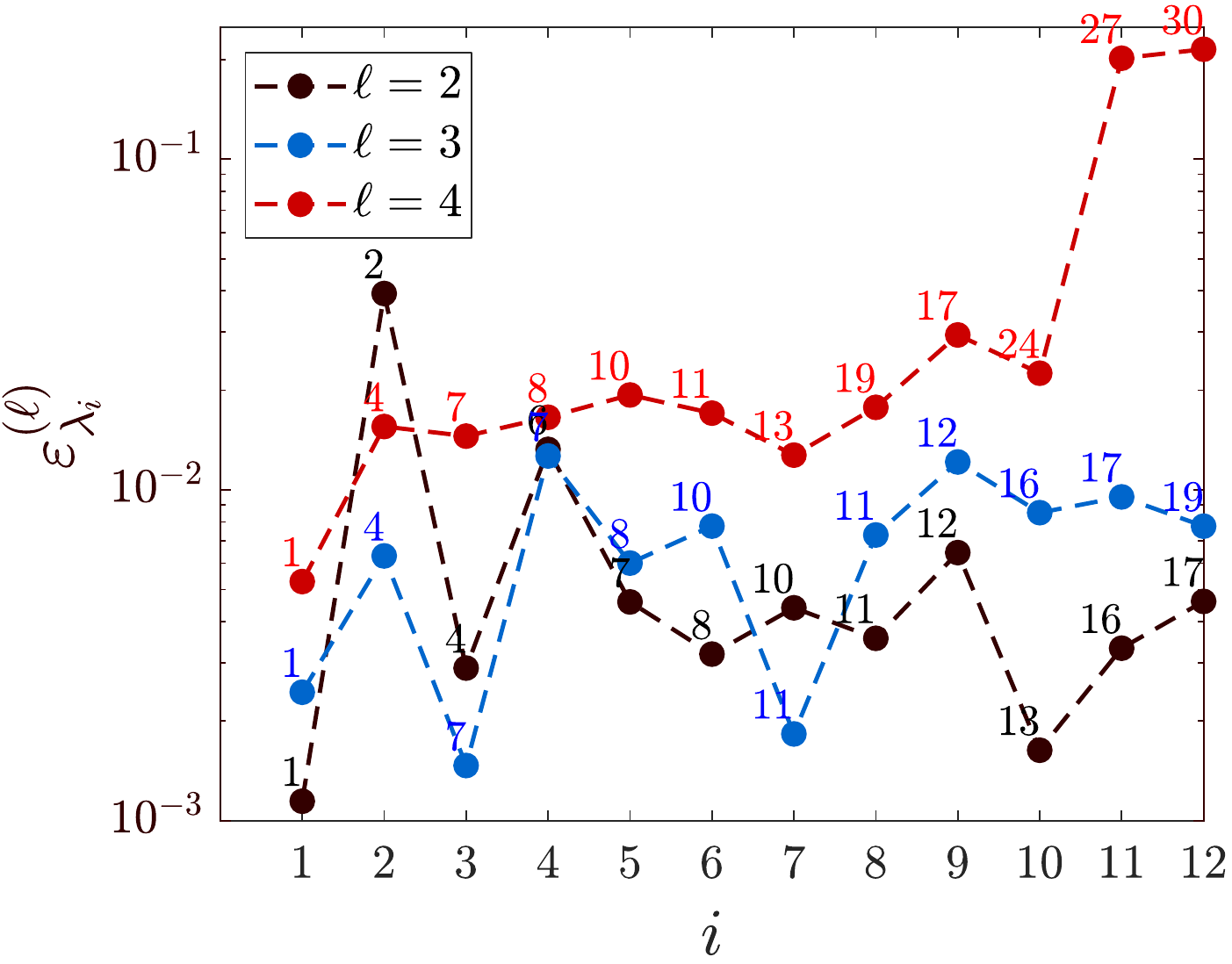}}
  \caption{Behavior of the error measure defined in \eqref{eq:relativeBLFandCplErrorMeasures} as the multilevel procedure starts from different coarse levels $\ell$. In plots (a,c) all the computed BLFs are considered and sorted according to magnitude when applying \eqref{eq:relativeBLFandCplErrorMeasures}. In plots (b,d) only the BLFs of modes that can be paired through the MAC ($c_{jr} \geq 0.95$) are shown and the corresponding sorting is used when applying \eqref{eq:relativeBLFandCplErrorMeasures}.}
 \label{fig:dataPostProcessing_2BT_minCpl_M4}
\end{figure}

Let us start analysing the design of \autoref{fig:NumericalTest-2DExample}(b). The active buckling modes for the BW design, computed again starting from $\ell = 3$, are shown in the top row of \autoref{fig:NumericalTest-2DExample-modesComparison-M4}, and we remark that the fundamental buckling mode remains unchanged after the projection operation. Numerical values of the lowest 12 BLFs are listed in \autoref{tab:tableA} (see  column $\ell = 3$). The value of $\lambda_{1}$ is kept essentially the same, whereas the higher BLFs all increase by $0.7$ to $22\%$. Also, the coefficients $\delta_{i}$ are substantially increased, in the order of $10-21\%$; therefore gaps between eigenvalues are widened.

Then we analyze the design starting from different coarse grids up to $\ell = 1$, which means that the LBA is directly performed on the fine scale $\Omega_{1}$. In order to keep track of the modification or switching of buckling modes between grids, we compute the following coefficient $c_{jr} \in [0, 1]$

\begin{equation}
 \label{eq:macCorrelationIndex}
  c_{jr} = \tilde{\boldsymbol{\varphi}}^{T}_{j}
  K[\mathbf{x}]\boldsymbol{\varphi}_{r}
\end{equation}
which, considering the normalization $\tilde{\boldsymbol{\varphi}}^{T}_{j}K[\mathbf{x}]\tilde{\boldsymbol{\varphi}}_{j} = 1$ $\forall \: j$, is a variant of the well--known Modal Assurance Criterion (MAC) \cite{allemang-brown_82a} accounting for the modal strain energy \cite{brehm-etal_10a}. Equation \eqref{eq:macCorrelationIndex} measures the degree of similarity of mode $\tilde{\boldsymbol{\varphi}}_{j}$, computed starting from a coarse scale $\Omega_{\ell}$, to a mode $\boldsymbol{\varphi}_{r}$ on the fine scale $\Omega_{1}$. We have $c_{jr} = 0$ for two orthogonal modes and $c_{jr} \rightarrow 1$ as $\boldsymbol{\varphi}_{j}$ resembles $\boldsymbol{\varphi}_{r}$ more closely. It is emphasized that, in order for the MAC to be reliable, only values very close to 1 (e.g. $c_{jr} > 0.95$) should be accepted \cite{allemang_02a}. Coefficients $c_{jr}$ computed for the current example are shown in \autoref{fig:NumericalTest-2DExample-macPlot-M4}(a, b).

Combining the information of the plots with the buckling modes represented in \autoref{fig:NumericalTest-2DExample-modesComparison-M4} we can observe the following
\begin{enumerate}
 \item Some very localized modes appear when performing the LBA directly on the finest level $\ell = 1$ (see $\boldsymbol{\varphi}_{3}$, $\boldsymbol{\varphi}_{5}$, $\boldsymbol{\varphi}_{6}$ and $\boldsymbol{\varphi}_{8}$ in \autoref{fig:NumericalTest-2DExample-modesComparison-M4}), whereas these are not found with the multilevel strategy starting from $\ell = 3$. Also the coefficients $c_{jr}$ in \autoref{fig:NumericalTest-2DExample-macPlot-M4} (a) show how these very localized modes do not correlate with any of the modes computed starting from the coarse level;
 \item On the intermediate level $\ell =2$ there are still some localized effects (e.g. $\boldsymbol{\varphi}_{6}$ and $\boldsymbol{\varphi}_{9}$) which can be correlated with some of those on the fine level, but with a lower level of confidence ($c_{jr} \in [0.85, 0.95]$). This is another indication of the mesh dependency of these modes, which therefore should not be trusted as physically meaningful;
 \item We emphasize that, if $\Omega_{\ell}$ is choosen coarse enough, with the multilevel procedure we can compute considerably less buckling modes and still represent the physical and global ones. This is clear from \autoref{fig:NumericalTest-2DExample-macPlot-M4} (a,c), where the 12 modes computed from the coarse scale match with global and physically meaningful fine scale modes up to the 28th and 24th, respectively.
\end{enumerate}

The same post--processing operations have been carried out for the design of \autoref{fig:NumericalTest-2DExample} (c), and results are reported in \autoref{tab:tableB}. Similar conclusions can be drawn regarding the effect of the BW projection, and the computed buckling modes are shown in \autoref{fig:NumericalTest-2DExample-modesComparison-M8}. Again, we recognize some extremely localized deformations among the modes computed on the finest scale (see $\boldsymbol{\varphi}_{2}$, $\boldsymbol{\varphi}_{3}$, $\boldsymbol{\varphi}_{5}$ and $\boldsymbol{\varphi}_{9}$). In addition to these, now there are other local modes associated with the failure of single bars (see $\boldsymbol{\varphi}_{6}$, $\boldsymbol{\varphi}_{10}$ etc.) and also their occurence is shifted to higher modes as we increase the coarse level $\ell$. Again, the coefficients $c_{jr}$ shown in \autoref{fig:NumericalTest-2DExample-macPlot-M4} (c,d) tell us that very localized modes on $\Omega_{1}$ do not match with any modes on coarser grids. On the other hand the local, single--bar failures have a match, with a low $c_{jr}$ value (see blue points in \autoref{fig:NumericalTest-2DExample-macPlot-M4} (c,d)) and therefore they are very sensitive to the mesh fineness. From \autoref{fig:NumericalTest-2DExample-macPlot-M4} (c) we notice that the fine scale mode $\boldsymbol{\varphi}_{12}$ is matched (with low confidence) simultaneously by two modes ($\boldsymbol{\varphi}_{7}$ and $\boldsymbol{\varphi}_{10}$) computed through the multilevel procedure. This has to be expected, especially for modes involving a single bar, as many of these might appear and have similar BLFs (i.e. similar strain energy).

We can discuss the accuracy of the BLFs approximations referring to the error measure

\begin{equation}
 \label{eq:relativeBLFandCplErrorMeasures}
  \varepsilon^{(\ell)}_{\lambda_{i}} =
   1 - \tilde{\lambda}^{(\ell)}_{i}/\lambda_{i}
\end{equation}
where $\tilde{\lambda}^{(\ell)}_{i}$ are the BLFs approximations given by the Rayleigh quotient \eqref{eq:FineScaleRQ}, with modes computed starting from a certain coarse scale, and $\lambda_{i}$ are those directly computed on the fine scale. The trend of this quantity is shown in \autoref{fig:dataPostProcessing_2BT_minCpl_M4}. We immediately appreciate how $\lambda_{1}$, which in this context is the safety measure for the obtained design, possesses a very good accuracy ($\varepsilon^{(\ell)}_{\lambda_{1}} \approx 0.2- 0.5\%$ for both designs). The accuracy seems to rapidly deteriorate for higher buckling modes, especially if comparing them by keeping the original sorting according to their magnitude (see plots (a,c)). However, if we compare only the modes which can be paired by the MAC ($c_{jr} \geq 0.95$), sorting them accordingly, which is much more meaningful from a physical point of view, we have a good approximation of all the BLFs associated with global modes (see \autoref{fig:dataPostProcessing_2BT_minCpl_M4} (b,d)). From \autoref{fig:dataPostProcessing_2BT_minCpl_M4} (b) we notice how the accuracy of BLFs computed starting from $\ell = 4$ is considerably lower. Also, many relevant modes are missed (see the gap between $\boldsymbol{\varphi}_{8}$ and $\boldsymbol{\varphi}_{50}$). This means that a too coarse scale has been selected to adequately represent the behavior of the inner frame of bars. Therefore, a tradeoff between accuracy and computational efficiency must be expected, at some level, when applying the multilevel procedure of \autoref{sSec:multilevelBucklingModes}.

We can conclude that by using the multilevel approximation strategy we can not only capture the most meaningful and global buckling modes, but also associate them with BLFs approximations with suitable accuracy.

%% PLOT SHOWING THE EFFECT OF DISCRETIZATION COARSENING ON THE STRESS DISTRIBUTION. NOW IT COULD BE THE PLOT OF THE "ISOTROPIC" MEASURE OF THE STRESS STIFFNESS MATRIX.
\begin{figure}[tb]
 \centering
  \subfloat[$\Omega_{1}$, $\zeta_{1} = 1.93$]{
   \includegraphics[scale = 0.45, keepaspectratio]
   {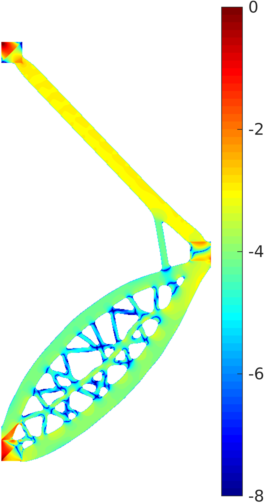}} \qquad
  \subfloat[$\Omega_{2}$, $\zeta_{2} = 0.83$]{
   \includegraphics[scale = 0.45, keepaspectratio]
   {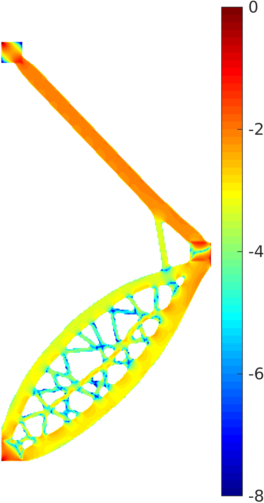}} \qquad
  \subfloat[$\Omega_{3}$, $\zeta_{3} = 0.41$]{
   \includegraphics[scale = 0.45, keepaspectratio]
   {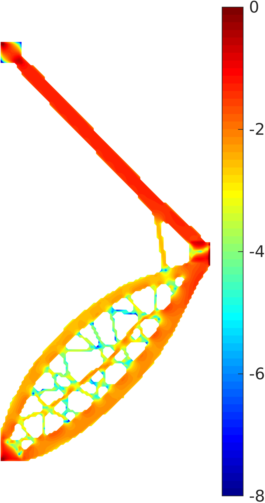}}
 \caption{Log--modulus plots of the scalar measure of the stress stiffness matrix defined by \eqref{eq:isomeasG} (normalized w.r.t. the maximum value) on the three discretization levels for the structure of \autoref{fig:NumericalTest-2DExample}(b)}
 \label{fig:isotropicMeasureStress}
\end{figure}

% ILLUSTRATION OF THE LOCAL POST-REINFORCEMENT PROCESS ON THE 2D EXAMPLE ON THE FINE SCALE
\begin{figure}[t]
 \centering
  \subfloat[$|\nabla \pi_{6}|$]{
   \includegraphics[scale = 0.4, keepaspectratio]
   {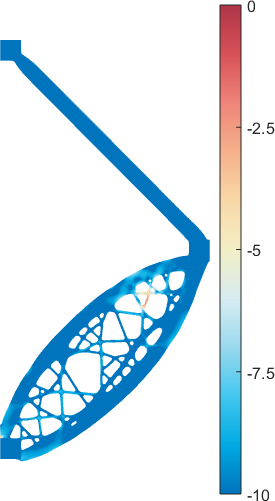}} \quad
  \subfloat[$|\nabla \pi_{10}|$]{
   \includegraphics[scale = 0.4, keepaspectratio]
   {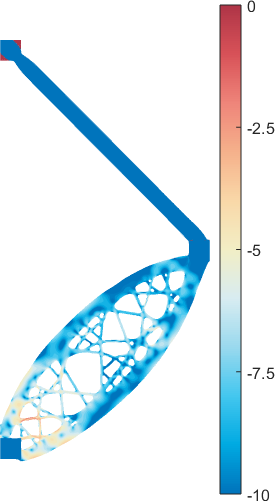}} \quad
  \subfloat[$\mathcal{N}_{e}$]{
   \includegraphics[scale = 0.275, keepaspectratio]
   {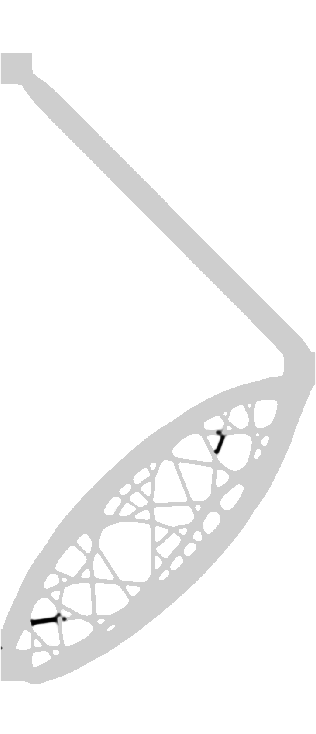}} \quad
  \subfloat[$\boldsymbol{\varphi}_{6}$]{
   \includegraphics[scale = 0.225, keepaspectratio]
   {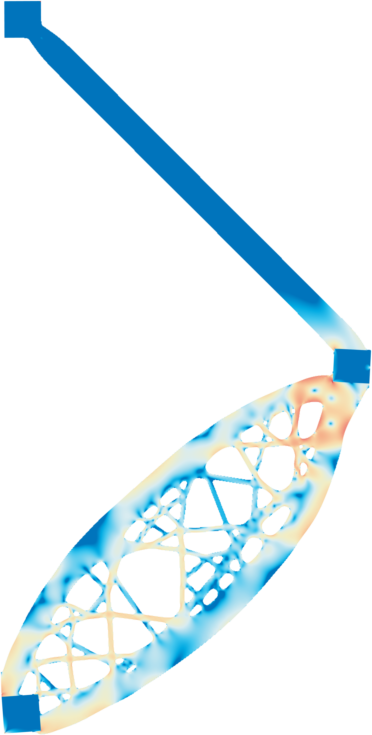}} \quad
  \subfloat[$\boldsymbol{\varphi}_{10}$]{
   \includegraphics[scale = 0.225, keepaspectratio]
   {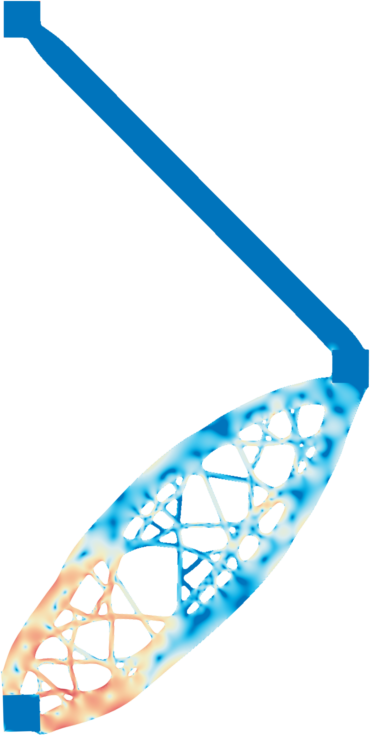}}
 \caption{Illustration of the local thickening post--processing operation. (a,b) show the log--scale variation of $|\nabla \pi|$ associated with the two modes representing a single--bar failure for the design of \autoref{fig:NumericalTest-2DExample} (c). The set $\mathcal{N}_{e}$ corresponding to the threshhold $|\nabla \pi_{j^{\ast}|e}| \geq 10^{-2} \max_{e} \nabla \pi_{j^{\ast}|e}$ is shown in (c) and (d,e) show the two modes for the reinforced structure}
 \label{fig:postProcessingFilterwithTV}
\end{figure}

\subsection{Discussion on localized buckling modes}
 \label{sSecSpuriousLocalizedModes}

If the interpolations for $E_{\kappa}$ and $E_{\sigma}$ are not such that the Rayleigh quotient \eqref{eq:minBLFcharRayleigh} remains bounded as $x_{e} \rightarrow 0$, artificial buckling modes are likely to happen in low density regions \cite{neves-etal_95a}. With the choice of \eqref{eq:interpolationSIMP2-SIMP}, originally proposed by \cite{book:bendsoe-sigmund_2004}, we have not found this issue. According to the criterion proposed in \cite{gao-ma_15a}, based on the ratio between the strain energy density associated with low density regions and the overall one, some high order modes occasionally become spurious as the optimization proceeds. However, the evolution of the fundamental BLF is generally smooth, and jumps are associated with the increase of the penalization factor (see \autoref{fig:NumericalTest-2DExample-convergenceCurves-M4}).

The localized buckling modes discussed in the previous section have nothing to do with grayscales, as they appear in solid regions of the pure BW design. With the goal of obtaining a preliminary design, dealing with these localized modes in the optimization process does not seem a good strategy, at least for the following two reasons

\begin{enumerate}
 \item Extremely localized deformations, such as $\boldsymbol{\varphi}_{3}$, $\boldsymbol{\varphi}_{5}$ and $\boldsymbol{\varphi}_{6}$ in \autoref{fig:NumericalTest-2DExample-modesComparison-M4} or $\boldsymbol{\varphi}_{2}$, $\boldsymbol{\varphi}_{3}$ and $\boldsymbol{\varphi}_{5}$ in \autoref{fig:NumericalTest-2DExample-modesComparison-M8}, are artifacts due to stress concentrations and/or singularities \cite{book:achenbach1973}, linked to sharp geometric variations, such as corners or sharp boundaries (see \autoref{fig:isotropicMeasureStress}(a)). In reality, the critical condition in those regions will be the reaching of a limit stress and material failure, before geometric instability. We underline that such stress concentrations are by no means a prerogative of density--based TO. Other researchers, using alternative parameterizations (e.g. level--set \cite{dunning-etal_16a}), have experienced similar artifacts;
 \item As an intrinsic trend, compliance or mass--optimized designs may show many thin bars, especially for fine discretizations and/or low volume fractions while building up a hierarchical structure with extreme buckling response. Thus, many local modes involving single bars may appear, and taking all of them into account would overly increase the number of eigenpairs to be computed.
\end{enumerate}

Therefore, we on purpose overlook local modes in the optimization process, achieving a design with ``global" stability more easily. An exhaustive discussion of local buckling \cite{book:bazant-cedolin2010, rozvany_96a} is beyond the scope of this work. However, we point out that a global, ``averaged" approach has since long been recognized as being meaningful for studying the geometric stability of continua \cite{book:achenbach1973}, as local effects are not soundly defined.

Concerning Point 1, we have already discussed how local artifacts can be alleviated by the multilevel approximation strategy when $\Omega_{\ell}$ is set coarse enough. This is due to the intrinsic filtering effect of the Galerkin projection, which applied to the fine--grid stress stiffness matrix $G[\mathbf{x}, \mathbf{u}(\mathbf{x})]$ essentially corresponds to a smoothing of stresses. This is visualized in \autoref{fig:isotropicMeasureStress}, showing the following scalar representation of the stress stiffness matrix
\begin{equation}
 \label{eq:isomeasG}
  \mathcal{G}^{\ell}_{i} = \sqrt{(G^{\ell}_{2i})^{2} + (G^{\ell}_{2i-1})^{2}}
\end{equation}
on three levels $\ell$. Labeling $G^{\ell}_{k} = G^{\ell}_{kk}$, the diagonal coefficients of the stress stiffness matrix corresponding to the $\ell$--th level, \eqref{eq:isomeasG} associates to each node $i$ of the discretization an equivalent value accounting for the contribution of the two DOFs. On $\Omega_{1}$, this quantity appears very localized in regions where localized modes appear (see \autoref{fig:NumericalTest-2DExample-modesComparison-M4}), wherease it becomes progressively more spread on coarser levels. Moreover, the parameter $\zeta_{\ell} = \max_{i} \mathcal{G}^{\ell}_{i}/\mathcal{K}^{\ell}_{i}$, where $\mathcal{K}^{\ell}_{i}$ is defined as in \eqref{eq:isomeasG} corresponding to the stiffness matrix, descreases while shifting on coarser grids, attesting the decrease of local effects.

The optimized design may be reinforced against buckling of single bars in a post--processing step, marginally affecting the structural volume. The procedure may build on first identifying such local buckling by computing the Strain Energy Density ($\pi_{j|e} = \boldsymbol{\varphi}^{T}_{j|e}K_{e}[x_{e}]\boldsymbol{\varphi}_{j|e}$, $e = 1,\ldots, m$) for each mode, and the Total Variation of this quantity \cite{jordan_881a}

\begin{equation}
 \label{eq:TVsed}
 \mathfrak{L}_{j} = \mathfrak{L}\left( \pi_{j}, \Omega_{\geq} \right) = 
 \int\limits_{\Omega_{\geq}}
 |\nabla \pi_{j}| \: {\rm d}\Omega
\end{equation}
where $\nabla(\cdot)$ here is the spatial gradient and $\Omega_{\geq} := \{ \Omega^{e} \in \Omega_{1} \, \mid \: x_{e} \geq \bar{x} \}$ identifies the solid domain on the fine discretization. We may choose $\bar{x} = 1$ for a pure BW design. If $\mathfrak{L}_{j} = 0$, the quantity $\pi_{j}$ is constant on $\Omega_{\geq}$, and the corresponding $\boldsymbol{\varphi}_{j}$ is a global mode. On the other hand $\boldsymbol{\varphi}_{j}$ becomes more and more localized as $\mathfrak{L}_{j} \rightarrow \infty$.

Once identified the localized modes, say $j^{\ast} \in\mathcal{B}$, a local dilation--like operator \cite{sigmund_07a} is introduced

\begin{equation}
 \label{eq:postProcessingDilation}
  \mathcal{R}[ x_{e} ] = \max_{k\in\mathcal{N}_{e}} x_{k}
\end{equation}
where $\mathcal{N}_{e}$ is the set of neighboring elements $\Omega^{k}$ s.t. ${\rm dist}(\Omega^{e}, \Omega^{k}) \leq r_{\rm th}$ and $|\nabla \pi_{j^{\ast}|e}|$ exceeds a fixed value $\forall j^{\ast}\in\mathcal{B}$. \autoref{eq:postProcessingDilation} is used to perform a local thickening of the critical bars, according to the width $r_{\rm th}$.

For example, for the structure of \autoref{fig:NumericalTest-2DExample} (c), modes $\boldsymbol{\varphi}_{6}$ and $\boldsymbol{\varphi}_{10}$ each involve the buckling of a single thin bar and are associated with values $\mathfrak{L}_{6} = 41.7$ and $\mathfrak{L}_{12} = 13.1$, which are very high if normalized w.r.t to the one of the fundamental buckling mode $\mathfrak{L}_{1} = 1.0$. Based on the distribution of $| \nabla \pi_{j^{\ast}} |$ for $j^{\ast} = 6, 12$, the set $\mathcal{N}_{e}$ is built (here considering elements where $|\nabla \pi_{j^{\ast}|e}| \geq 10^{-2} \max_{e} \nabla \pi_{j^{\ast}|e}$) and the operator $\mathcal{R}$ performs the local thickening as shown in \autoref{fig:postProcessingFilterwithTV} (d). In the end, these local modes are avoided within the set $\mathcal{B}$, with a minute increase of about $0.36\%$ in the structural volume. The variations of the BLFs are reported in the last column of \autoref{tab:tableB}, and we notice that all of them increase a few percent, with those corresponding to global modes barely affected by this local reinforcement.

%% GEOMETRICAL SKETCH OF THE 3D CANTILEVER BEAM WITH THE TWO SECTIONS SHOWING THE LONGITUDINAL AND TRANSVERSAL VIEWS
\begin{figure}[tb]
 \centering
  \includegraphics[scale = 0.225, keepaspectratio]
   {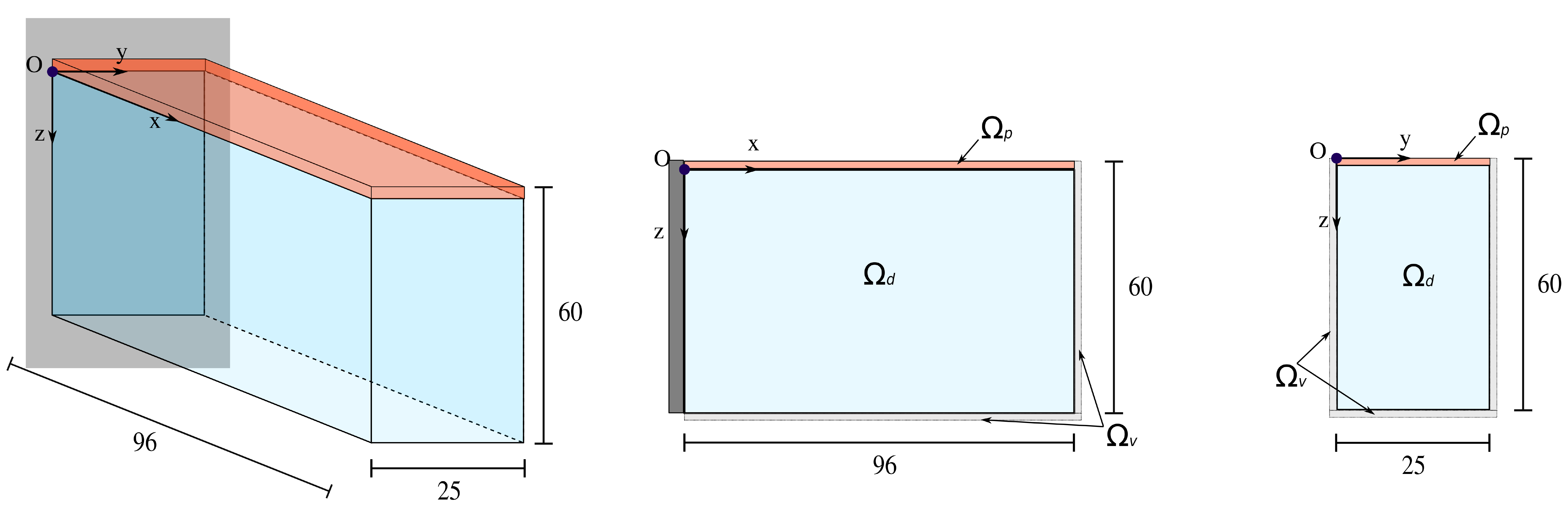}
 \caption{Setting for the 3D cantilever example. The beam is fully clamped at the built in end and a uniform load with total magnitude $\| \mathbf{f} \| = 3.6\cdot 10^{4}$ acts on the top face. A layer of passive, solid elements with thickness $t = e_{z}/30$, is at the top face (see $\Omega_{p}$ shown in red), while the shaded area denoted by $\Omega_{v}$ is the extended domain used for padding the density filter (see \cite{clausen-andreassen_17a} for details)}
 \label{fig:test3Dcantilever}
\end{figure}

%% COMPLIANCE DESIGN AND FUNDAMENTAL BUCKLING MODE
\begin{figure}[tb]
 \centering
  \subfloat[]{
   \includegraphics[scale = 0.065, keepaspectratio]
   {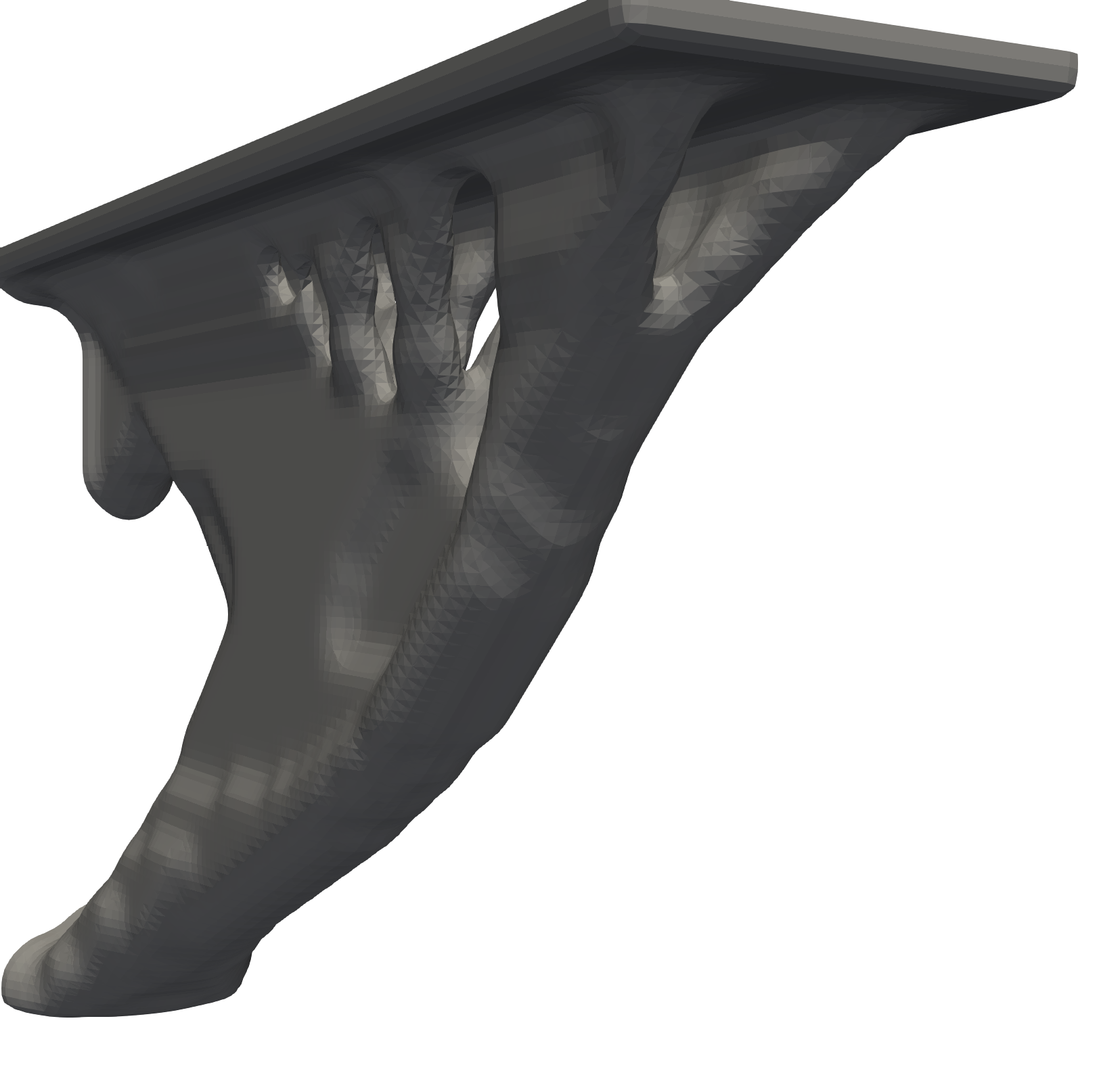}}
  \subfloat[]{
   \includegraphics[scale = 0.065, keepaspectratio]
   {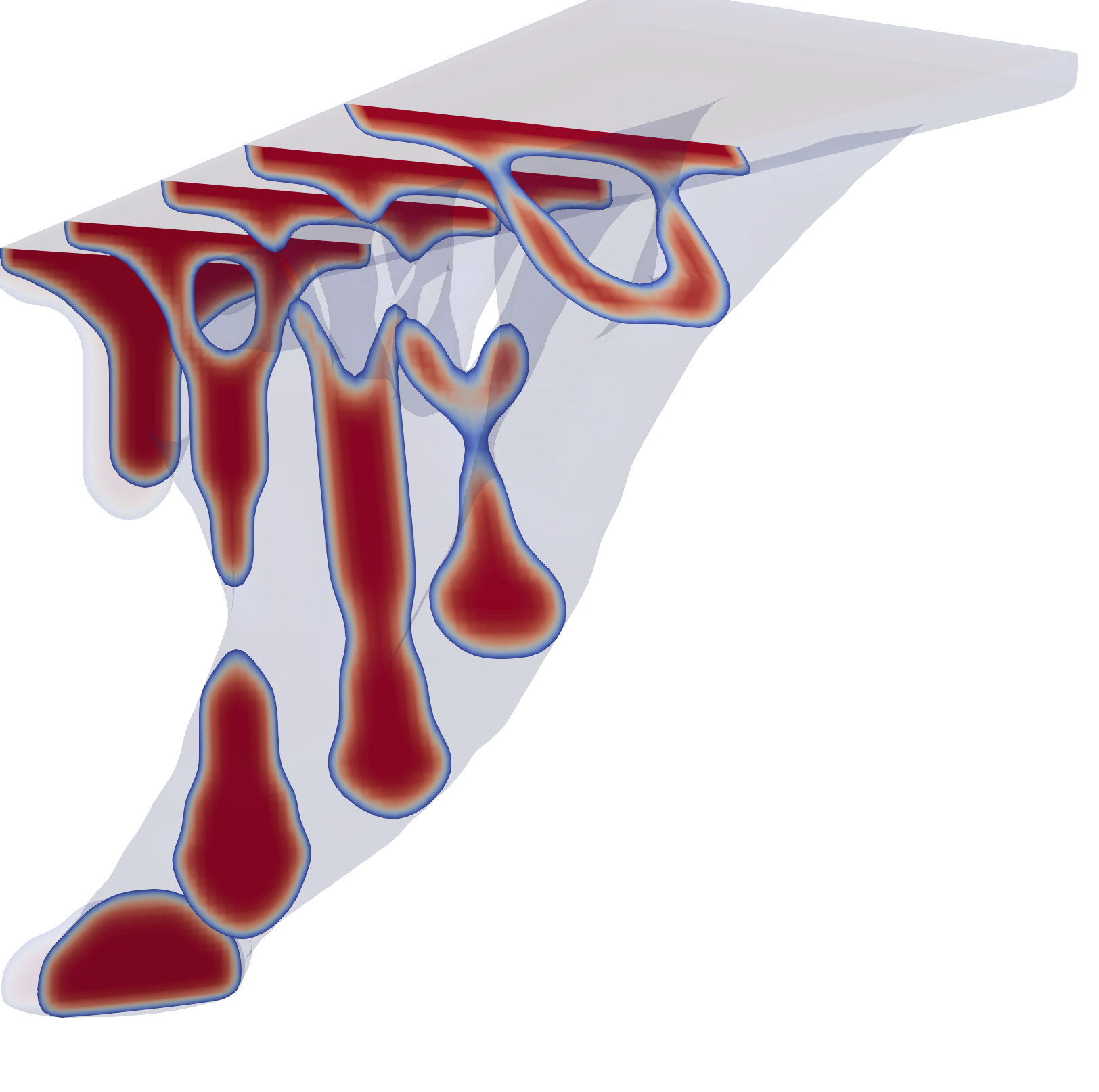}}
  \subfloat[$\boldsymbol{\varphi}_{1}$]{
   \includegraphics[scale = 0.065, keepaspectratio]
   {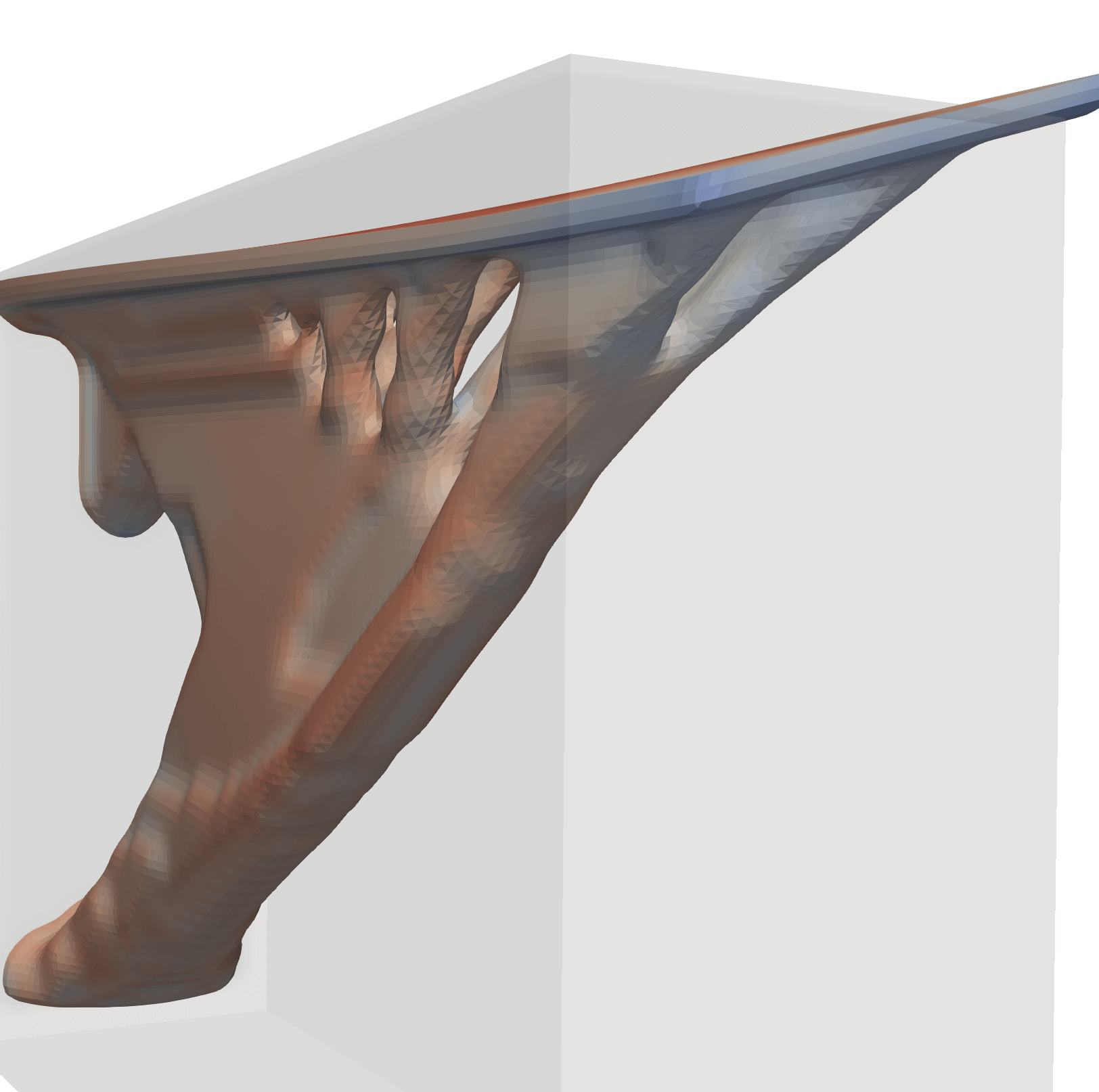}}
  \subfloat[$\boldsymbol{\varphi}_{2}$]{
   \includegraphics[scale = 0.065, keepaspectratio]
   {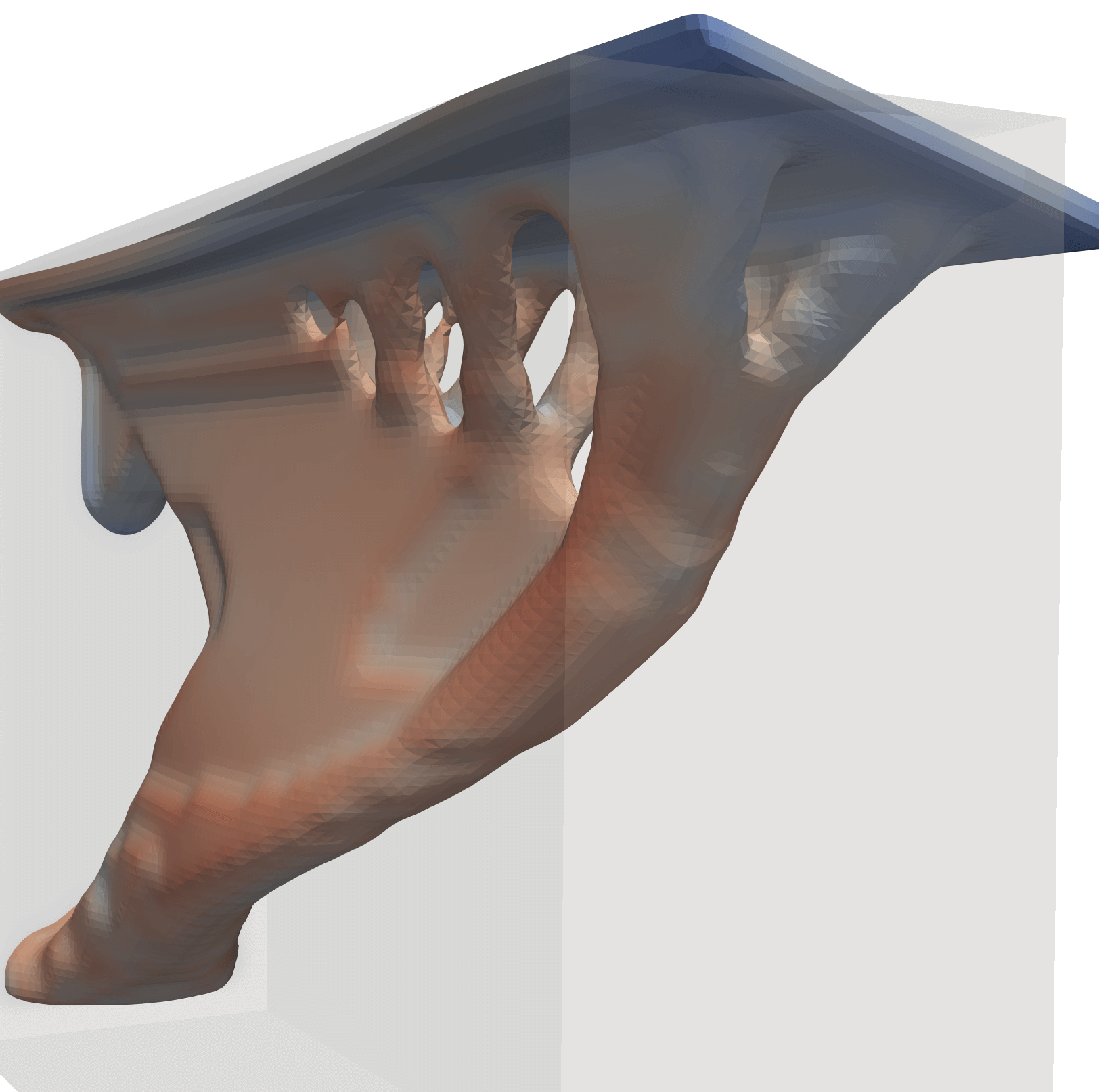}}
 \caption{Compliance design obtained from \eqref{eq:optProblemMinMass3D} for $\bar{\lambda} = 0$ (a). The value of the volume fraction is $f = 0.1295$. In (b) we can see some cross sections illustrating the concentration of material in a single strut at the center and (c,d) show the lowest two buckling modes, and the distribution of strain energy associated with these deformations}
 \label{fig:3DcantileverComplianceDesign}
\end{figure}

% DESIGNS CORRESPONDING TO SOME PRESCRIBED BUCKLING LOAD FACTORS
\begin{figure}[tb]
 \centering
  \subfloat[$\bar{\lambda} = 0.25$, $f = 0.135$]{
   \includegraphics[scale = 0.065, keepaspectratio]
   {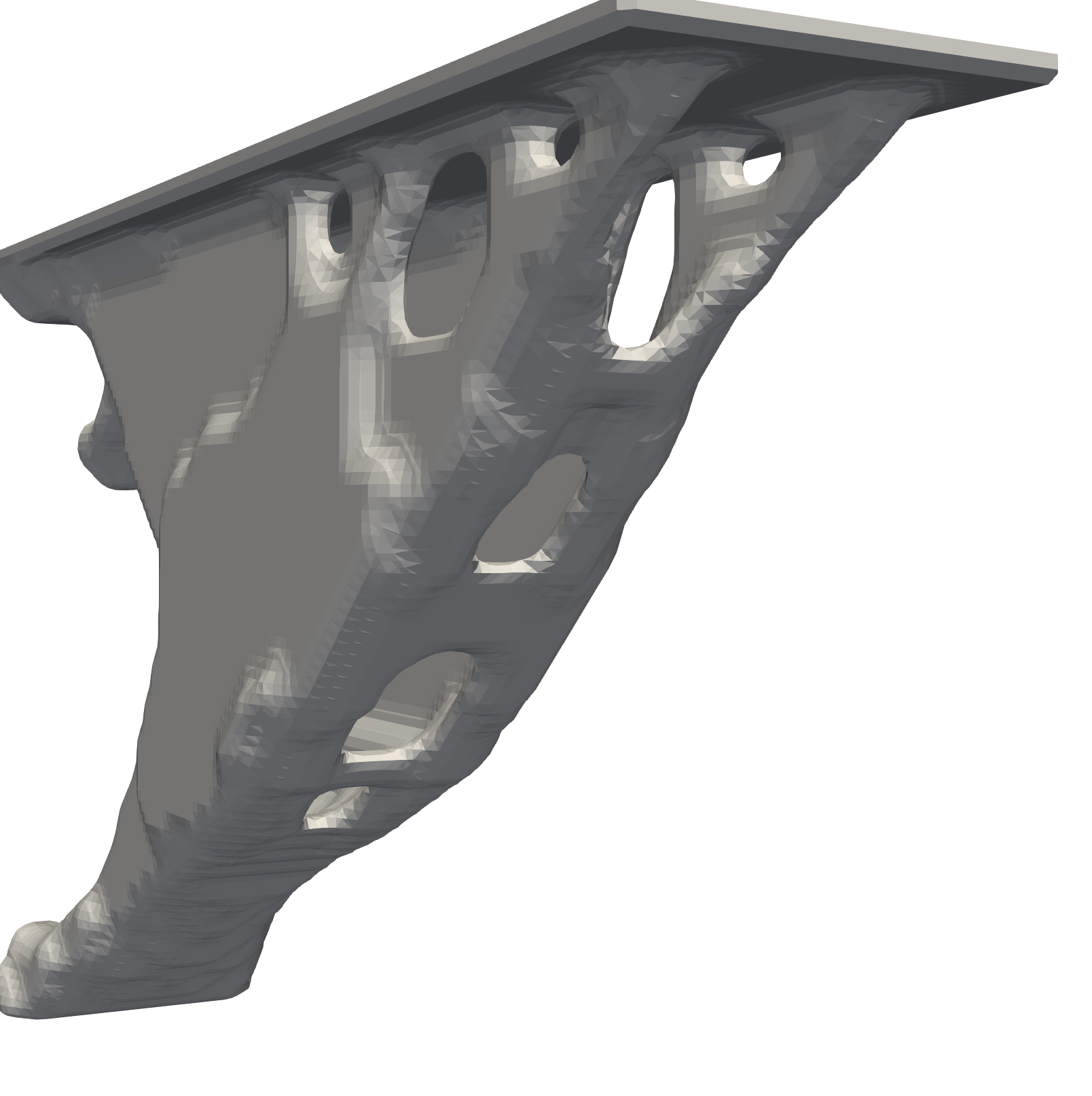}}
  \subfloat[$\bar{\lambda} = 0.5$, $f = 0.145$]{
   \includegraphics[scale = 0.065, keepaspectratio]
   {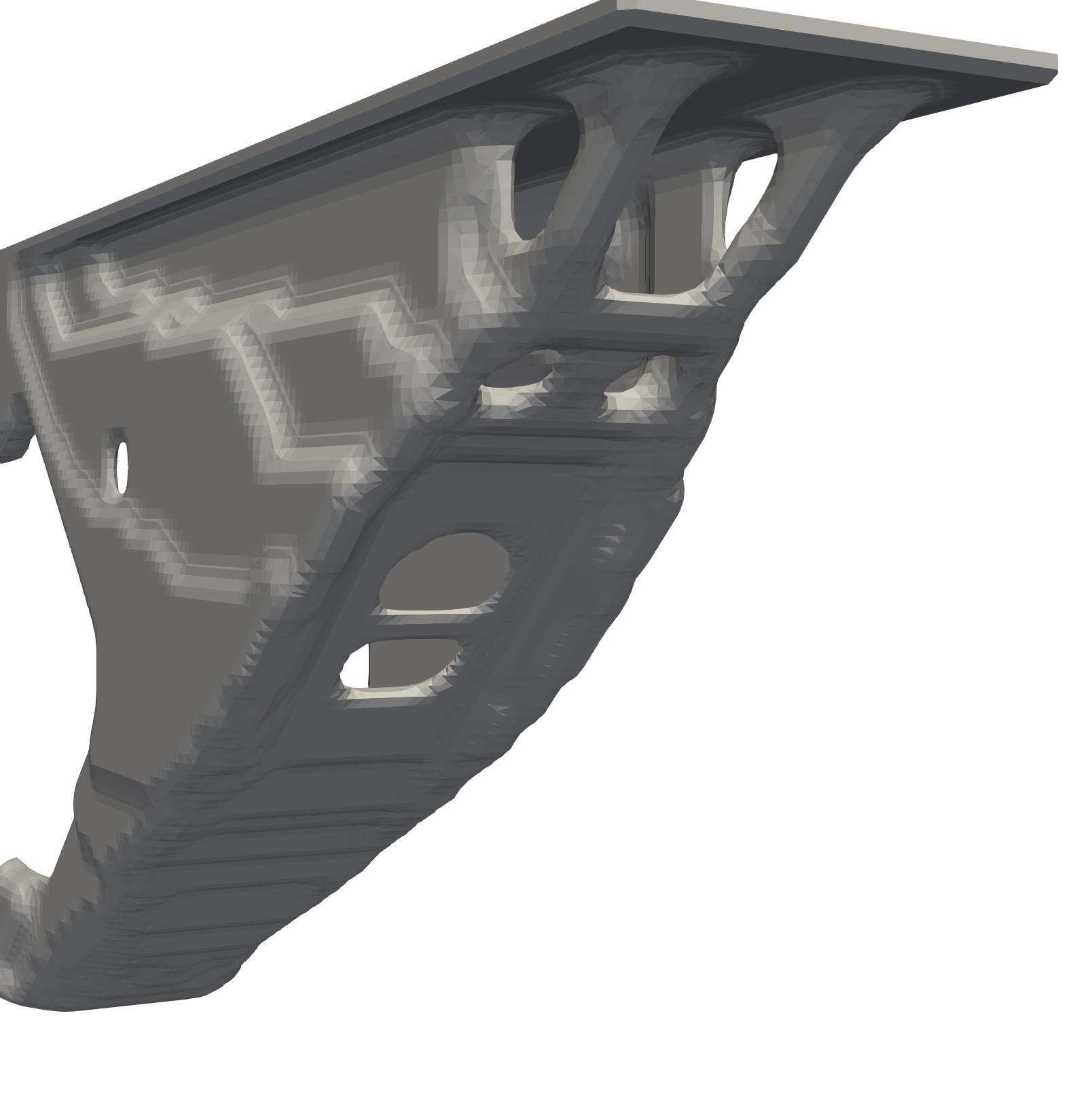}}
  \subfloat[$\bar{\lambda} = 0.75$, $f = 0.172$]{
   \includegraphics[scale = 0.065, keepaspectratio]
   {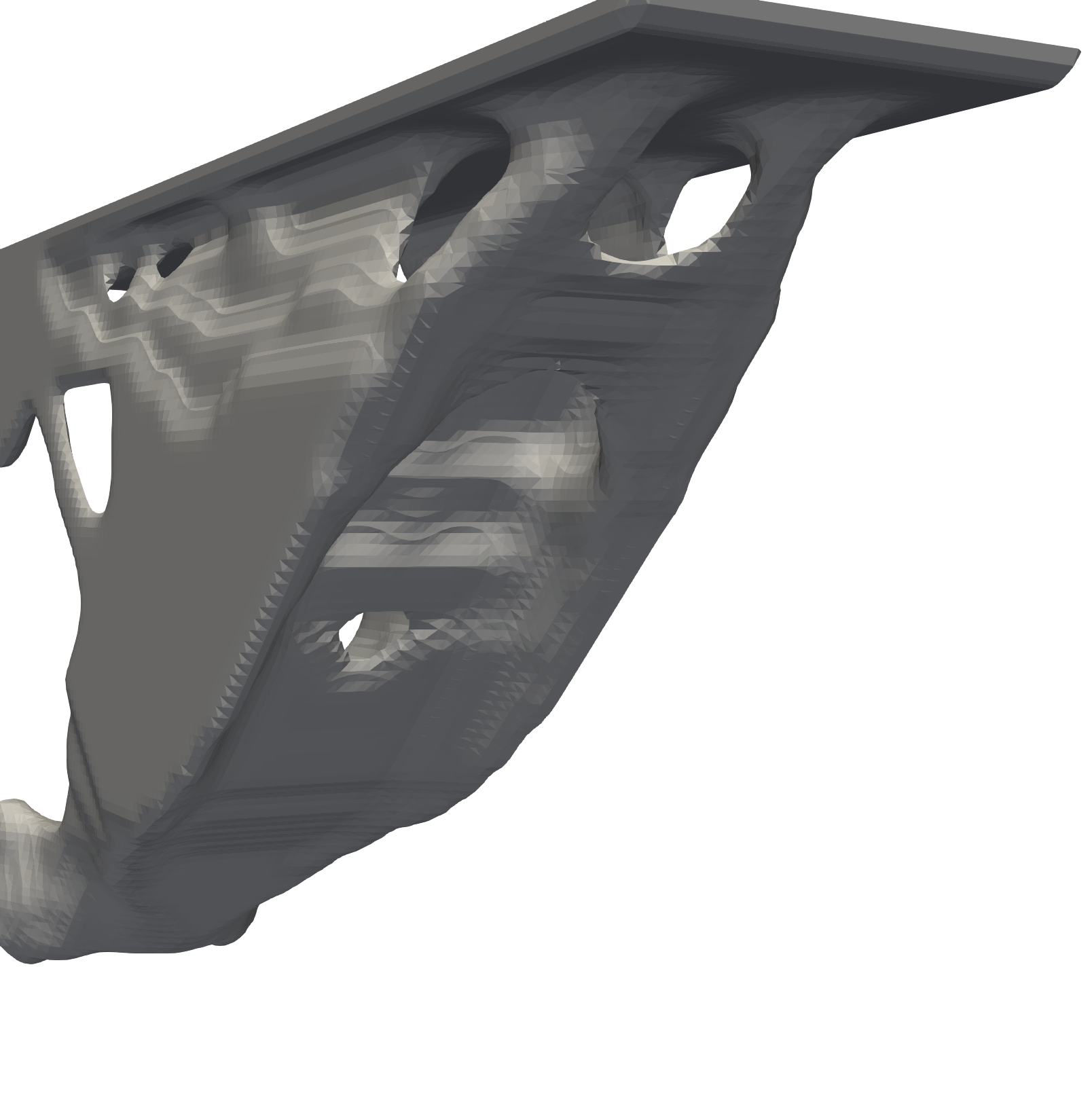}}
  \subfloat[$\bar{\lambda} = 1.0$, $f = 0.221$]{
   \includegraphics[scale = 0.065, keepaspectratio]
   {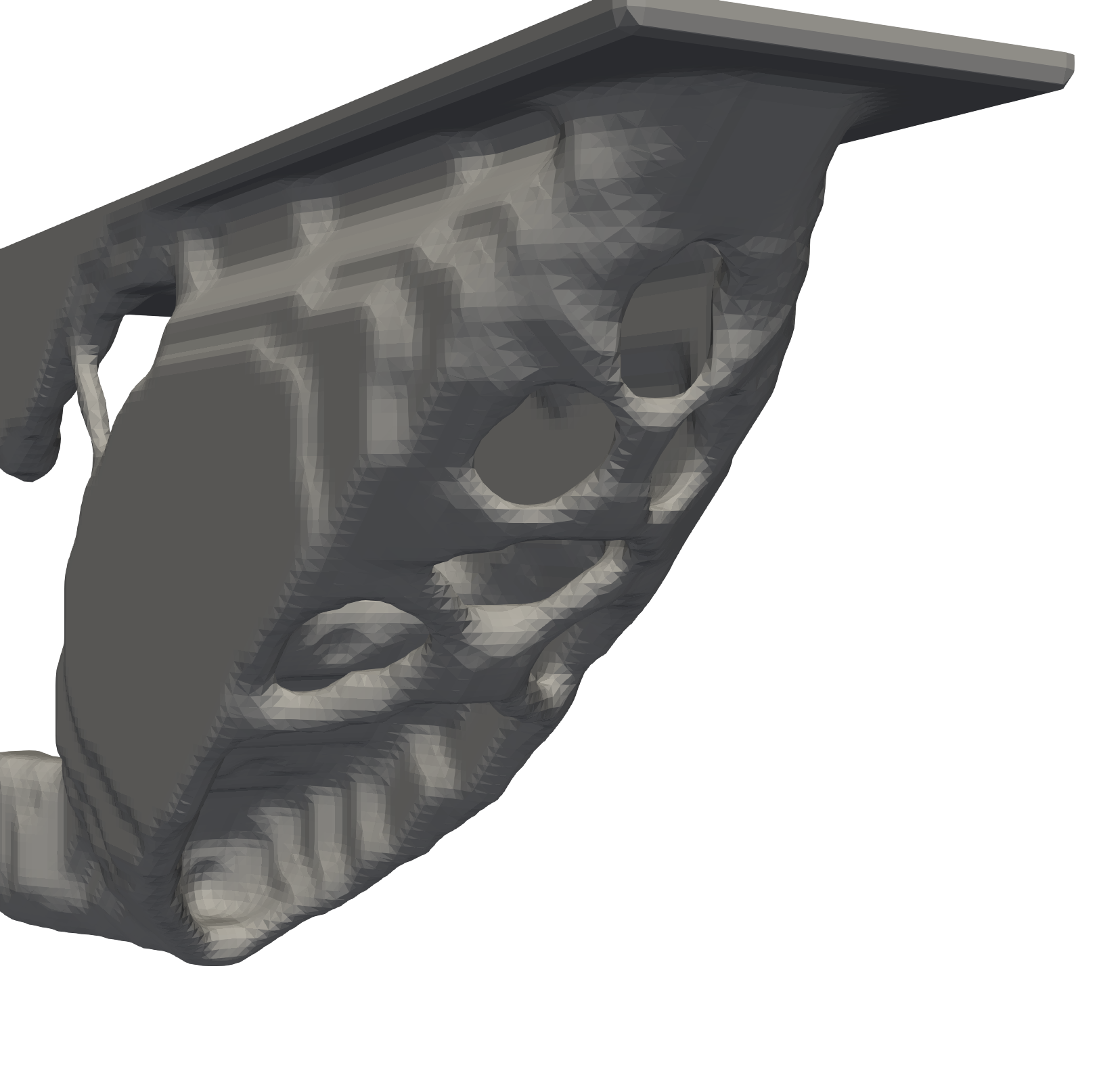}} \\
  \subfloat[$\bar{\lambda} = 0.25$]{
   \includegraphics[scale = 0.065, keepaspectratio]
   {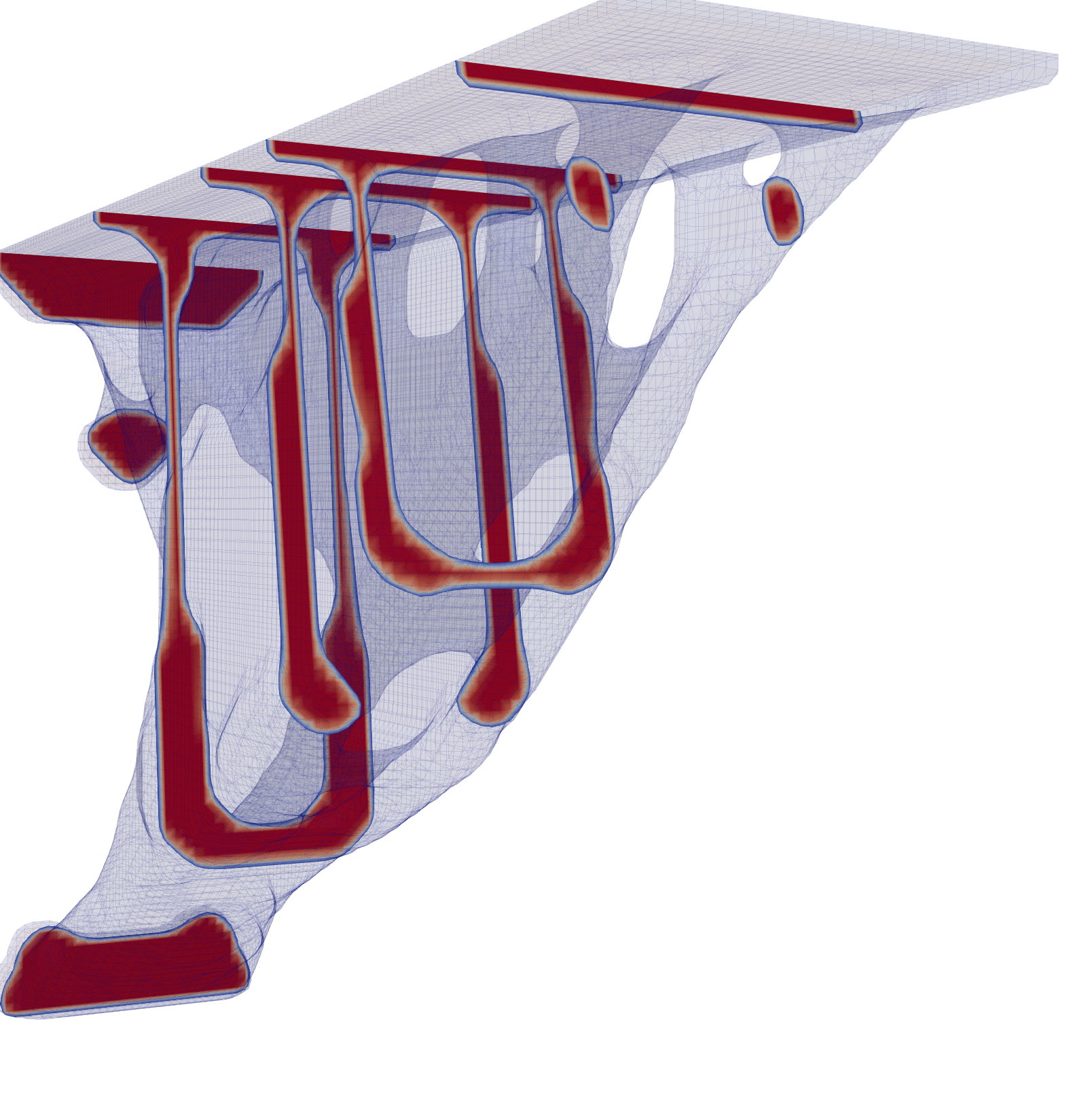}}
  \subfloat[$\bar{\lambda} = 0.5$]{
   \includegraphics[scale = 0.065, keepaspectratio]
   {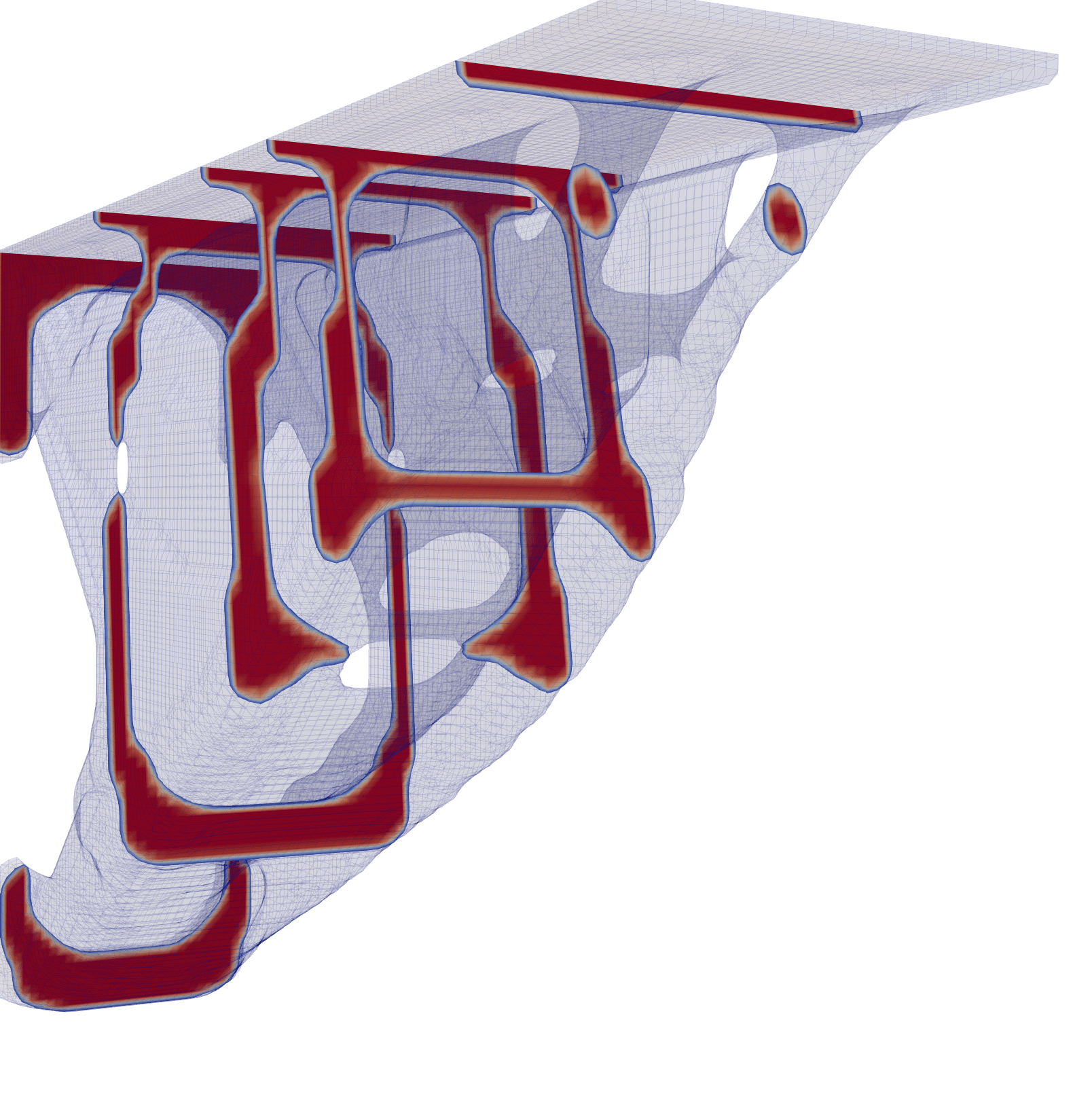}}
  \subfloat[$\bar{\lambda} = 0.75$]{
   \includegraphics[scale = 0.065, keepaspectratio]
   {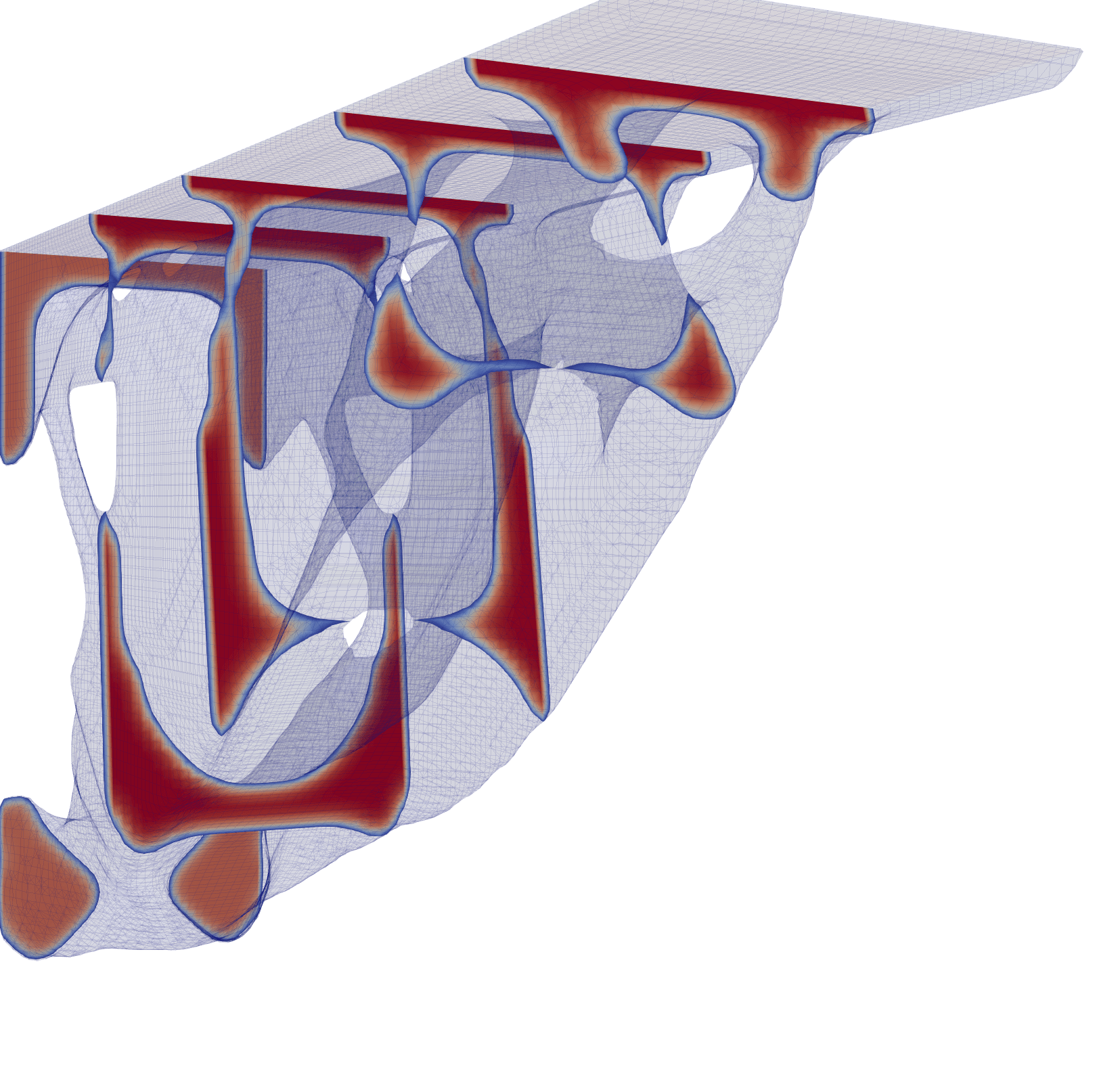}}
  \subfloat[$\bar{\lambda} = 1.0$]{
   \includegraphics[scale = 0.065, keepaspectratio]
   {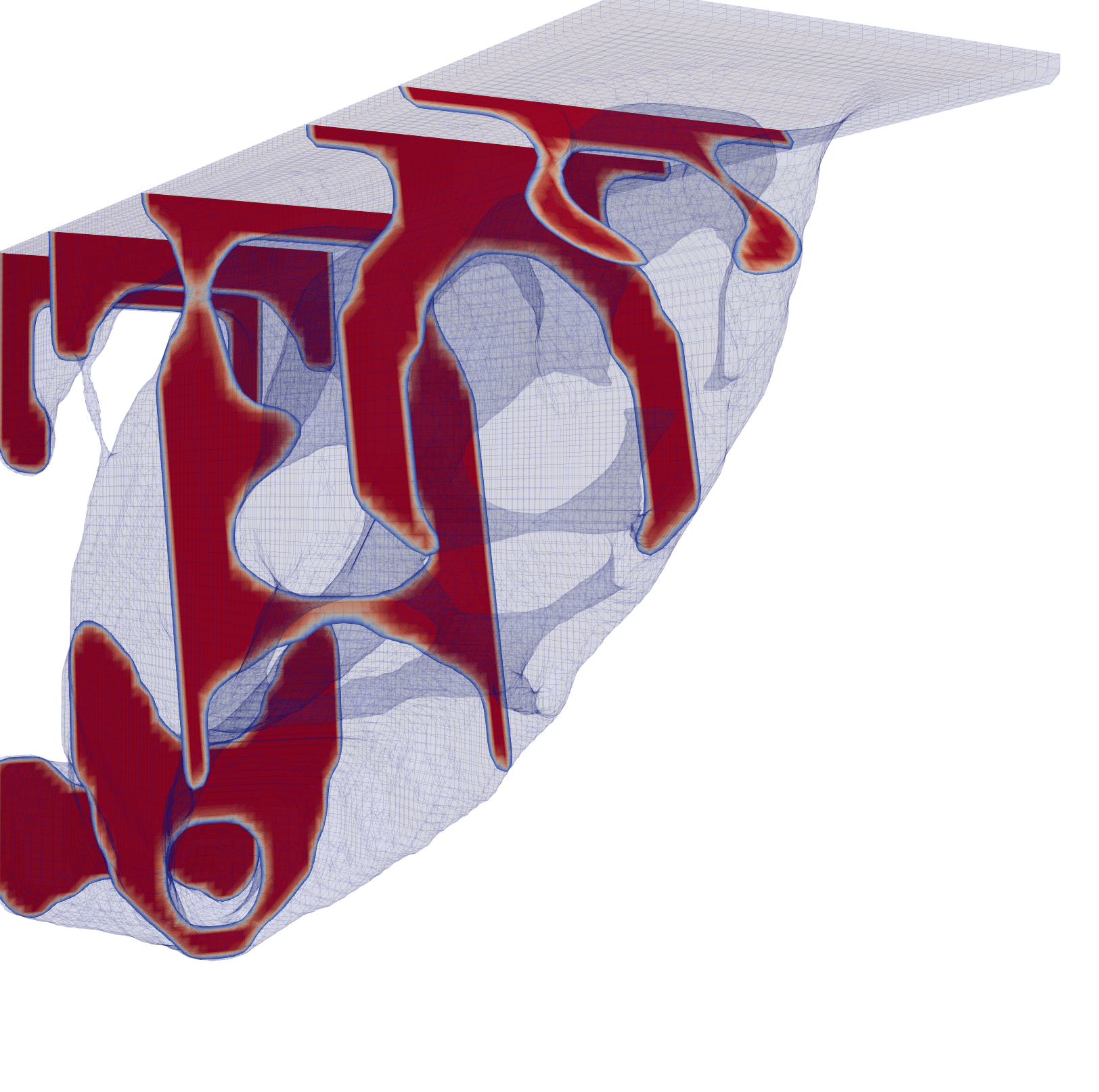}}
 \caption{Some designs obtained for increasing values of the lower bound on the fundamental BLF (a--d) and corresponding views with sections (e--h). Physical densities $x_{e} \geq 0.9$ are plotted}
 \label{fig:3DcantileverOptTop}
\end{figure}

% SECTION VIEW OF THE DESIGN CORRESPONDING TO BLF = 1 AND SOME BUCKLING MODES
\begin{figure}[tb]
 \centering
  \includegraphics[scale = 0.165, keepaspectratio]
   {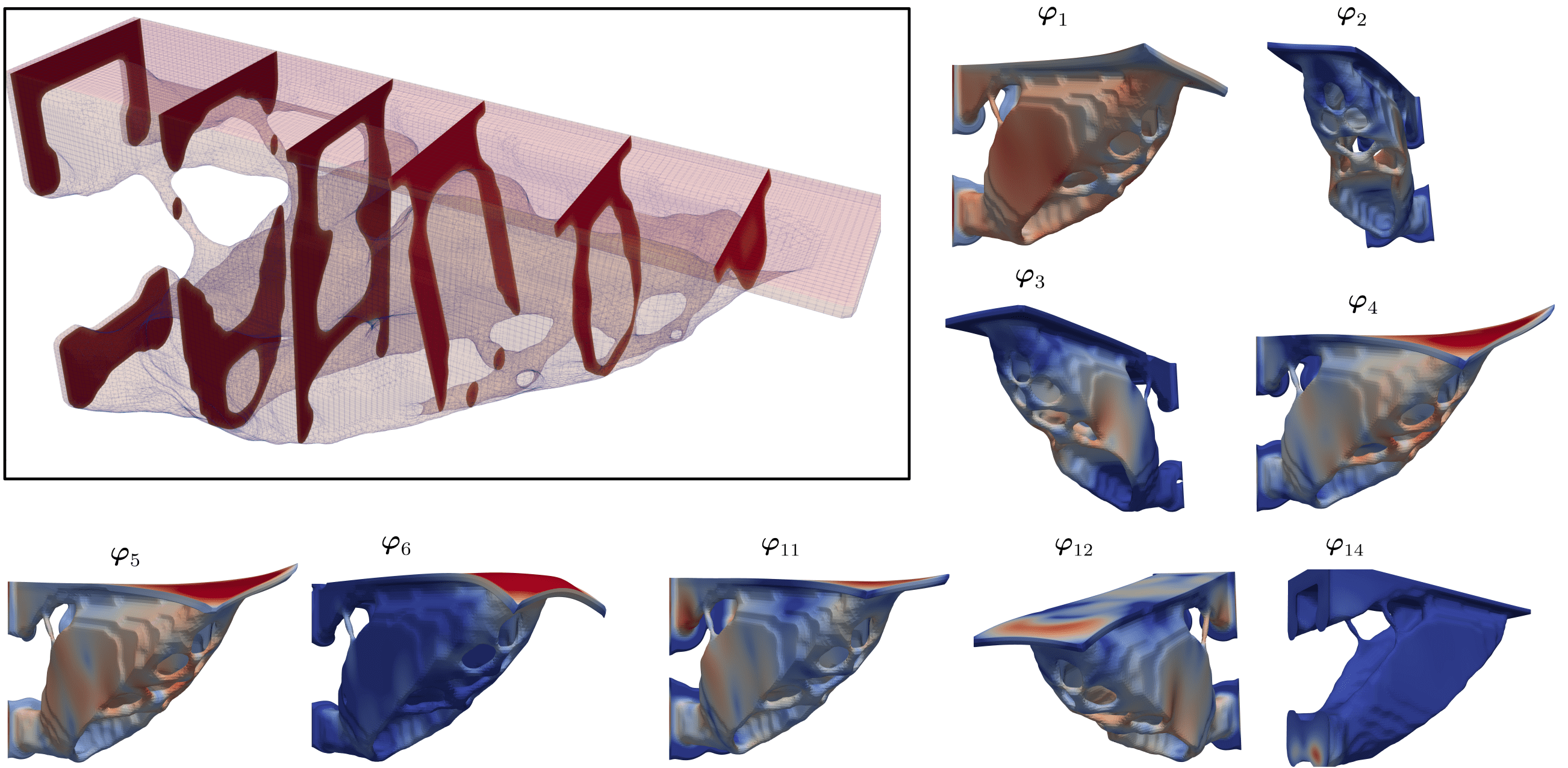}
 \caption{Details of the 3D design obtained for $\bar{\lambda}=1$ (left) and some of the buckling modes. The colormap (blue to red) on buckling modes refers to the logarithm of the strain energy density (low to high)}
 \label{fig:3DcantileverOptTopPc1}
\end{figure}

%% PLOT SUMMARIZING THE OPTIMIZATION HISTORY AND THE RELATIONSHIP BETWEEN MASS, BLFs AND THEIR COALESCENCE
\begin{figure}[tb]
 \centering
  \subfloat[]{
   \includegraphics[scale = 0.35, keepaspectratio]
   {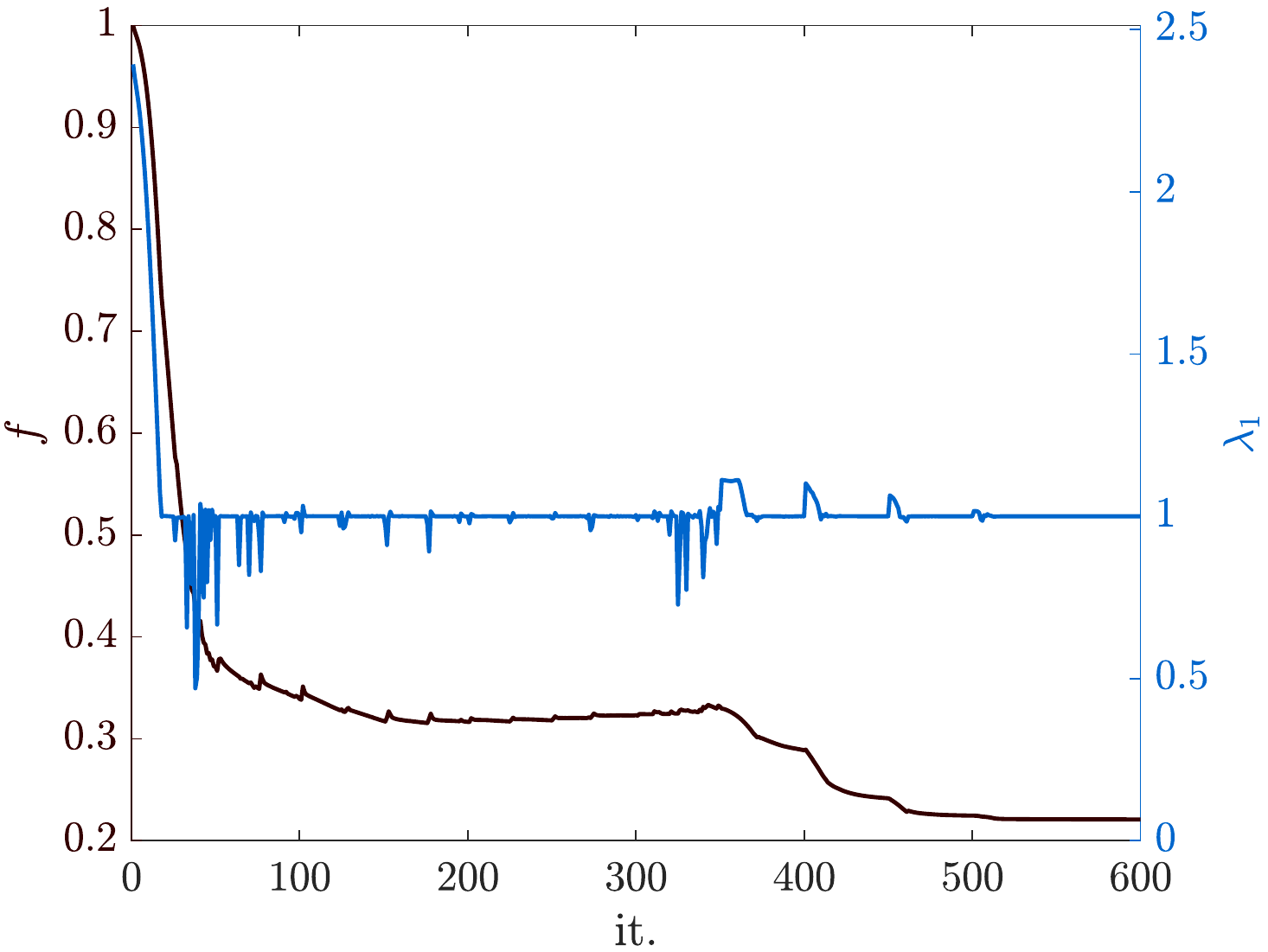}}
  \subfloat[]{
   \includegraphics[scale = 0.35, keepaspectratio]
   {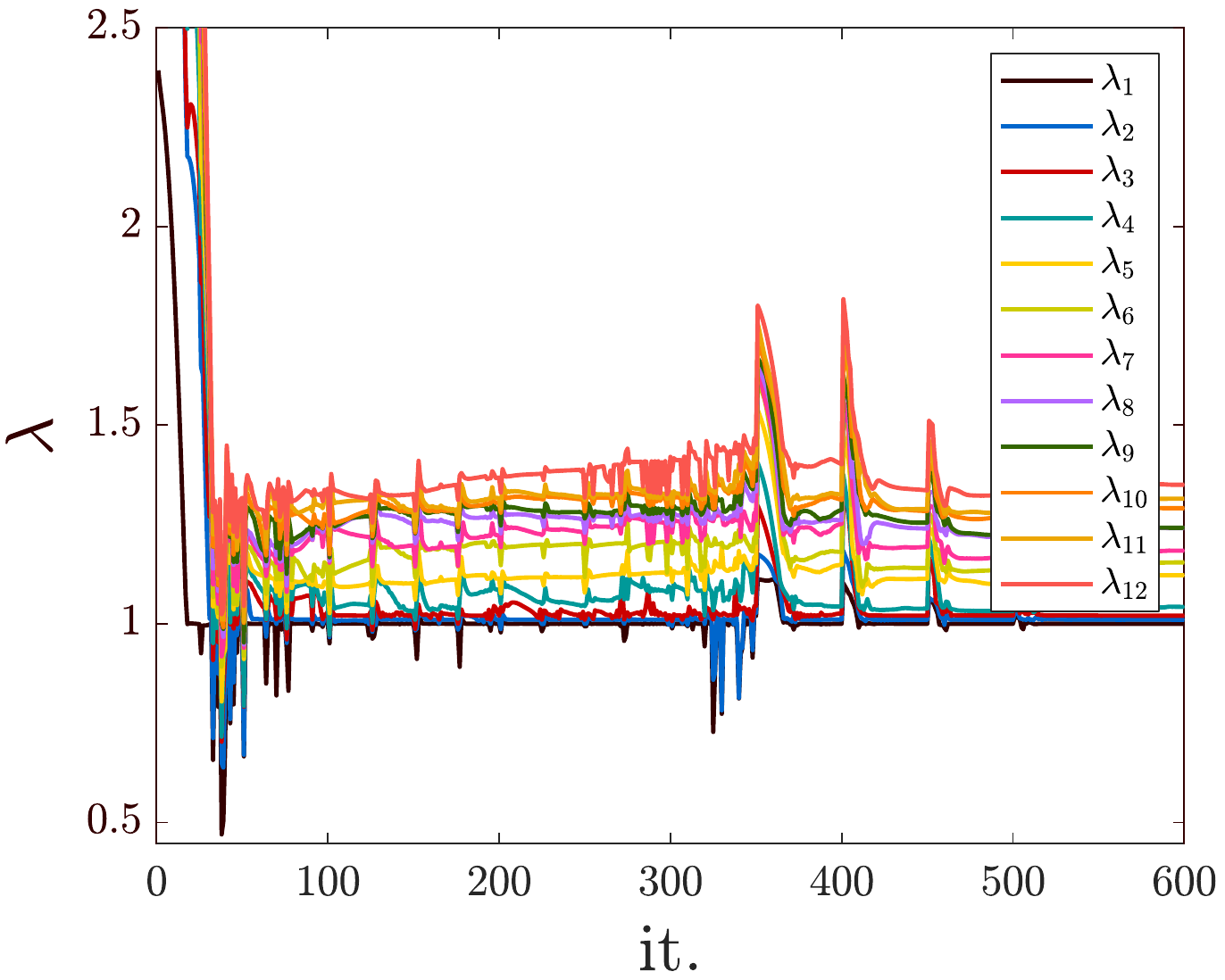}}
  \subfloat[]{
   \includegraphics[scale = 0.35, keepaspectratio]
   {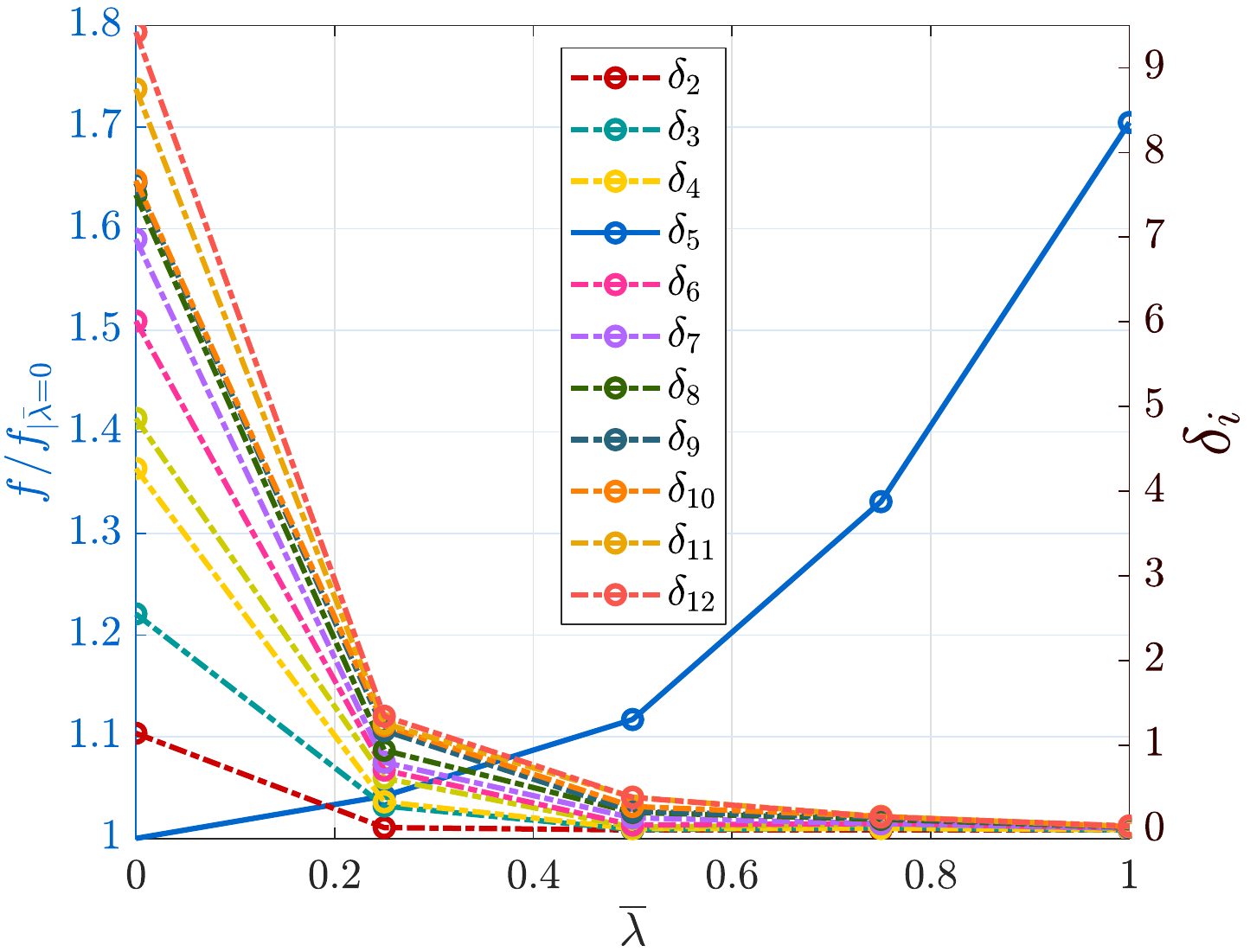}}
 \caption{(a) shows the evolution of the volume fraction (black curve) and of the fundamental BLF (blue curve) and (b) the evolution of the BLFs corresponding to the constrained modes for problem \eqref{eq:optProblemMinMass3D} with $\bar{\lambda} = 1$. (c) shows the relationships between $\bar{\lambda}$ and the volume fraction of the obtained designs, normalized with respect to that of the design corresponding to $\bar{\lambda} = 0$ (blue curve, plotted against the left axis), and the eigenvalue separation parameters $\delta_{i}$ (plotted against the right axis)}
 \label{fig:3DcantileverConvCurves}
\end{figure}

\section{Mass minimization of a 3D cantilever beam}
 \label{Sec:3DcantileverMassMinimization}

We now consider the mass minimization problem for the 3D structure sketched in \autoref{fig:test3Dcantilever}, inspired by the example presented in \cite{dunning-etal_16a}

\begin{equation}
 \label{eq:optProblemMinMass3D}
 \mathcal{P}_{2}
 \begin{cases}
  & \min\limits_{\hat{\mathbf{x}}\in [ 0, 1 ]^{m}}
  f = V\left(\mathbf{x}\right) / | \Omega_{h} | \\
  {\rm s.t.} & \min\limits_{i\in\mathcal{B}}\lambda_{i}
  \geq \bar{\lambda} \\
  & \mathbf{u}^{T}K\left[\mathbf{x}\right] \mathbf{u} \leq 6\bar{J}(\mathbf{x} = 1)
 \end{cases}
\end{equation}
for some values of $\bar{\lambda} \in [0, 1]$. The buckling load factor of the fully solid design ($f = 1$) is $\lambda_{1} = 2.393$; therefore we are always starting within the feasible set.

The fine scale discretization is set to $\Omega_{1} = \left( 134 + 2, 44 + 4, 94 + 2 \right)$, corresponding to 626,688 design variables and 1,953,483 DOFs, about 4.5 times more than in \cite{dunning-etal_16a}. The multilevel procedure is built with $\ell = 3$ and on $\Omega_{\ell}$ we just have 34,125 DOFs. The density filter radius is set to $r_{\rm min} = 4\sqrt{3}$ and, in order to alleviate boundary effects, we adopt the domain extension strategy suggested in \cite{clausen-andreassen_17a}. Therefore, in the above we referred to the fine discretization with the convention $\Omega_{1} = \left( e_{x} + {\rm d}_{x}, e_{y} + {\rm d}_{y}, e_{z} + {\rm d}_{z} \right)$, where ${\rm d}_{j}$ is the number of elements extending the domain for each direction. For this example, 24 modes are computed and the lowest 12 BLFs are constrained.

The optimization is run for a total of 600 steps and, in order to obtain efficient designs with low grayscale, the following continuation strategy is used for $p$ and $\beta$. The optimization is started with $p = 3$ and $\beta = 1$, then $p$ is raised with increments $\Delta p = 0.25$ each 25 steps, up to the value $p_{\rm max} = 6$. As this penalization value is reached, $\beta$ is doubled each 50 steps up to the value $\beta_{\rm max} = 32$.

\autoref{fig:3DcantileverComplianceDesign} (a) shows the optimized design when only considering the compliance constraint ($\bar{\lambda} = 0$) and we see how the material is mainly localized near the centerline of the structure. Near the tip, the top face is supported by some thin members, converging in a single strut connecting with the built-in end at the foot. It can be clearly recognized how the strut is taking an I-shaped configuration while progressing towards the built-in end, then splitting into two regions, optimizing the bending response. The volume fraction for this design is $f = 0.1295$, and the compliance constraint is active. The fundamental BLF is $\lambda_{1} = 0.067$; thus, from a practical point of view the design is worthless, as it would buckle under the external load. The fundamental buckling mode, shown in \autoref{fig:3DcantileverComplianceDesign} (c) resembles a global twist of the structure in the $y-z$ plane, which is poorly restrained by the supporting structure.

Intuitively, material should be deployed far from the centerline, to resist rotations in the $y-z$ plane when including the buckling constraint. Designs obtained for higher values of the lower bound $\bar{\lambda}$ are displayed in \autoref{fig:3DcantileverOptTop}. As expected, the material is progressively moved away from the centerline and for $\bar{\lambda} \geq 0.5$ two distinct shear plates appear. An internal distribution of reinforcing beams is also emerging as $\bar{\lambda}$ is increased further (see \autoref{fig:3DcantileverOptTop}(c)).

A detailed view of the design corresponding to $\bar{\lambda} = 1$ is given in \autoref{fig:3DcantileverOptTopPc1}, where some of the buckling modes are also displayed. All the constrained modes $\boldsymbol{\varphi}_{1}-\boldsymbol{\varphi}_{12}$ represent global, physically meaningful deformations, either involving a twisting of the structure or a warping of the two shell--like struts. Some of the modes associated with higher BLFs still represent very localized deformations (see $\boldsymbol{\varphi}_{14}$); however these are remarkably shifted outside the considered range of BLFs by the multilevel approach.

Optimization histories for the design corresponding to $\bar{\lambda} = 1$ are shown in Figure 14 (a, b). Again, jumps in the compliance and BLFs correspond to increases of the penalization or projection parameters. \autoref{fig:3DcantileverConvCurves} (c) shows the relationship between $\bar{\lambda}$, the volume fraction of the optimized design and the $\delta_{i}$ parameters, defined in \eqref{eq:coalescingMeasure}. As we expected \cite{ferrari-sigmund_19a}, with a maximum prescribed compliance, the volume fraction of the optimized design increases as $\bar{\lambda}$ is raised. Also, more and more buckling modes become simultaneously active: for the design corresponding to $\bar{\lambda} = 1$ we have $\delta_{2}$ and $\delta_{3} < 10^{-5}$ and $\delta_{i} < 5 \cdot 10^{-2}$ up to $\delta_{10}$.

The average computational time for performing the LA (see \ref{Sec:App-NumericalProcedures} for details) is about $185s$. The time spent for the calculation of the 24 buckling modes can be split into $38s$ for Steps 1--2 in \autoref{sSec:multilevelBucklingModes} and about $375s$ for solving \eqref{eq:FineScaleBVP} on the fine discretization (i.e. $\approx 420s$ for the overall EA). Therefore, the ratio $\texttt{tEA}/\texttt{tLBA} = 0.69$ indicates again the efficiency of the multilevel procedure. Comparing to \cite{dunning-etal_16a}, reporting a computational times of about $120s$ when computing 25 buckling modes for a much smaller 3D problem ($\approx 4.6\cdot 10^{5}$ DOFs) and using a parallel algorithm, the presented multilevel approach seems to enhance the efficiency considerably. This also considering that no parallelization (which is possible for all the methods presented) was considered for the present examples.

\section{Concluding discussion}
 \label{Sec:Conclusions}

The goal of this work was to cut the complexity and computational cost of structural Topology Optimization accounting for buckling, to make it feasible for large--scale problems. The results presented in \autoref{Sec:2Dexample-2BarsTruss} and \autoref{Sec:3DcantileverMassMinimization} indicate that a multilevel strategy for selecting buckling modes and to approximate the corresponding buckling load factors makes this goal achievable. The computational effort for obtaining buckling modes and load factors is reduced to a fraction of that required from a full scale eigenvalue analysis; moreover, it scales approximately as the cost of a multi--load linear compliance problem.

We also discussed how the multilevel strategy alleviates some artifacts due to stress concentrations and filters out some local buckling modes. Within this context these are seen as very positive effects, as motivated in \autoref{sSecSpuriousLocalizedModes}. Basically, minimum mass or compliance optimization inherently produces hierarchical layouts with many thin bars prone to undergo local buckling. Taking into account all of these modes would be computationally unfeasible and also unneccessary to the goals of achieving a preliminary design meeting global stability.

The authors believe that a rational and effective approach to tackle local buckling within topology optimization of continua is still an open question, which definitely requires further research efforts. Nevertheless, the procedure we have proposed here represents a first effective method for improving the overall geometric stability of large--scale topology optimized designs.

\section*{Acknowledgement}
This project is supported by the Villum Foundation through the Villum Investigator project ``InnoTop".

\appendix
\begin{algorithm}[tb]
 \caption{Linearized Buckling Analysis by the multilevel iterative method}
  \label{alg:MLeigenvalueAlgorithm}
   \begin{algorithmic}[1]
   		\State \textbf{set} $\ell$ and \textbf{build} $\Omega_{1} \supset \Omega_{2} \supset \ldots \Omega_{j} \supset \ldots \supset \Omega_{\ell}$
   		\Comment{Build nested discretizations}
        \State \textbf{assemble} $K\left[ \mathbf{x} \right]$ on $\Omega_{1}$
        \Comment{Fine scale stiffness matrix}
        \State \textbf{solve} $\mathbf{f} - K\left[ \mathbf{x} \right]\mathbf{u} = \mathbf{r}$ on $\Omega_{1}$
        \Comment{LA performed by \texttt{mgPCG}}
        \State \textbf{assemble} $G\left[ \mathbf{x}, \mathbf{u}\left( \mathbf{x} \right) \right]$ on $\Omega_{1}$
        \Comment{Fine scale stress stiffness matrix}
        \State \textbf{restrict} $K^{\ell} = I^{1}_{\ell}K I^{\ell}_{1}$ and $G^{\ell} = I^{1}_{\ell}G I^{\ell}_{1}$ 
        \Comment{Galerkin projection of matrices on $\Omega_{\ell}$}
   		\State \textbf{compute} $(\lambda^{\ell}_{i}, \boldsymbol{\varphi}^{\ell}_{i})$, $i = 1, \ldots, q$
   		\Comment{Solve coarse scale eigenvalue problem}
      	\State \textbf{set} $\Psi^{\ell} = \{\boldsymbol{\varphi}^{\ell}_{i}\}^{q}_{i=1}$, $\Lambda^{\ell} = {\rm diag}\{ \lambda^{\ell}_{i}\}$
   		\For{$j = \ell - 1, \ldots, 0$}
   		\State $\Psi^{j} \leftarrow I^{j}_{j+1} \Psi^{j+1}$
   		\Comment{Project modes on the next finer grid $\Omega_{j}$}
		\State $\tilde{\lambda}^{j}_{1}
		 \leftarrow \min_{i}\tilde{\Lambda}^{j+1}_{i}$
		\Comment{Set current shift on $\Omega_{j}$}
   		\State $Y^{j} = (K^{j} + \tilde{\lambda}^{j}_{1} G^{j}) \Psi^{j}$
   		\Comment{Compute residuals $Y^{j}=[\mathbf{y}^{j}]^{q}_{i=1}$}
   		\State $\Psi^{j} \leftarrow \mathcal{S}[Y^{j}]$
   		\Comment{Smooth modes $\Psi^{j}$ iterating on the residuals}
   		\State $\tilde{\Lambda}^{j} \leftarrow R[ \Psi^{j} ]$
   		\Comment{Ritz projection to compute approax. $\tilde{\lambda}_{i}$ on $\Omega_{j}$}
   		\EndFor
   		\State \textbf{solve} $K\left[ \mathbf{x} \right] \tilde{\Phi} = G\left[ \mathbf{x}, \mathbf{u}(\mathbf{x}) \right]\Psi$
   		\Comment{Solve one fine scale system by \texttt{mgBPCG}} 
   		\State $\tilde{\Lambda} = R(\tilde{\Phi})$
   		\Comment{Estimate BLFs associated with fine scale modes}
   \end{algorithmic}
\end{algorithm}

%% FIGURE SHOWING THE ACCURACY OBTAINED BY SOLVING THE LINEAR SYSTEM AND THE ADJOINT PROBLEM BY MEANS OF THE MULTILEVEL METHOD
\begin{figure}[tb]
 \centering
  \subfloat[]{
   \includegraphics[scale = 0.45, keepaspectratio]
   {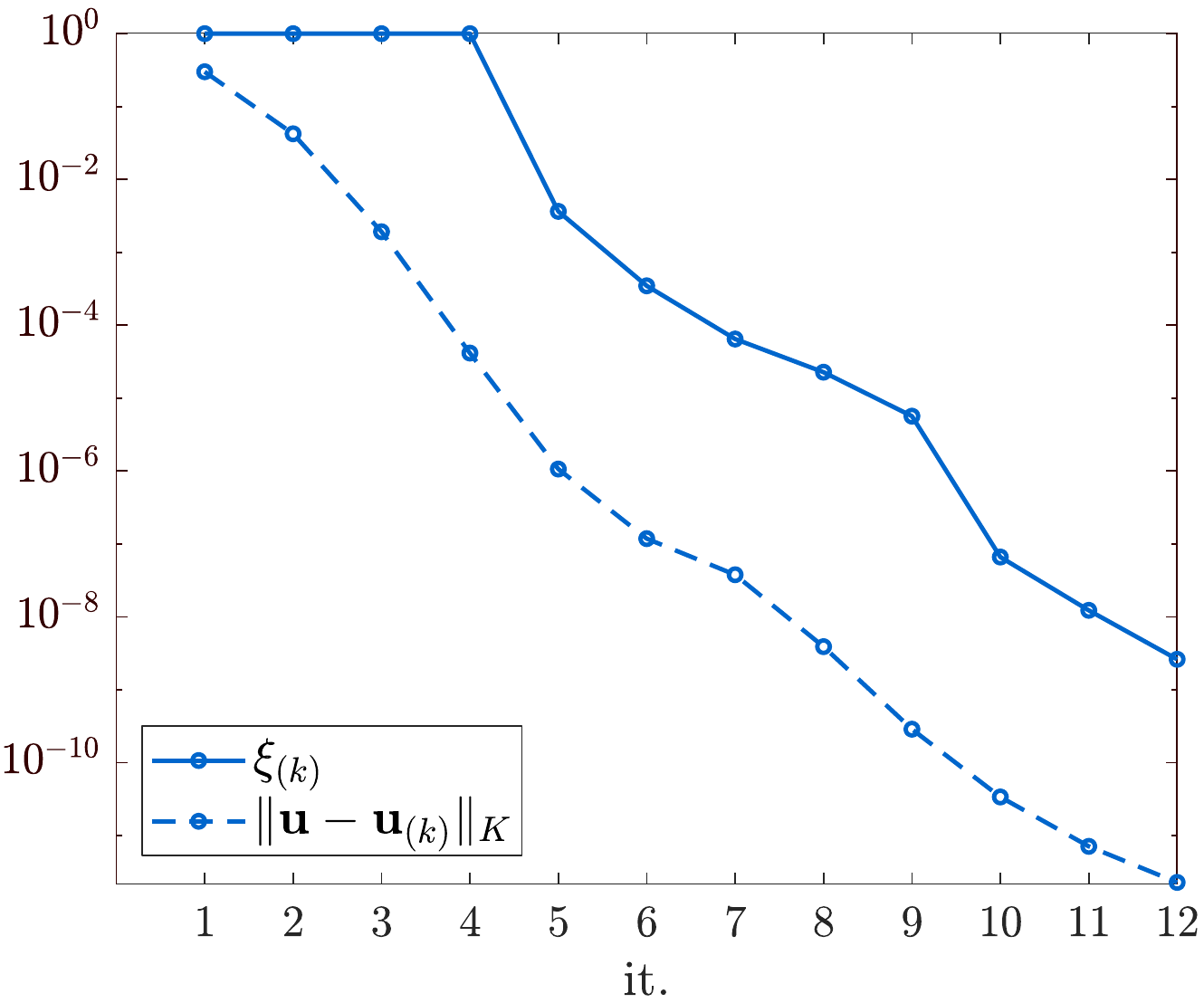}} \qquad
  \subfloat[]{
   \includegraphics[scale = 0.45, keepaspectratio]
   {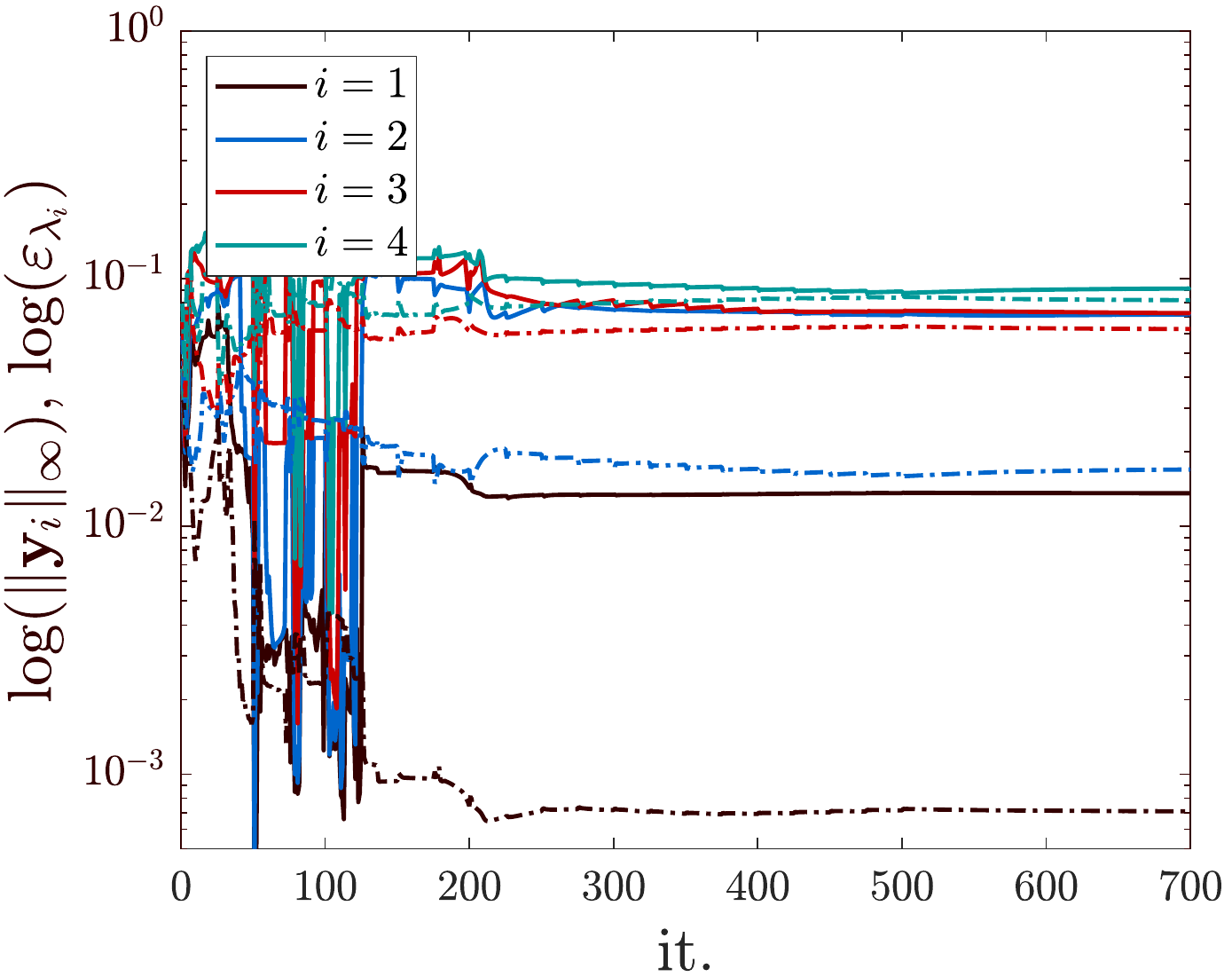}}
 \caption{Some convergence plots corresponding to example \autoref{fig:NumericalTest-2DExample} (c). (a) shows the convergence of the \texttt{mgPCG} iteration, where $\xi_{(k)}$ (solid line) estimates the energy norm error (dashed line), as proposed in \cite{arioli_04a}. (b) shows the measure of the residuals of the eigenvalue equation defined by \autoref{eq:eigResidual} (solid line) and the error in the BLFs (dash--dot line), as the optimization progresses}
 \label{fig:App_accuracyPCGandBPCG}
\end{figure}

\section{Further details on numerical procedures}
 \label{Sec:App-NumericalProcedures}

\autoref{alg:MLeigenvalueAlgorithm} shows the detailed steps for performing the LBA by using iterative solvers and multilevel preconditioners.

The LA (Line 3) is performed by the Conjugate Gradient method, preconditioned by a geometric multigrid built on the set of nested discretizations $\{ \Omega_{j} \}^{\ell}_{j=1}$ \cite{benzi_02a}. The convergence criterion adopted is the one proposed by Arioli \cite{arioli_04a}, which is based on a direct estimate of the energy norm of the error (i.e. $\epsilon_{K} = \|\mathbf{u}_{(k)} - \mathbf{u}\|_{K}$ where $\mathbf{u}_{(k)}$ is the approximation at the $k$--th CG step) making use of quantities which are already computed by the CG iteration (see Eq. 20 in \cite{arioli_04a} and the relative discussion). By using this criterion, convergence has been observed to occurr in 5 to 13 iterations for the 2D examples, and in 6 to 25 for the 3D one (see \autoref{fig:App_accuracyPCGandBPCG}(a)).

The adjoint equation \eqref{eq:AdjointSystem}
is solved by the block version of the \texttt{mgPCG} \cite{oleary_80a, feng-etal_95a} and it is our experience that this may require more iterations. This is reasonable, because the right hand sides $\boldsymbol{{\varphi}}^{T}_{j}(\nabla_{\mathbf{u}}G)\boldsymbol{{\varphi}}_{j}\in\mathbb{R}^{n}$ are now distributed over the whole domain; therefore the error propagates slower.

The multilevel steps for computing an approximation to the fine scale modes (Lines 4--14 in \autoref{alg:MLeigenvalueAlgorithm}) have already been outlined in \autoref{sSec:multilevelBucklingModes}; here we just point out that
\begin{itemize}
 \item The coarse scale eigenvalue problem (Line 6) is solved by the Krylov--Schur algorithm \cite{stewart_02a} and $(\lambda^{\ell}, \boldsymbol{\varphi}^{\ell})_{i}$ are the \emph{only} ``true'' eigenpairs in the overall process; i.e. the only computed running the algorithm to convergence (measured by a tolerance $\tau = 10^{-8}$ on the residual);
 \item On each grid $\Omega_{j}$, the projected modes $\Psi^{j}$ are smoothened by iterating on the associated residuals, in order to filter high frequencies introduced by the projection. Since the matrix $A^{j} := K^{j} + \tilde{\lambda}_{1}G^{j}$ is not positive definite (as $\tilde{\lambda}^{j}_{1} > \lambda^{j}_{1}$), a Kaczmarz iteration \cite{kaczmarz_37a} (denoted as $\mathcal{S}[\cdot]$ in \autoref{alg:MLeigenvalueAlgorithm}) is used, which amounts to the Gauss--Seidel iteration for $(A^{j})^{T}A^{j}$ \cite{wittum_89a};
\end{itemize}

On the fine grid $\Omega_{1}$ the buckling modes approximation can be improved by taking one step of inverse iterations (Line 15--16). Again, this set of linear equations is solved by a preconditioned iteration which, for each vector $\tilde{\boldsymbol{\varphi}}_{i}$, has the form \cite{golub-ye_00a}

\begin{equation}
 \label{eq:PInvIt-GeneralInteration}
  \tilde{\boldsymbol{\varphi}}_{i(k+1)} = 
  \tilde{\boldsymbol{\varphi}}_{i(k)} - P^{-1} (K +
  \tilde{\lambda}_{1} G)\tilde{\boldsymbol{\varphi}}_{i(k)}
\end{equation}
where $P$ is the preconditioner and $\tilde{\boldsymbol{\varphi}}_{i(0)} = \boldsymbol{\psi}_{i}$. \autoref{eq:PInvIt-GeneralInteration} is a gradient method for the minimization of the Rayleigh quotient $R_{i} = R[\tilde{\boldsymbol{\varphi}}_{i}]$ \cite{dyakonov-orekhov_80a, dyakonov_83a}, with descent direction

\begin{equation}
 \label{eq:PInvIt-GradientRayleighQuotient}
 \nabla_{P} R_{i} =
 P^{-1} ( K + \tilde{\lambda}_{1} G)\tilde{\boldsymbol{\varphi}}_{i}
\end{equation}

Referring to \cite{neymeyr_01a} for details and proofs, we recall that, if $R_{i(0)} \in ( \lambda_{i}, \lambda_{i+1} )$, for each step \eqref{eq:PInvIt-GeneralInteration} we have either $R_{i(k+1)} < \lambda_{i}$ (if $i>1)$ or $R_{i(k+1)} \in [R_{i(k)}, \lambda_{i}]$ and the following estimate applies \cite{knyazev-enymeyr_03a}
 \begin{equation}
  \label{eq:errorEstimate}
  \frac{R_{i(k+1)} - \lambda_{i}}
  {\lambda_{i+1} - R_{i(k+1)}} \leq 
  c^{2}
  \frac{R_{i(k)} - \lambda_{i}}
  {\lambda_{i+1} - R_{i(k)}}
 \end{equation}
where $c^{2} \propto \left( 1 - \lambda_{i}/\lambda_{i+1}\right)$.

From \eqref{eq:errorEstimate}, the sequence $\{R_{i(k)}\}_{k}$ monotonically converges to the fine scale eigenvalue $\lambda_{i}$; and by the definition of $c^{2}$ we see that as the BLFs start to coalesce together the ratio $\lambda_{j}/\lambda_{i+1} \rightarrow 1$, so does $c^{2}$ and the convergence rate deteriorates.

Finally, the consistent expression of the sensitivity for the BLFs $\tilde{\lambda}_{i}$, accounting for the residual $\mathbf{y}_{i}$ defined in \eqref{eq:eigResidual}, is

\begin{equation}
 \label{eq:consistentSensitivity}
  \frac{\partial\tilde{\lambda}_{i}}{\partial x_{e}} =
  \frac{1}{1 - \mathbf{p}^{T}_{i}G\tilde{\boldsymbol{\varphi}}_{i}}\left[\tilde{\boldsymbol{\varphi}}^{T}_{i}\left( \frac{\partial K}{\partial x_{e}} + \tilde{\lambda}_{i} \frac{\partial G}{\partial x_{e}} \right)\tilde{\boldsymbol{\varphi}}_{i} - \tilde{\lambda}_{i}\mathbf{z}^{T}_{i}\frac{\partial K}{\partial x_{e}}\mathbf{u} + \mathbf{p}^{T}_{i}\left( \frac{\partial K}{\partial x_{e}} + \tilde{\lambda}_{i} \frac{\partial G}{\partial x_{e}} - \frac{\partial \mathbf{y}_{i}}{\partial x_{e}}\right)\tilde{\boldsymbol{\varphi}}_{i} \right]
\end{equation}
where $\mathbf{p}_{i} = - 2 (K + \tilde{\lambda}_{i}G)^{-1}\mathbf{y}_{i}$ is the adjoint variable associated with \eqref{eq:eigResidual}. \autoref{fig:App_accuracyPCGandBPCG}(b) shows the evolution of $\|\mathbf{y}_{i}\|_{\infty}$ for the lowest four modes, in the optimization progresses for the design \autoref{fig:NumericalTest-2DExample}(c). Even if this measure is not generally small, the corresponding BLFs approximations (represented by dashed lines) are still accurate, especially for $\tilde{\lambda}_{1}$. Therefore, we treat $(\tilde{\lambda}_{i}, \tilde{\boldsymbol{\varphi}}_{i})$ as an approximation of the ``true" eigenpair, and apply the sensitivity expression \eqref{eq:SensitivityLambda}.

% _______________________________ Bibliography

\begin{small}
\section*{References}
\bibliographystyle{spmpsci} % mathematics and physical sciences
\bibliography{main.bib}

\begin{thebibliography}{10}
\providecommand{\url}[1]{{#1}}
\providecommand{\urlprefix}{URL }
\expandafter\ifx\csname urlstyle\endcsname\relax
  \providecommand{\doi}[1]{DOI~\discretionary{}{}{}#1}\else
  \providecommand{\doi}{DOI~\discretionary{}{}{}\begingroup
  \urlstyle{rm}\Url}\fi

\bibitem{aage-etal_17a}
Aage, N., Andreassen, E., Lazarov, B.S., Sigmund, O.: Giga--voxel computational
  morphogenesis for structural design.
\newblock Nature \textbf{550}(7674), 84--86 (2017)

\bibitem{book:achenbach1973}
Achenbach, J.D.: Wave propagation in elastic solids.
\newblock North--Holland (1973).
\newblock Section 13.1

\bibitem{achtziger_99a}
Achtziger, W.: Local stability of trusses in the context of topology
  optimization, {P}art {I}: exact modelling.
\newblock Structural Optimization \textbf{17}(235--246) (1999)

\bibitem{achtziger_99b}
Achtziger, W.: Local stability of trusses in the context of topology
  optimization part ii: {A} numerical approach.
\newblock Structural optimization \textbf{17}(4), 247--258 (1999).
\newblock \doi{10.1007/BF01207000}

\bibitem{allemang-brown_82a}
Allemagn, R.J., Brown, D.L.: A correlation coefficient for modal vector
  analysis.
\newblock In: International Modal Analysis Conference (1982)

\bibitem{allemang_02a}
Allemang, R.J.: The {M}odal {A}ssurance {C}riterion--twenty years of use and
  abuse.
\newblock Journal of Sound and Vibration  (2002)

\bibitem{amir-etal_14a}
Amir, O., Aage, N., Lazarov, B.S.: On multigrid--{CG} for efficient topology
  optimization.
\newblock Structural and Multidisciplinary Optimization \textbf{49}(5),
  815--829 (2014).
\newblock \doi{10.1007/s00158-013-1015-5}

\bibitem{andreassen-etal_17a}
Andreassen, E., Ferrari, F., Sigmund, O., Diaz, A.: Frequency response as a
  surrogate eigenvalue problem in topology otpimization.
\newblock Internal Journal for Numerical Methods in Engineering
  \textbf{113}(8), 1214--1229 (2017)

\bibitem{arioli_04a}
Arioli, M.: A stopping criterion for the conjugate gradient algorithm in a
  finite element method framework.
\newblock Numerische Mathematik \textbf{97}(1), 1--24 (2004).
\newblock \doi{10.1007/s00211-003-0500-y}

\bibitem{book:bazant-cedolin2010}
Bazant, Z.D., Cedolin, L.: Stability of Structures.
\newblock World Scientific (2010)

\bibitem{bendsoe_89a}
Bends\o{}e, M.P.: Optimal shape design as a material distribution problem.
\newblock Structural optimization \textbf{1}(4), 193--202 (1989).
\newblock \doi{10.1007/BF01650949}

\bibitem{book:bendsoe-sigmund_2004}
Bends\o{}e, M.P., Sigmund, O.: Topology Optimization: Theory, Methods and
  Applications.
\newblock Springer (2004)

\bibitem{benzi_02a}
Benzi, M.: Preconditioning techniques for large linear systems: A survey.
\newblock Journal of Computational Physics \textbf{182}(2), 418 -- 477 (2002).
\newblock \doi{http://dx.doi.org/10.1006/jcph.2002.7176}

\bibitem{bian-feng_17a}
Bian, X., Feng, Y.: Large--scale buckling--constrained topology optimization
  based on assembly--free finite element analysis.
\newblock Advances in Mechanical Engineering \textbf{9}(9), 1--12 (2017)

\bibitem{book:deborst2012}
de~Borst, R., Crisfield, M.A., Remmers, J.J.C., Verhoosel, C.V.: Non--Linear
  Finite Element Analysis of Solids and Structures, second edition edn.
\newblock John Wiley \& Sons (2012)

\bibitem{bourdin_01a}
Bourdin, B.: Filters in topology optimization.
\newblock International Journal for Numerical Methods in Engineering
  \textbf{50}(9), 2143--2158 (2001).
\newblock \doi{10.1002/nme.116}

\bibitem{brehm-etal_10a}
Brehm, M., Zabel, V., Bucher, C.: An automatic mode pairing strategy using an
  enhanced modal assurance criterion based on modal strain energies.
\newblock Journal of Sound and Vibration \textbf{329}(25), 5375--5392 (2010).
\newblock \doi{https://doi.org/10.1016/j.jsv.2010.07.006}

\bibitem{book:briggs00}
Briggs, W., Henson, V., McCormick, S.: A Multigrid Tutorial: Second Edition.
\newblock Other Titles in Applied Mathematics. Society for Industrial and
  Applied Mathematics (2000)

\bibitem{bruyneel-etal_08a}
Bruyneel, M., Colson, B., Remouchamps, A.: Discussion on some convergence
  problems in buckling optimisation.
\newblock Structural and Multidisciplinary Optimization \textbf{35}(2),
  181--186 (2008)

\bibitem{chin-kennedy_16a}
Chin, T.W., Kennedy, G.J.: Large--scale compliance--minimization and buckling
  topology optimization of the undeformed common research model wing.
\newblock In: AIAA SciTechForum (2016)

\bibitem{clausen-etal_16a}
Clausen, A., Aage, N., Sigmund, O.: Exploiting additive manufacturing infill in
  topology optimization for improved buckling load.
\newblock Engineering \textbf{2}(2), 250--257 (2016).
\newblock \doi{10.1016/J.ENG.2016.02.006}

\bibitem{clausen-andreassen_17a}
Clausen, A., Andreassen, E.: On filter boundary conditions in topology
  optimization.
\newblock Structural and Multidisciplinary Optimization \textbf{56}(5),
  1147--1155 (2017)

\bibitem{dunning-etal_16a}
Dunning, P.D., Ovtchinnikov, E., Scott, J., Kim, A.: Level--set topology
  optimization with many linear buckling constraints using and efficient and
  robust eigensolver.
\newblock International Journal for Numerical Methods in Engineering  (2016)

\bibitem{dyakonov_83a}
D'yakonov, E.G.: Iteration methods in eigenvalue problems.
\newblock Mathematical notes of the Academy of Sciences of the USSR
  \textbf{34}(6), 945--953 (1983).
\newblock \doi{10.1007/BF01157412}

\bibitem{dyakonov-orekhov_80a}
D'yakonov, E.G., Orekhov, M.Y.: Minimization of the computational labor in
  determining the first eigenvalues of differential operators.
\newblock Mathematical notes of the Academy of Sciences of the USSR
  \textbf{27}(5), 382--391 (1980).
\newblock \doi{10.1007/BF01139851}

\bibitem{feng-etal_95a}
Feng, Y., Owen, D., Peri\'{c}, D.: A block conjugate gradient method applied to
  linear systems with multiple right-hand sides.
\newblock Computer Methods in Applied Mechanics and Engineering
  \textbf{127}(1), 203 -- 215 (1995).
\newblock \doi{http://dx.doi.org/10.1016/0045-7825(95)00832-2}

\bibitem{ferrari-etal_18a}
Ferrari, F., Lazarov, B.S., Sigmund, O.: Eigenvalue topology optimization via
  efficient multilevel solution of the {F}requency {R}esponse.
\newblock International Journal for Numerical Methods in Engineering
  \textbf{115}(7), 872--892 (2018)

\bibitem{ferrari-sigmund_19a}
Ferrari, F., Sigmund, O.: Revisiting topology optimization with buckling
  constraints.
\newblock Structural and Multidisciplinary Optimization \textbf{59}(5),
  1401--1415 (2019).
\newblock \doi{10.1007/s00158-019-02253-3}

\bibitem{gao-ma_15a}
Gao, X., Ma, H.: Topology optimization of continuum structures under buckling
  constraints.
\newblock Computers \& Structures \textbf{157}, 142--152 (2015)

\bibitem{golub-ye_00a}
Golub, G.H., Ye, Q.: Inexact inverse iteration for generalized eigenvalue
  problems.
\newblock BIT Numerical Mathematics \textbf{40}(4), 671--684 (2000).
\newblock \doi{10.1023/A:1022388317839}

\bibitem{jordan_881a}
Jordan, C.: Sur la s\'{e}rie de fourier.
\newblock Comptes Rendus Hebdomadaires des Seances de l'Academie des Sciences
  \textbf{92}, 228--230 (1881)

\bibitem{kaczmarz_37a}
Kaczmarz, S.: Angen\"{a}herte aufl\"{o}sung von systemen linearer gleichungen.
\newblock Bulletin de l'Academie Polonaises des Sciences et Lettres
  \textbf{A35}, 355--357 (1937)

\bibitem{khot-etal_76a}
Khot, N.S., Venkayya, V.B., Berke, L.: Optimum structural design with stability
  constraints.
\newblock International Journal for Numerical Methods in Engineering
  \textbf{10}(5), 1097--1114 (1976)

\bibitem{knyazev_01a}
Knyazev, A.V.: Toward the optimal preconditioned eigensolver: Locally optimal
  block preconditioned conjugate gradient method.
\newblock SIAM Journal on Scientific Computing \textbf{23}(2), 517--541 (2001).
\newblock \doi{10.1137/S1064827500366124}

\bibitem{knyazev-enymeyr_03a}
Knyazev, A.V., Neymeyr, K.: A geometric theory for preconditioned inverse
  iteration {III}: {A} short and sharp convergence estimate for generalized
  eigenvalue problems.
\newblock Linear Algebra and its Applications \textbf{358}(1), 95 -- 114
  (2003).
\newblock \doi{http://dx.doi.org/10.1016/S0024-3795(01)00461-X}

\bibitem{kreisselmeier-steinhauser_79a}
Kreisselmeier, G., Steinhauser, R.: Systematic control design by optimizing a
  vector performance index.
\newblock IFAC Proceedings Volumes \textbf{12}(7), 113 -- 117 (1979).
\newblock IFAC Symposium on computer Aided Design of Control Systems, Zurich,
  Switzerland, 29-31 August

\bibitem{lazarov-etal_16a}
Lazarov, B.S., Wang, F., Sigmund, O.: Length scale and manufacturability in
  density-based topology optimization.
\newblock Archive of Applied Mechanics \textbf{86}(1), 189--218 (2016).
\newblock \doi{10.1007/s00419-015-1106-4}

\bibitem{ma-etal_93a}
Ma, Z.D., Kikuchi, N., Hagiwara, I.: Structural topology and shape optimization
  for a frequency response problem.
\newblock Computational Mechanics \textbf{13}(3), 157--174 (1993).
\newblock \doi{10.1007/BF00370133}

\bibitem{neves-etal_95a}
Neves, M.M., Rodrigues, H., Guedes, J.M.: Generalized topology design of
  structures with a buckling load criterion.
\newblock Structural optimization \textbf{10}(2), 71--78 (1995)

\bibitem{neymeyr_01a}
Neymeyr, K.: A geometric theory for preconditioned inverse iteration {I}:
  {E}xtrema of the {R}ayleigh quotient.
\newblock Linear Algebra and its Applications \textbf{322}(1), 61 -- 85 (2001).
\newblock \doi{http://dx.doi.org/10.1016/S0024-3795(00)00239-1}

\bibitem{oleary_80a}
O'Leary, D.P.: The block conjugate gradient algorithm and related methods.
\newblock Linear Algebra and its Applications \textbf{29}, 293 -- 322 (1980).
\newblock \doi{http://dx.doi.org/10.1016/0024-3795(80)90247-5}

\bibitem{ovtchinnikov_08a}
Ovtchinnikov, E.E.: Computing several eigenpairs of hermitian problems by
  conjugate gradient iterations.
\newblock Journal of Computational Physics \textbf{227}(22), 9477--9497 (2008)

\bibitem{pedersen_00a}
Pedersen, N.: Maximization of eigenvalues using topology optimization.
\newblock Structural and Multidisciplinary Optimization \textbf{20}(1), 2--11
  (2000).
\newblock \doi{10.1007/s001580050130}

\bibitem{pian-cheng_82a}
Pian, T.H.H., Chen, D.P.: Alternative ways for formulation of hybrid stress
  elements.
\newblock International Journal for Numerical Methods in Engineering
  \textbf{18}(11), 1679--1684 (1982)

\bibitem{pian-sumihara_84a}
Pian, T.H.H., Sumihara, K.: Rational approach for assumed stress finite
  elements.
\newblock International Journal for Numerical Methods in Engineering
  \textbf{20}(9), 1685--1695 (1984)

\bibitem{book:pian-wu2005}
Pian, T.H.H., Wu, C.C.: Hybrid and incompatible finite element methods.
\newblock Taylor \& Francis (2005)

\bibitem{rodrigues-etal_95a}
Rodrigues, H.C., Guedes, J.M., Bends\o{}e, M.P.: Necessary conditions for
  optimal design of structures with a nonsmooth eigenvalue based criterion.
\newblock Structural Optimization \textbf{9}, 52--56 (1995)

\bibitem{rozvany_96a}
Rozvany, G.: Difficulties in topology optimization with stress, local buckling
  and system stability constraints.
\newblock Structural Optimization \textbf{11}, 213--217 (1996)

\bibitem{sigmund_07a}
Sigmund, O.: Morphology--based black and white filters for topology
  optimization.
\newblock Structural and Multidisciplinary Optimization \textbf{33}(4),
  401--424 (2007)

\bibitem{stewart_02a}
Stewart, G.: A {K}rylov--{S}chur algorithm for large eigenproblems.
\newblock SIAM Journal on Matrix Analysis and Applications \textbf{23}(3),
  601--614 (2002)

\bibitem{svanberg_87a}
Svanberg, K.: The method of moving asymptotes - {A} new method for structural
  optimization.
\newblock International Journal for Numerical Methods in Engineering
  \textbf{24}(2), 359--373 (1987)

\bibitem{thompson-hunt_74a}
Thompson, J.M.T., Hunt, C.W.: Dangers of structural optimization.
\newblock Engineering Optimization \textbf{1}(2), 99--110 (1974).
\newblock \doi{10.1080/03052157408960580}

\bibitem{thomsen-etal_18a}
Thomsen, C.R., Wang, F., Sigmund, O.: Buckling strength topology optimization
  of 2{D} periodic materials based on linearized bifurcation analysis.
\newblock Computer Methods in Applied Mechanics and Engineering \textbf{339},
  115--136 (2018)

\bibitem{wang-etal_11a}
Wang, F., Lazarov, B., Sigmund, O.: On projection methods, convergence and
  robust formulations in topology optimization.
\newblock Structural and Multidisciplinary Optimization \textbf{43}(6),
  767--784 (2011)

\bibitem{book:washizu}
Washizu, K.: Variational Methods in Elasticity and Plasticity, {S}econd edn.
\newblock Pergamon Press (1975)

\bibitem{wilson-ibrahimbergovic_90a}
Wilson, E.L., Ibrahimbergovic, A.: Use of incompatible displacement modes for
  the calculation of element stiffness or stresses.
\newblock Finite Elements in Analysis and Design \textbf{7}, 229--241 (1990)

\bibitem{wilson-etal_73a}
Wilson, E.L., Taylor, R.L., Doherty, W., Glaboussi, J.: Incompatible
  displacement models, pp. 41--57.
\newblock Academic Press (1973)

\bibitem{wittum_89a}
Wittum, G.: Multi-grid methods for stokes and {N}avier--{S}tokes equations.
\newblock Numerische Mathematik \textbf{54}(5), 543--563 (1989)

\bibitem{zhou-saad_08a}
Zhou, Y., Saad, Y.: Block {K}rilov--{S}chur mehtod for large symmetric
  eigenvalue problems.
\newblock Numerical Algorithms \textbf{47}(4), 341--359 (2008)

\end{thebibliography}
\end{small}
\end{document}